\documentclass[11pt,leqno]{article}
\usepackage{amsthm,amsfonts,amssymb,epsfig,graphics,amsmath,oldgerm}
 \usepackage[latin1]{inputenc}\relax

\newcommand{\range}{\text{\rm range}}

\newcommand{\trans}{\text{\rm tr}}
\newcommand{\diag}{\text{\rm diag}}

\newcommand{\RR}{{\mathbb R}}
\newcommand{\NN}{{\mathbb N}}

\newcommand{\CC}{{\mathbb C}}
\newcommand{\EE }{{\mathbb E}}

\newcommand{\FF}{{\mathbb F}}

\newcommand\cA{{\cal  A}}
\newcommand\cB{{\cal  B}}

\newcommand\cU{{\cal  U}}

\newcommand\cV{{\cal  V}}
\newcommand\cW{{\cal  W}}
\newcommand\cC{{\cal  C}}

\newcommand\cG{{\cal  G}}

\newcommand\cL{{\cal  L}}

\newcommand\cE{{\cal  E}}

\newcommand\cP{{\cal  P}}

\newcommand\cX{{\cal X}}

\newcommand\cM{{\mathcal M}}

\newcommand\cS{{\mathcal S}}

\newcommand{\rmh}{{\mathrm{h } }}

\newcommand{\rmp}{{\mathrm{p } }}

\newcommand{\rC}{{\mathrm{C} }}

\newcommand{\bdiag}{{\text{\rm block-diag} }}

\newcommand{\cs}{\check \sigma}
\newcommand{\cz}{\check \zeta}
\newcommand{\cg}{\check \gamma}
\newcommand{\ch}{\check \eta}
\newcommand{\ct}{\check \tau}

\newcommand{\chH}{\check H}
\newcommand{\cx}{\check \xi}

\newcommand{\ccD}{\check D}

\newcommand{\ccSi}{\check \Sigma}

\newcommand{\up}{\underline p} 
 
\newcommand{\utau}{\underline \tau}
\newcommand{\uxi}{\underline \xi}

\newcommand{\uzeta}{\underline \zeta}

\newcommand{\oA}{\overline{A}}
\newcommand{\oB}{\overline{B}}
\newcommand{\oL}{\overline{L}}

\newcommand{\ls}{{\ \lesssim \ }}

\newcommand{\eps}{\varepsilon }
\newcommand{\vp}{\varphi}
\newcommand{\D}{\partial }
\newcommand\adots{\mathinner{\mkern2mu\raise1pt\hbox{.}
\mkern3mu\raise4pt\hbox{.}\mkern1mu\raise7pt\hbox{.}}}

\newcommand{\Id}{{\rm Id }}
\newcommand{\im}{{\rm Im }\, }
\newcommand{\re}{{\rm Re }\, }

\renewcommand{\div}{{\rm div}}

\newcommand{\curl}{{\rm curl}}
\newcommand{\na}{{\nabla}}

\newcommand{\la}{\langle }
\newcommand{\ra}{\rangle }

\newcommand{\PG}{{\mathrm{P}\Gamma }} 

\newcommand{\mez}{{\frac{1}{2}}}

\newtheorem{theo}{Theorem}[section]
\newtheorem{prop}[theo]{Proposition}
\newtheorem{cor}[theo]{Corollary}
\newtheorem{lem}[theo]{Lemma}
\newtheorem{defi}[theo]{Definition}
\newtheorem{ass}[theo]{Assumption}

\newtheorem{exam}[theo]{Example}
\newtheorem{rem}[theo]{Remark}

\newtheorem{cexam}[theo]{Counter example}
\newtheorem{nota}[theo]{Notations}

\numberwithin{equation}{section}

\title{Viscous Boundary Value Problems for 
Symmetric Systems with Variable Multiplicities}

 \author{Olivier Gues, Guy M\'etivier, Mark Williams, Kevin Zumbrun}

\begin{document}

\maketitle

\begin{abstract}
Extending investigations of M\'etivier\&Zumbrun in the hyperbolic
case, we treat stability of viscous shock and boundary layers
for viscous perturbations of multidimensional hyperbolic systems
with characteristics of variable multiplicity,
specifically the construction of symmetrizers in the low-frequency
regime where variable multiplicity plays a role.
At the same time, we extend the boundary-layer theory to ``real'' or partially
parabolic viscosities, Neumann or mixed-type parabolic boundary conditions,
and systems with nonconservative form, in addition proving a more fundamental
version of the Zumbrun--Serre--Rousset theorem, 
valid for variable multiplicities, 
characterizing the limiting hyperbolic system
and boundary conditions as a nonsingular limit of a reduced viscous system.
The new effects of viscosity are seen to be surprisingly subtle; in particular,
viscous coupling of crossing hyperbolic modes may 
induce a destabilizing effect.
We illustrate the theory with applications to magnetohydrodynamics.
  \end{abstract}

\tableofcontents

\bigbreak

\section{Introduction}

This work  is motivated by the stability analysis of boundary value problems and shock waves
for viscous perturbations of multidimensional systems of conservation laws. 
In this analysis, three main steps are present. 
 In the first step, one constructs simple waves
$w((\nu\cdot x - \sigma t)/ \eps)$, 
or ``profiles'',
which are exact solutions of the viscous equation, with viscosity of order
$\eps$. This amounts to solving
 an ordinary differential equation (the profile equation);
the solutions describe the fast transition  between  the hyperbolic solution and 
the parabolic  boundary conditions (boundary layers) or between two smooth 
hyperbolic solutions (shock layers). 
Next, given a profile,   
formal 
plane wave or spectral   analysis yields 
necessary  stability conditions in terms of 
  Evans functions. 
The second 
step
is to compute explicitly this function on specific examples and 
  check the stability conditions.   
The third 
step
is to prove   
the   linear and nonlinear stability of solutions, 
assuming that the suitable Evans-Lopatinski  condition is satisfied, in particular for curved fronts or boundaries
and non-piecewise  constant hyperbolic solutions. 
This paper deals with the third 
step,
with specific applications to magneto-hydrodynamics.  
The first and second steps are discussed for shock and
boundary layers, respectively, in companion papers \cite{GMWZ7}
and \cite{GMWZ5}.

We concentrate on the construction of symmetrizers for the linearized equations, and more specifically in the so called low-frequency regime,
 as they are the key point in the proof of stability estimates which 
eventually yield short-time existence and nonlinear stability theorems; 
see \cite{Ma1} \cite{Ma2} for hyperbolic shocks and 
\cite{MZ1, GMWZ2,GMWZ3,GMWZ4} for viscous perturbations. 
In these papers, 
it is proved that strong stability estimates hold,  under the natural  uniform Lopatinski condition, 
or Evans' condition, 
provided  that the equations satisfy  a  structural condition, called  
\emph{the block structure condition} (see  \cite{Ma1, MaOs} in the hyperbolic 
case and \cite{MZ1} for the viscous case).  This condition is in some sense necessary 
 for the construction of Kreiss' symmetrizers which are used 
to   prove  the stability estimates.
It is satisfied  in the 
case of  inviscid Euler's equations of gas dynamics (\cite{Ma1}),  but 
does not hold in other interesting examples  such as  the equations of magneto-hydrodynamics (MHD).  
 So there is a real need for an extension of the analysis beyond the class 
 of systems satisfying the block structure condition. 
 This is done in  \cite{MZ2} for hyperbolic systems and the main goal of this 
 paper is to extend the analysis to viscous systems,  in view of  applications 
 to MHD.

We carry out in passing several other useful generalizations
of the basic boundary-layer analyisis of \cite{MZ1, Met3},
extending the theory 
to ``real'' or partially parabolic viscosities, 
Neumann or mixed-type parabolic boundary conditions,
and systems with nonconservative form.
In addition, we prove a more fundamental
version of the Zumbrun--Serre--Rousset theorem, 
valid for variable multiplicities, 
characterizing the limiting hyperbolic system
and boundary conditions as a nonsingular limit of a reduced viscous system
as frequency goes to zero.
Extensions to the shock case are given in \cite{GMWZ7}.

 Consider boundary value problems for hyperbolic systems 
 \begin{equation}
\label{hypeq}
\D_t u + \sum_{j=1}^d  A_j \D_j u 
\end{equation}
 on $\{ x_d \ge 0\}$ with boundary conditions on $\{x_d = 0\}$ which  is assumed 
 to be noncharacteristic. The plane wave analysis of such system leads to consider 
 ordinary differential system in $ z \ge 0$, depending on Fourier-Laplace 
 frequencies $\zeta = (\tau, \eta, \gamma)$,  with $ \gamma > 0$, 
 \begin{equation}
\label{sysh}
\D_z - H_0 (\zeta) , \quad  H_0 (\zeta) := - A_d^{-1} \Big( (i \tau + \gamma) 
+ \sum_{j=1}^{d-1} \eta_j A_j \Big) 
\end{equation}
 Viscous perturbations of \eqref{hypeq} are systems of the form 
 \begin{equation}
\label{viseq}
\D_t u + \sum_{j=1}^d A_j \D_j u  - \eps \sum_{j, k = 1}^d \D_j \big ( B_{j, k} \D_k u \big)
\end{equation}
with natural  structural conditions  which are recalled 
below. 
The \emph{low-frequency plane wave analysis} of such systems lead to consider 
perturbations of \eqref{sysh}: 
   \begin{equation}
\label{ode}
\D_z - H (\zeta, \rho) ) , \quad  H  (\zeta, \rho ) :=  H_0(\zeta) + \rho H_1 (\zeta, \rho)  
\end{equation}
depending smoothly on an additional parameter $\rho \ge 0$
(see section 2 below or \cite{MZ1, Z1, Met3}). 

 In this paper, our main concern is the construction of symmetrizers
  $\Sigma(\zeta, \rho)$ for \eqref{ode}. The precise conditions 
we impose   on $\Sigma$ are given in   Section \ref{S3}. 
In particular, we focus on \emph{smooth symmetrizers}, as they serve as 
symbols for pseudodifferential symmetrizers in the variable coefficient analysis. 
 
 When $\rho = 0$, $\Sigma_0(\zeta) = \Sigma (\zeta, 0)$ is a symmetrizer for $H_0(\zeta)$. 
 Such symmetrizers were  constructed first for strictly hyperbolic systems \eqref{hypeq}
 by Kreiss (\cite{Kr}) (see also \cite{CP}). Strict hyperbolicity is used at only one place: 
 it implies  that   $H_0$  can be put in a normal form, 
which is called  block structure in \cite{Ma1,MaOs}.  Therefore Kreiss' construction of symmetrizers extends immediately to systems which satisfy 
this \emph{block structure condition}. 
In \cite{MZ2}, it is proved that  this  condition  is satisfied
if and only if the symbol $A(\xi) := \sum \xi_j A_j$ of \eqref{hypeq} is smoothly diagonalizable 
for $\xi \ne 0$, recovering known examples  
such as Euler's equation of gas dynamics or Maxwell's equations.  
The second important result in \cite{MZ2} 
 it that  the construction of symmetrizers  is extended to a   class of 
symmetric systems which are not smoothly diagonalizable for 
all $\xi \ne0$: we demand that the  ``bad''   multiple modes are 
\emph{totally incoming  or totally outgoing} (see the Definitions~\ref{def22} and 
\ref{def23} below),  and this applies 
  to inviscid MHD. 
 
 In the small viscosity case, the construction of symmetrizers is performed
 in \cite{MZ1}, with application to the analysis of shocks in \cite{GMWZ1, GMWZ2, 
 GMWZ3, GMWZ4}, assuming that the eigenvalues of $A(\xi)$ 
have constant multiplicity  
for $\xi \ne 0$.  As mentioned above, this assumption 
rules out the case of MHD. 
 
The main objective of this paper is to start the analysis of 
\eqref{viseq} or \eqref{ode} when the constant multiplicity assumption is relaxed
and in particular to investigate the construction of symmetrizers. 
   It turns out that the influence of the viscosity 
is much more  subtle  than expected near multiple modes.  In some cases, it may induce destabilizing effects.  Let us list several new phenomena which can occur when  
 there are multiple modes with nonconstant multiplicity.

$\bullet$ Smooth diagonalization of $A_0(\xi)$ implies a smooth block reduction 
for $H_0(\zeta)$. The perturbation $\rho H_1$ in general couples the different blocks 
associated to a multiple eigenvalues (and this occurs for MHD). 
If the  crossing  eigenvalues do not have the same behavior with respect to the boundary
(typically if they are not all incoming  nor all outgoing), the spectral negative space 
$\EE^- (\zeta, \rho)$ \emph{is not continuous} (in general) at $\rho = 0$. 
This happens for slow shock waves in MHD. 
This phenomena is excluded  when the eigenvalues have constant 
multiplicities; see \cite{MZ3}
(in this case, since crossing eigenvalues are equal, they have  
 the same behavior with respect to the boundary). 
 Recall from \cite{MZ2} 
that the continuity of $\EE^-$ is a \emph{necessary condition} 
 for Kreiss'  construction of \emph{smooth} symmetrizers, more precisely  
 for the existence of what is called below, smooth K-families of symmetrizers. 
 In any case, \emph{the discontinuity of $\EE^-$ is a major  difficulty in the 
 construction of symmetrizers}.

$\bullet$  As a consequence of the previous phenomenon, the 
\emph{Evans function can be discontinuous} at $\rho = 0$. 
In the shock problem, the usual Evans function is in every case singular 
at $\rho = 0$ (see \cite{ZS}), but the remark applies to the modified 
(or desingularized) Evans function introduced in \cite{GMWZ3, GMWZ4} (see below).

$\bullet$ Because of the lack of continuity of the Evans function, it may happen that
\emph{the strong Lopatinsky  stability condition for the hyperbolic problem (at $\rho = 0)$ 
 does not imply   the strong Evans stability condition for small $\rho$}.  
This is in sharp contrast with the known results obtained in the constant 
multiplicty case (\cite{ZS, Z1, Rou,    MZ1, GMWZ3, GMWZ4}). 
This is illustrated by an example in Section \ref{S4} and this can occur for MHD, 
for some ad hoc boundary condition. 
An interesting question is to know whether this can happen or not for 
physical boundary conditions, in particular for slow 
MHD 
shocks.

\medbreak

On the other hand, we prove in this paper  the existence of 
smooth symmetrizers under a natural 
\emph{generalized block structure condition} for 
\eqref{ode}.  We also provide a geometrical characterization of this condition
on the matrices $A$ and $B$ occurring in \eqref{viseq}.  
Moreover, modes  that are totally incoming or totally outgoing do not cause trouble 
in the analysis of $\EE^-$ nor of  the Evans function. They are easily handled 
in the case of symmetric systems
as in \cite{MZ2}. 
 For instance, an important outcome of the present paper is the following result. 
We refer to the next sections for precise definitions. 

\begin{theo} 
\label{th11}
Suppose that the full system \eqref{viseq} is symmetric. Suppose in addition that  the eigenvalues of the hyperbolic 
system \eqref{hypeq} are either semi-simple with constant multiplicity or totally nonglancing
in the sense of Definition~$\ref{def23}$.   
Then, there are K-families of  symmetrizers for the  
associated reduced system \eqref{ode}, 
for $\rho\ge 0$ sufficiently small. 

\end{theo}

 As recalled in the next section, K-families of  symmetrizers provide Kreiss symmetrizers 
 for boundary value problems which satisfy a uniform Lopatinski stability condition. 
 One important application and motivation is the following 

\begin{exam}
\label{ex1}
Fast  Lax' shocks for MHD satisfy the  assumptions of Theorem~$\ref{th11}$. 
\end{exam}

But we also have the following 

\begin{cexam}
\label{ex2}
Slow  Lax' shocks for MHD do not satisfy the  assumptions of 
Theorem~$\ref{th11}$. 

\end{cexam}

\begin{rem}
\textup{Fast shocks  with small magnetic field are perturbations 
of acoustic shocks of gas dynamics, whose stability has been studied by A.Majda (\cite{Ma1}). 
Therefore, there are good reasons to think that the uniform Evans-Lopatinski 
condition is satisfied for Fast  Lax' shocks for MHD, at least for 
perfect gases state laws and small magnetic field. 
 }
\end{rem}

\begin{rem}
\textup{When the assumptions of Theorem~\ref{th11} are not satisfied, or more generally
when the generalized block structure fails, one could try to construct 
nonsmooth symmetrizers as in \cite{MZ2}. The counterexample~\ref{ex2} 
would be a good motivation for that. However, nonsmooth symmetrizer
would require  much more sophisticated pseudodifferential  tools to handle 
  variable coefficients. 
Moreover, slow shocks are not so closely related to acoustic shocks, 
and it is not known whether the uniform Lopatinski condition is likely to be satisfied or not.  }
\end{rem}

When all the eigenvalues have constant multiplicity, Theorem~\ref{th11} is proved in 
\cite{MZ1} (see also \cite{Met3}).\footnote{
The reduction to \eqref{ode} is carried out for strictly parabolic
viscosities in \cite{MZ1, Met3} and for partial viscosities in \cite{GMWZ4}.
However, the form of $H_1$ is the same in each case (a 
consequence of Kawashima's genuine coupling condition, 
Assumption (H4) below).
}
The construction is based on a reduction of \eqref{ode} 
to a suitable block diagonal form.  Blocks which correspond to totally nonglancing modes
(incoming or outgoing) are treated using the symmetry of the system as in \cite{MZ2}. 
For other blocks, 
 we discuss in detail in Section 4 the   
\emph{generalized block structure condition} which is needed for the construction of 
Kreiss symmetrizers.

The symmetrizers are used in \cite{MZ1, GMWZ1, GMWZ2, GMWZ3, GMWZ4} 
to prove maximal  stability estimates for boundary value problems.  
  The  Fourier multipliers  $\Sigma(p, \zeta)$ serve as 
  symbols  for  pseudo-differential symmetrizers.  
 All the other steps in these papers, linearization, paralinearization, 
separation of frequencies, the high- and medium-frequency analysis, 
the conversion of the plane wave or symbolic calculus into an operator 
calculus via the use of a paradifferential calculus, are independent  of the 
constant multiplicity assumption  which was assumed there as a sufficient condition 
for the generalized block structure condition. 
Therefore, all these analyses are valid under the assumptions of Theorem~\ref{Kfam}.

As already mentioned, the main novelty of this paper with respect to previous works
of the authors  is the consideration of systems with variable multiplicity. 
To lighten the presentation, we will now 
 we focus on \emph{boundary layers} for noncharacteristic boundary value problems.  
The extension to classical, conservative Lax-type shocks requires 
only to incorporate the ideas already 
explained in detail in (for instance) \cite{GMWZ3, GMWZ4}. 
(The treatment of nonconservative 
and or undercompressive shocks involves new issues,
and is carried out in \cite{GMWZ7}.)
Similarly, we will concentrate only on the symbolic analysis for 
\emph{constant-coeffient
equations} and the construction of  smooth Fourier-Laplace multipliers. 
The passage from these multipliers to linear and nonlinear stability estimates
for variable coefficients is already performed in previous works (see \cite{MZ1, GMWZ3, GMWZ4})
and  can be used as an independent black box.


\section{Spectral stability }

In this section, we recall the main steps of the spectral stability analysis  of noncharacteristic  boundary layers, refereeing to 
  \cite{Gue, GG,  MZ1, Met3, Z1, Zhandbook,   GMWZ3, GMWZ4} for 
 details and  further references and applications to the similar analysis of  shock profiles.  In particular, we give a new proof 
 of the Zumbrun-Serre Lemma (\cite{ZS, Rou}) which allows for variable multiplicities. 
 Moreover,  not only does it provide a comparison between the Evans function of the viscous
 equation and the Lopatinski determinant of the inviscid system, but it also shows 
 the link between the equations themselves: for low frequencies, the viscous boundary value problem  decouples into two boundary value problems, one of them being a \emph{nonsingular
 perturbation}  of the limiting hyperbolic boundary value problem. We will also recall from 
 \cite{GMWZ4} the main arguments for the high-frequency regime.


 \subsection{Structural assumptions}
Consider a system of equations
\begin{equation}
\label{visceq}
\cL_ \eps  (u) := 
A_0(u) \D_t u + \sum_{j=1}^d   A_j(u)\D_j u   -
 \eps \sum_{j,k= 1}^d \D_j \big( B_{j,k}(u) \D_k u \big) = 0. 
\end{equation}
When $\eps = 0$, $\cL_0$ is first order and assumed to hyperbolic; $\eps$ plays the role 
of a non-dimensional  viscosity and for 
$\eps > 0$, the system is assumed to be parabolic or at least
partially parabolic.  Classical examples are the Navier-Stokes  
equations of gas dynamics, or the equations of magneto-hydrodynamics (MHD). 

The form of the equations is preserved 
under a  change of unknowns 
$  u = \Phi (\tilde u)$  or  multiplication on the left by 
  a \emph{constant}  invertible matrix. 
To cover the case of partial viscosity and 
motivated by the examples of Navier-Stokes equations and MHD, we make 
the following assumption:  
  \begin{ass}\label{ass1}
  \textup{(H0)} The matrices  $A_j$ and $B_{j,k}$ are $C^\infty$ $N \times N$ real matrices of 
  the variable $u \in \cU^* \subset \RR^N$.  Moreover, for all
   $u \in \cU^*$,
   the matrix $A_0(u)$ is invertible.

   \textup{(H1)} Possibly after a change of variables $u$ and multiplying the system
   on the left by an invertible  constant-coefficient matrix,   there is $N' \in \{ 1, \ldots,  N \} $
   and there are coordinates
   $u = (u^1, u^2) \in \RR^{N-N'} \times \RR^{N'}$
   and  $f = (f^1, f^2) \in \RR^{N-N'} \times \RR^{N'}$ such that the following 
block structure condition   is satisfied : 
  \begin{equation}
  \label{struc1}
  A_0(u) := f'_0(u) =  \begin{pmatrix}A_0^{11}&0 \\A_0^{21}&A_0 ^{22}\end{pmatrix},
  \quad
  B_{jk} (u) =\begin{pmatrix}0&0 \\0 &B_{jk}^{22}\end{pmatrix},
  \end{equation}
   \end{ass}

We refer to \cite{GMWZ4} or  \cite{Zhandbook} for further comments and explanations.  
  From now on we work with variables $u = (u^1, u^2) \in \cU^*$ such
that \eqref{struc1}   holds. We set
\begin{align}\label{y2}
A_j=f'_j,\quad \oA_j=A^{-1}_0A_j,\quad \oB_{j,k}=A^{-1}_0B_{j,k},
\end{align}
and systematically use the notation $M^{\alpha \beta} $ for the
sub-blocks of a matrix $M$ corresponding to the splitting $u =
(u^1, u^2)$.
Note that 
  \begin{equation}
  \label{struc2}
  \oB_{j, k}(u) := A_0(u)^{-1}
  B_{jk} (u) =\begin{pmatrix}0&0 \\0 &\oB_{jk}^{22}(u)
  \end{pmatrix}.
  \end{equation}
The triangular form of the equations also reveals the importance
of the $(1,1)$ block which plays a special role in the analysis : 
 \begin{equation}
 \label{11block}
L^{11} (u,\D) = \sum_{j=0}^d A^{11}_j(u) \D_j \,, \quad
\mathrm{or} \quad \oL^{11} (u,\D)= \big( A^{11}_0(u)\big)^{-1}
L^{11} (u,\D) . 
\end{equation}
In this spirit, \emph{the high-frequency principal part} of the equation is 
\begin{equation}
\label{princpart}
\left\{ \begin{array}{l}
     \overline L^{11} (u, \D)  u^1     \\
     \D_t  u^2 - \eps \overline B^{22} (u, \D) u^2
\end{array}\right. 
\end{equation}
with $\oB^{22}(u, \xi) =
\sum^d_{j,k=1} \xi_j \xi_k \oB^{2,2}_{j,k}(u) $.  We refer to Lemma~\ref{lemspec2} 
for a more detailed account of this notion of principal part. 
The first natural hypothesis is that 
$L^{11}(u, \D) $ is hyperbolic and $\D_t -  \overline B^{22} (u, \D)$ is parabolic 
in the direction $dt$. 

\begin{ass}\label{ass2}

\textup{(H2)}  There is $c > 0$ such that for all $u \in \cU^*$
and  $\xi \in \RR^d$, the eigenvalues of   $\oB^{22}(u, \xi)$  satisfy $\re \mu
\ge c \vert \xi \vert^2 $.

\textup{(H3)} For all  $u  \in \cU^*$ and   all $\xi \in \RR^d
\backslash\{0\}$,   $\oA^{11}(u, \xi) = \sum^d_{j=1}
\xi_j \bar A^{11}_j (u) $ has only  real  eigenvalues. 
\end{ass}

For the applications we have in mind such as Navier-Stokes and MHD, 
the operator $\oL^{11}$ is a transport field and (H3) is trivially satisfied.

Next we assume that the inviscid equations are hyperbolic and  that  
Kawashima' s \emph{genuine coupling} 
condition  is satisfied  for $u$,
 in some open subdomain $\cU \subset \cU^*$.  Let
\begin{align}\label{y3}
\oA (u, \xi) = \sum^d_{j=1} \xi_j \oA_j (u) \quad \text{ and
}\quad \oB(u,\xi)=\sum^d_{j,k=1}\xi_j\xi_k\oB_{j,k}(u).
\end{align}

\begin{ass}\label{ass3}

\textup{(H4)}  There is $c > 0$ such that for $u \in \cU$ and  $\xi \in \RR^d $,
the eigenvalues of $ i \oA (u, \xi)  + \oB (u, \xi ) $
satisfy
\begin{equation}\label{sdiss}
\re \mu \ge  c   \frac{ \vert \xi \vert^2  }{ 1 + \vert \xi\vert^2} \,.
\end{equation}
\end{ass}

\begin{rem}
\label{remhyp}
\textup{(H4) implies \emph{hyperbolicity} of the inviscid equation: 
for all $u \in \cU$ and  $\xi \in \RR^d\backslash\{0
\}$ the eigenvalues of $\oA (u, \xi)$ are real. 
The set $\cU$ may be thought of as the ``hyperbolic set'' where
interior, inviscid solutions are be constructed, and the larger
$\cU^*$ as the ``hyperbolic--parabolic'' set where exterior,
boundary layer solutions are to be constructed, matching
$\cU$ to boundary values in a multi-scale expansion.
In contrast with \cite{GMWZ4} and \cite{Zhandbook}, 
we do not assume here that the eigenvalues of $\oA$ have constant multiplicity. 
It is precisely the aim of this paper to substitute weaker conditions, 
allowing us to treat the case of MHD. 
}

\end{rem} 

Symmetric systems play an important role, and symmetry will be an important 
assumption in some of our results.  In particular, the Assumption  (H4) is satisfied when 
the following conditions are satisfied (see \cite{KaS1, KaS2}): 
\begin{defi}
\label{defsymm}
The system $\eqref{visceq}$ is said to be symmetric dissipative if there exists a real matrix $S(u)$, which depends smoothly 
on $u \in \cU $, such that for all $u \in \cU$ and all $\xi \in \RR^d \backslash\{0\}$, 
the matrix $S (u) A_0(u) $ is symmetric definite positive,  $S(u)  A (u, \xi) $ is symmetric and 
the symmetric matrix 
 $\re  S(u)  B(u, xi) $ is non negative with kernel of dimension $N-N'$. 
\end{defi} 

\bigbreak

We consider a boundary value problem for \eqref{visceq} and 
the model case of a  half space, 
  which is given by 
 $\{ x > 0 \}$,   in some  coordinates  
 $(y_1, \ldots, y_{d-1}, x)$ for  the space variables.    
  We assume that the boundary is not characteristic 
  both for the viscous and the inviscid equations. 
 The principal term of the viscous equation  is block diagonal  
 as indicated in \eqref{princpart}   The $B^{22}$ block is noncharacteristic by 
  (H2).  
  Restricting $\cU^*$ to a component where the profiles will take 
their values, the condition for the $\oA^{11}$ block reads

\begin{ass}
\label{assnoncar2}
$\cU^*$ is connected and for all $u \in \cU^*$, $\det    A^{11}_d  (u)  \ne 0 $.
\end{ass} 

 For the inviscid equation, 
 restricting $\cU$ to  the  component where the hyperbolic solutions will take their value, 
  the condition reads
 
 \begin{ass}\label{assnoncar} $\cU$ is connected and 
 for all $u \in \cU$, $\det \big( A_d (u) \big) \ne 0 $. 
 \end{ass} 
  
 By Assumption  (H3)  and Remark~\ref{remhyp}, 
$\overline A^{11}_d(u)$ and $\overline A_d(u)$ have only real eigenvalues, which by   Assumptions \ref{assnoncar} and  \ref{assnoncar2} never vanish. This leads  to two important 
indices :

\begin{nota}\label{defNi}
 \textup{With assumptions as above, $N_+$  denotes the number of positive eigenvalues of 
 $\overline A_d(u)$ for 
$u \in \cU$ and 
$N^1_+$   the number of positive eigenvalues of $\overline A^{11}_d(u) $ for 
$u \in \cU^*$.  We also set $N_b = N' + N^1_+ $. } 
\end{nota}

  The block structure \eqref{princpart} suggests that  $N_b   $   is the correct 
number of boundary conditions for the well posedness of  \eqref{visceq}, 
for solutions with values in $\cU^*$.  
Indeed, the high-frequency decoupling \eqref{princpart} suggests 
$N'$ boundary conditions for $u^2$  and $N^1_+ $ boundary conditions for $u^1$.
On the other hand, $N_+$ is the correct number of boundary 
conditions for the inviscid equation for solutions with values in $\cU$.  
Thus we supplement \eqref{visceq} with boundary conditions 
\begin{equation}
\label{viscbc}
\Upsilon  (u, \eps \D_y u^2, \eps  \D_x u^2 )_{ | x = 0 }   = 0 . 
\end{equation}
 Without pretending to maximal generality, we assume that they 
 decouple into zero-order boundary conditions for $u^1$ and zero-order 
and first-order conditions for $u^2$: 
 \begin{equation}
\label{viscbcd}
\left\{\begin{array}{l}
\Upsilon_1   (u^1 )_{ | x = 0 }   = 0 , \\  
\Upsilon_2   (u^2 )_{ | x = 0 }   = 0 , \\  
\Upsilon_3  (u,   \eps \D_y u^2,   \eps  \D_x u^2 )_{ | x = 0 }   = 0 .
\end{array}\right.
\end{equation}
 with 
 \begin{equation*}
\Upsilon_3  (u,  \D_y u^2,    \D_x u^2 ) =
  K_d  \D_x u^2 +  \sum_{j=1}^{d-1}  K_j (u)  \D_j u^2 . 
\end{equation*}
 
\begin{ass}
\label{assbc} 
$ \Upsilon_1 $, $\Upsilon_2$ and $\Upsilon_3$  are  smooth functions of their 
arguments  with values in  $\RR^{N^1_+}$,  $\RR^{N'  - N'' }$ and $\RR^{N''}$
respectively, where $N'' \in \{ 0, 1, \ldots, N'\}$. 
Moreover, $K_d$ has maximal rank $N''  $ and for all 
$u \in \cU^*$ the Jacobian matrices 
$\Upsilon'_1(u^1)$  and $\Upsilon'_2(u^2)$ have  maximal rank $N^1_+$
and $N'- N''$ respectively. 

\end{ass}


\subsection{Profiles and  inviscid boundary conditions}

To match constant solutions $\underline u$ of the inviscid problem to 
solutions satisfying the boundary conditions, one looks 
for  exact solutions  of \eqref{visceq} \eqref{viscbc} of the form: 
\begin{equation}
\label{layprof} 
u_\eps(t,y,x) = w \Big( \frac{x  }{\eps} \Big)   ,  
\end{equation}
such that  
\begin{equation}
\label{endstate}
 \lim_{z \to + \infty} w(z) = \underline  u \,.
\end{equation}
The equation for $w$ reads 
\begin{equation} \label{layeq}
\left\{ \begin{aligned}
&  A_d (w) \D_z w   - \D_z \big( B_{d,d}(w) 
  \D_z w  \big) = 0, \quad   z \ge 0, \\
&  \Upsilon (w, 0, \D_z w^2) _{ | z = 0 }  =  0   .
  \end{aligned}  \right. 
\end{equation}
Solutions are called \emph{layer profiles}. 
This equation can be written as a first order system for
$ U = ( w, \D_z w^2)$, which is nonsingular if and only if 
$A^{11}_d $  is invertible (this  indicates the strong link between 
Assumption~\ref{assnoncar2} and the ansatz \eqref{layprof}):  
\begin{equation} 
\label{layeq2}
\begin{aligned}
   \D_z w^1   &  =   -  ( A_d^{11} )^{-1} A^{12}_d  w^3  , \\
   \D_z w^2   &  =   w ^3   , \\ 
 \D_z \big( B_{d,d}   w ^3 )  & =
 \big( A_d^{22} - A_d^{21} ( A_d^{11})^{-1} A_d^{12} \big) w^3 , 
\end{aligned}
\end{equation}
and the matrices are evaluated at $w = (w^1, w^2)$. 

The natural limiting boundary conditions for the inviscid problem read
\begin{equation}
\label{hypbc}
u_{| x = 0}  \in \cC, 
\end{equation}
where  $\cC$ denotes the set of end points $\underline u$ such that there is a layer profile  
$w \in C^\infty(\overline \RR_+ ; \cU^*) $ satisfying  \eqref{endstate} \eqref{layeq}.   
The properties of the set $\cC$ as well as the stability analysis of 
\eqref{layeq} depend on the spectral properties of the linearized equations from \eqref{layeq}
near $w(z)$.  In particular we will discuss the notion of 
\emph{transversality}  for the profile $w$ (see \cite{MZ1, Met3}). 
However, to avoid repetitions and prepare the multidimensional stability analysis, 
we enlarge the framework and consider the  multidimensional linearized equations 
from the full system \eqref{visceq} 
near  solutions \eqref{layprof}.

For further use, it is convenient to enlarge the class of functions $w$: 
consider a function  $C^\infty(\overline \RR_+ ; \cU^*) $  
which converges at an exponential rate  to and end state $\underline u \in \cU$: 
 there is 
$\delta > 0$ such that for all $k \in \NN $
\begin{equation}
\label{expconv}
 e^{\delta z  } \big|  \D_z^k ( w(z) - \underline u) \big |   \in L^\infty (\overline \RR_+). 
\end{equation}
We refer to such a function as a \emph{profile}; it need not be a solution of 
\eqref{layeq}, though it will be in applications. Note that solutions of 
\eqref{layeq} \eqref{endstate} satisfy the exponential convergence above. 

Consider the linearized equations from \eqref{visceq} \eqref{viscbc}  around $u_\eps = w(x/ \eps) $:
\begin{equation}
\label{linq}
\cL'_{u_{\eps}} \dot u  = \dot f , \quad 
 \Upsilon' (  \dot u, \eps \D_y \dot u, \eps \D_x \dot u)  _{|
x = 0 } = \dot g. 
\end{equation}
Here $\Upsilon'$ is the differential of $\Upsilon $ at $(w(0), 0, \D_z w(0))$. 
  $\cL'_{u_{\eps}}$ is a differential operator with coefficients that 
  are smooth functions 
 of   $z := x/ \eps$. Factoring out $\eps^{-1}$  
 it also appears as an operator in $\eps \D_t, \eps \D_y, \eps \D_x$: 
 \begin{equation}
\label{linq2}
\cL'_{u_{\eps} }= \frac{1}{\eps} 
L  \Big( \frac{x}{\eps}, \eps \D_t, \eps \D_y, \eps \D_ x \Big).
\end{equation}   
It  has constant coefficients 
in $(t,y)$, and following the usual theory of constant-coefficient evolution equations,
one performs a  Laplace-Fourier  transform in $(t, y)$, with frequency variables 
denoted by $\tilde \gamma + i \tilde \tau$ and $\tilde \eta$ respectively, yielding the 
systems 
 $$
 \frac{1}{\eps} L \Big(   \frac{x}{\eps},   \eps(\tilde \gamma + i \tilde \tau) , i \eps \tilde \eta, \eps \D_x \Big).
 $$ 
 Next, we introduce explicitly the fast variable $z = x/ \eps$, rescale 
 the frequency variables as 
$  \zeta =  (\tau, \eta, \gamma) 
= \eps (\tilde \tau, \tilde \eta, \tilde \gamma)  $, 
and multiply the equation by $\eps$, revealing the 
 equation 
 \begin{equation}
\label{linq3}
L( z,   \gamma+ i   \tau, i    \eta, \D_z)    u =    f, \quad 
\Upsilon'  (   u,  i \eta   u, \D_z   u ) _{| z = 0} = g, 
\end{equation}
\begin{equation}
\label{linL}
L      = -  \cB (z)   \D^2 _z   +  \cA (z, \zeta) \D_z + \cM (z, \zeta),     
\end{equation}
with in particular, $\cB(z) = B_{d,d} (w(z))$ and 
$\cA^{11} (z, \zeta) = A^{11}_d(w(z))$. We do not give here the explicit form of 
$\cA$ and $\cM$. Using   (H2) and  Assumption~\ref{ass2}, the equation is written as a first order
system     
 \begin{equation}
\label{linq4} 
\D_z U  = \cG(z, \zeta) U + F  , \quad   \Gamma (\zeta)  U _{| z = 0} = g, 
\end{equation} 
where
 \begin{equation}
 \label{defU} 
U = {}^t ( u , \D_z u^2)  = (u^1, u^2, \D_z u^2)   \in \CC^{N+ N'},   
 \end{equation}
  \begin{equation}
 \label{defF}
F  =    \big((\cA^{11}(z))^{-1} f^1, 0 ,  (\cB^{22}(z))^{-1} (-  f^2+ \cA^{21}(z)(\cA^{11}(z))^{-1} f^1)  \big) .      
 \end{equation}

 The analysis of this equation depends on the size of the frequencies 
 $\zeta$.  
When $\zeta$ is large, the character of the equations is dominated
by the high-frequency principal part \eqref{princpart},
and we use a slowly-varying-coefficients analysis (related to
the ``tracking lemmas'' of \cite{Z1, Zhandbook}) based on the relatively
slow rate of change of coefficients compared to the size of the frequency;
see \cite{GMWZ3, GMWZ4} and Section \ref{S4} below.
For small or bounded frequencies $\zeta$, we use   the conjugation lemma of 
\cite{MZ1}. The condition \eqref{expconv} implies that there is $\delta > 0$ and an end state matrix
$ G(\underline u, \zeta) $, depending on the endstate $\underline u $  of $w$, such that  
\begin{equation}
\label{expconvG}
\D_z^k ( \cG( z, \zeta) - G(\underline u, \zeta) ) = O (e^{ - \delta z }). 
\end{equation}

\begin{lem}
\label{lemconj}
Given $\underline \zeta \in \RR^{d+1}$, 
there is a smooth invertible matrix  $\Phi (z, \zeta)$ for 
$z \in \overline \RR_+$ and $\zeta$ in a neighborhood of $\underline \zeta$, such that 
$\eqref{linq3}$ is equivalent to   
 \begin{equation}
 \label{linq5}
\D_z  \tilde U    = G(\underline u,   \zeta) \tilde U   +  \tilde F    , \quad  
\tilde  \Gamma (\zeta)  \tilde  U  _{| z = 0} = g. 
\end{equation}
with  $  U  = \Phi(z, \zeta) \tilde U  $,  $  F = \Phi(z, \zeta) \tilde F  $ and 
$\tilde \Gamma (\zeta)  = \Gamma(\zeta)  \Phi (0, \zeta)$. 
In addition, $\Phi $ and $\Phi^{-1}$ converge the identity  matrix at an exponential rate when  
$z \to \infty$. 

Moreover, if the coefficients of the operator and $w$ depend smoothly on extra parameters $p$ 
(such as the end state $\underline u$),  then $\Phi$ can also be chosen to depend smoothly on $p$, on a neighborhood of a given 
$\up$.  
\end{lem}

\begin{rem}
\label{rem0}
\textup{The linearized profile equations from \eqref{layeq} around $w$, are exactly \eqref{linq3} at the frequency
$\zeta = 0$. In particular, Lemma~\ref{lemconj} implies that these equations are conjugated 
to   constant-coefficient equations, via the conjugation by  $\Phi(\cdot , 0)$.}
\end{rem}

Next we investigate the spectral properties of the matrix $G$. 
Below, $\RR^{d+1}_+ $ denotes the open half space 
$\{ \zeta = (\tau, \eta, \gamma) :  \gamma > 0 \}$ and 
$\overline \RR^{d+1}_+ $ its closure $\{\gamma \ge 0\} $. 
We also introduce the matrices  
\begin{eqnarray}
\label{defP0}
P_0 (u) & :=&  (B^{22} )^{-1}  \big( A_d^{22} - A_d^{21} ( A_d^{11})^{-1} A_d^{12} \big),  
\\
  \label{defH0}
   H_0(u, \zeta)   &  :=  &- (A_d(u))^{-1}  \Big( ( i \tau  + \gamma ) A_0(u) +
\sum_{j=1}^{d-1} i \eta_j A_j(u) \Big) . 
\end{eqnarray}

 \begin{lem}
\label{lemspec}
i) For $u \in \cU$,   $P_0(u)   $ 
has no eigenvalue on the imaginary axis.  We denote  by  $N^2_-$   the number of its eigenvalues 
in $\{ \re \mu < 0 \}$.  

\quad ii)  For $u \in \cU$ and $\zeta \in \overline \RR^{d+1}_+ \backslash\{ 0 \}$,  $G(u, \zeta)$ has no eigenvalue on the imaginary axis.  The number of its eigenvalues, counted with their multiplicity, 
in $\{ \re \mu < 0 \}$  is equal to  $  N_+ + N^2_- = N_b := N' + N^1_+ $. 

\quad iii) For a given $\underline u\in \cU$, there are smooth matrices 
   $V(u, \zeta) $  on a neighborhood of $(\underline u, 0)$ such that 
 \begin{equation}
\label{block1}
V^{-1} G    V  = \left(\begin{array}{cc }
H  & 0  \\ 0 & P    
\end{array}\right)  
\end{equation}
with $H(u, \zeta)$ of dimension $N \times N$, 
 $P(u , \zeta)$ of dimension $N' \times N'$,  and 
 
 a)  the eigenvalues of
$P$ satisfy $\vert \re \mu \vert \ge c$ for some $c > 0$, 
 
 b) there holds 
 \begin{equation}
 \label{eqH0}
 H(u, \zeta) =  H_0(u, \zeta)   + O(\vert \zeta \vert^2)
 \end{equation}
   
  c ) at $\zeta = 0$,  $V$ has a triangular form
  \begin{equation}
\label{Vtri}
V(u, 0 ) = \begin{pmatrix}
     \Id  &    \overline V  \\
     0  &  \Id 
\end{pmatrix}. 
\end{equation}
\end{lem}

\begin{proof}  
  i)   Take  $ u \in \cU$. If $ v^2 \in \ker P_0(u) $, then $ {}^t\big( -  (A_{d}^{11})^{-1} A_d^{12} v^2 , v^2) \in \ker A_d$, implying that $0$ is not an eigenvalue of $P_0$. Similarly,  if  
$i \xi $ is an eigenvalue of $P$ then 
$0 $ is an eigenvalue of $i \xi \oA_d + \xi^2 \oB_d$, which is impossible by (H4) 
if $\xi  \ne 0$ is real.  

ii)   Direct computations show that 
$G(u, \zeta) = G_d(u, \zeta)^{-1} M(u, \zeta)$ with 
\begin{equation*}
G_d(u, \zeta) = \begin{pmatrix}
    -  \tilde A_d   &   \tilde B_d   \\
     J  &   0 
\end{pmatrix}, \quad 
M = \begin{pmatrix}
    \tilde M   &   0_{N \times N'}  \\
     0_{N' \times N}   &  \Id_{N' \times N'}
\end{pmatrix}
\end{equation*}
 with, in the splitting $u = (u^1, u^2)$,  
 \begin{equation*}
 \tilde B_d (u)  = \begin{pmatrix}
      0_{N - N' \times N'}    \\
      \oB^{22}_{d, d}   (u) 
\end{pmatrix} , \quad  
J = \begin{pmatrix}
     0_{N' \times N-N'}  & \Id_{N' \times N'}     
     \end{pmatrix}. 
\end{equation*}
and 
 $$
\left\{
  \begin{aligned}{}&
\tilde A (u,   \zeta )  =     A_d (u)   -
 \sum_{j=1}^{d-1}   i \eta_j
( B_{j,d} (u)   + B_{d,j} (u)  )
  \\
&
\tilde M (u ,  \zeta)  =
   ( i \tau  + \gamma) A_0(u)   +
\sum_{j=1}^{ d-1}   i \eta_j  A_j(u) 
    +
\sum_{j,k=1}^{d-1} \eta_j \eta_k  B_{j,k}(u) 
 \,.
\end{aligned}
\right.
$$
In particular,   $i \xi$ is an eigenvalue of $G(u, \zeta)$ if and only if 
$\gamma + i \tau$ is an eigenvalue of  $  i \oA (\eta, \xi) + \oB (\eta, \xi)$, which, by (H4),  implies 
either that $\gamma < 0$ if $\xi$ is real and $(\eta, \xi ) \ne 0$ or that 
$\zeta = 0$. 

Thus  $G(u, \zeta)$ has no eigenvalues on the imaginary axis and the number 
$\tilde N$ of eigenvalues in $\{ \re \mu < 0 \}$ is constant 
for 
$u \in \cU$ and $\zeta \in \overline \RR^{d+1}_+ \backslash\{0\}$. 
That this number is equal to $N_b = N^1_+ + N'$ is a consequence 
of  the high-frequency analysis in Lemma~\ref{lemspec2} below (see also  Lemma 1.7 in \cite{Zhandbook}).
 
iii) Because $\tilde M (u, 0) = 0$ and $\tilde A(u, 0) = A_d(u)$,  there holds
\begin{equation}
\label{G0}
G(u, 0 ) = \begin{pmatrix}
0_{N \times N}  & \begin{pmatrix}
             - (A_d^{11})^{-1} A_d^{12}
      \\ \Id_{N' \times N'}  
\end{pmatrix} \\ 0_{N' \times N}  & P_0(u)
\end{pmatrix}
\end{equation}
  Since  $P_0$ is invertible, $G$ can be smoothly conjugated to a   block diagonal matrix
   as in \eqref{block1}, with $V$ satisfying \eqref{Vtri} and $H(u, 0) = 0$.   
   More precisely, the matrix $\overline V$ is 
   \begin{equation}
\label{Vtri2}
\overline V =  \begin{pmatrix}
             - (A_d^{11})^{-1} A_d^{12} P_0^-1
      \\ P_0^{-1} \end{pmatrix} 
\end{equation}
The expansion \eqref{eqH0} can be easily obtained by standard perturbation expansions, 
and we refer to  Lemma~\ref{lem642b} below for 
 a more precise version.  
 
For $\zeta$ small, the number of eigenvalues of $P$ in $\{\re \mu < 0 \}$ is equal 
to $N^2_-$, and for $\gamma > 0$, the number of eigenvalues of 
$H_0 (u, \zeta)$ in the negative half space is constant, by hyperbolicity, and equal to $N_+$. 
This implies  that $\tilde N = N_+ + N^2_-$. 
\end{proof}

\bigbreak

Similarly, one considers the linearized equations from the inviscid hyperbolic 
problem $\cL_0 (u) = 0$ around the constant solution $\underline u$:
\begin{equation}
\label{linhypeq}
\cL'_{0, \underline u} \dot u = \dot f. 
\end{equation}
After performing a  Laplace-Fourier transform,     this  equation reads
\begin{equation}
\label{linhypeq2}
L_0 ( \underline u,  \gamma + i \tau, i \eta, \D_x) u  = f 
\end{equation}
or, with $H_0$ defined at \eqref{defH0},   
\begin{equation}
\label{linhypeq1}
\D_x u  = H_0( \underline u, \zeta) u +  A_d^{-1}(\underline u) f. 
\end{equation}

\bigbreak

An important property for profiles is the notion of 
\emph{transversality} (see \cite{MZ1} or \cite{Met3} for the case of total viscosity). 
It concerns  
   the linearized equations from \eqref{layprof} around  $w$. As mentioned  in Remark~\ref{rem0}, they correspond exactly to the 
  first order system \eqref{linq3} with $\zeta = 0$.  We abbreviate  the homogeneous 
   problem 
  as 
\begin{equation}
\label{linlayeq}
\left\{ \begin{aligned}
&  L (z, 0, \D_z)   \dot w  = 0 , \qquad   z \ge 0, \\
&  \Upsilon'   (\dot w, 0,  \D_z \dot w^2) _{ | z = 0 }  =  0    .
\end{aligned}\right. 
\end{equation}
   A corollary of Lemmas~\ref{lemconj} and \ref{lemspec} is that the   solutions of the homogeneous equation 
  $  L(z, 0, \D_z ) \dot w = 0$   form a space   of dimension $N + N' $, parametrized by 
  $(u_H, u_P) \in \CC^N \times \CC^{N'}$ : 
  \begin{equation}
\label{solborn}
\dot w(z)  = \Phi_H(z) u_H  + \Phi_P(z) e^{ z P_0(\underline u)} u_P 
\end{equation}
where  the matrices  $\Phi_H(z)$ and $\Phi_P(z)$ are smooth and bounded
   on $\RR_+$ and $ \Phi_H (z) \to \Id$ as $z \to \Id$. 
   The solution is bounded if and only if $u_P$ belongs 
   to  the negative  space $\EE^-(P_0(\underline u))$
   of $P_0(\underline u)$, that is the invariant space of $P_0(\underline u)$ associated 
   to the spectrum lying in $\{ \re \mu < 0 \}$; thus the space $\cS$ of bounded 
   solutions has dimension $N + N^2_-$. 
   The space  of solutions that tend to zero at infinity, denoted by $\cS_0$, has dimension
   $N^2_-$,  corresponding to the conditions   $u_H = 0$  and $u_P \in \EE^-(P_0(\underline u))$.

   The boundary conditions in \eqref{linlayeq} read 
   \begin{equation}
\label{newbclayeq}
 \underline \Gamma_{  H} u_H + 
\underline \Gamma_{  P} u_P  := \underline \Gamma  (\dot w,  \D_z \dot w^2) _{ | z = 0 } = 0 . 
\end{equation}

  \begin{defi}
  \label{deftrans}
  The profile $w$ is said to be transversal if 
  
  \quad i)  there is no nontrivial  solution  
  $ \dot w \in \cS_0$ which satisfies the boundary conditions
  $  \underline \Gamma (\dot w, \D_z \dot w^2)_{| z = 0} = 0$, 
  
  \quad ii)  the mapping $\dot w \mapsto \underline \Gamma (\dot w, \D_z \dot w^2)_{| z = 0}$
  from $\cS$ to $\CC^{N_b}$  has  rank $N_b$.

  \end{defi}
 
 Equivalently, it means that  
 $\ker \underline \Gamma_{ P} \cap \EE^- (P_0(\underline u)) = \{ 0 \} $ and that 
 the rank of the matrix $(\underline \Gamma_{ H}, \underline \Gamma_{ P})$ from 
 $\CC^N \times \EE^-(P_0(\underline u))$ to $\CC^{N_b}$  is $N_b$. 
 
  If the profile satisfies condition i),  there is a decomposition 
 \begin{equation}
\label{dedbc}
\CC^{N_b}  = \FF_H  \oplus \FF_{ P}  , \quad    \FF_P := \underline \Gamma_{ P}  \EE^-(P_0(\underline u))
\end{equation}
 with $\dim \FF_H = N_+$ and $\dim \FF_{0, P} = N^2_-$. 
 Denote by $\pi_H$ and $\pi_P$ the projections associated to this splitting.

 For $\dot w \in \cS$ given by \eqref{solborn}, one can eliminate $u_P$ from the boundary conditions \eqref{newbclayeq} and write   them 
 \begin{equation}
\label{renewbclayeq}
\underline \Gamma_{  red} u_H = 0, \quad  u_P = R_{0, P} u_H, 
\end{equation}
 with 
 \begin{equation}
\label{linhypbc1}
\underline \Gamma_{red} := \pi_H \underline \Gamma_{  H}, \quad 
R_{0, P} := - (\underline \Gamma_{ P})^{-1} \pi_P \underline \Gamma_{  H}
\end{equation}
 and $ (\underline \Gamma_{ P})^{-1}$ is the inverse of the mapping  $ \underline \Gamma_{ P} $ from 
 $\EE^-(P_0(\underline u))$ to $\FF_{0, P}$. 
 
 With these notations,  ii) means 
 that $\underline \Gamma_{  red}$ has rank $N_+$.  
 Its kernel $\ker \underline \Gamma_{  red}$ is 
 the space   of $\dot u  \in \RR^d$ such that 
 there is a solution of $\dot w  $ of \eqref{linlayeq} with end point $\dot u$. 
 It  has dimension $N - N_+ $. 
 
 \begin{rem}
 \label{remhypbc}
\textup{When $w$ is a layer profile, solution of \eqref{layeq}, the transversality 
condition implies that near the end point $\underline u$, the set $\cC$ in \eqref{hypbc}
which describes the limiting hyperbolic conditions is a smooth manifold
of dimension $N_- = N- N_+$ and  $\ker \underline \Gamma_{  red}$ is the tangent space 
to $\cC$ at 
 $\underline u $.  Therefore, the natural boundary condition for the linearized 
 hyperbolic equation, and in particular for  \eqref{linhypeq},  are }
 \begin{equation}
\label{linhypbc}
\underline \Gamma_{  red}   u  =  h. 
\end{equation}
\end{rem}

\bigbreak

\subsection{Evans functions and Lopatinski determinant}

For a given $\zeta \in \overline \RR^{d+1}_+\backslash\{0\}$, we now investigate the well-posedness of equation \eqref{linq3} or equivalently \eqref{linq4}
or \eqref{linq5}. 
Introduce the space 
$\EE^- (\zeta)$ of initial conditions  $(u(0), \D_z u^2(0))$   (or equivalently $U(0)$) such that 
the corresponding solution of $ L( z, \zeta, \D_z ) u = 0$ 
(or $\D_z U - \cG (z, \zeta) U = 0$) is exponentially decaying at $+ \infty$. 
Lemmas~\ref{lemconj} and \ref{lemspec}  show that 
\begin{equation}
\label{EtE}
\EE^-(\zeta)  = \Phi (0, \zeta)    \EE^- (G(\underline u, \zeta) )
\end{equation}
  where we use the following notations:
  
  \begin{nota}
\textup{Given  a square matrix $M$, $\EE^-(M)$ [resp. $\EE^+(M)$ denotes the invariant space of $M$ associated 
  to the spectrum of $M$ contained in $\{\re \mu  < 0 \}$ [resp $\{\re \mu  > 0 \}$]. }
\end{nota}

In particular,  by Lemma~\ref{lemspec},  $\EE^- (\zeta)$ is a smooth vector bundle for   $\zeta \in \overline \RR^{d+1}_+\backslash\{0\}$ and $\dim (\EE^- (\zeta) )  = N_b $. 

The problems \eqref{linq3}, \eqref{linq4} or \eqref{linq5} are well posed if and only if 
\begin{equation}
\label{weakstab}
\EE^- (\zeta) \cap \ker \Gamma(\zeta) = \{ 0 \}  \quad \mathrm{ or } 
\quad  \EE^- (G(\underline u, \zeta) ) \cap \ker \tilde  \Gamma(\zeta) = \{ 0 \}. 
\end{equation}
Note that, because the rank of $\Gamma$ is at most $N_b$ and the dimension 
of $\EE^-$ is $N_b$, this condition implies  and is equivalent to 
\begin{equation}
\label{weakstab2}
\CC^{N + N'} = \EE^- (\zeta) \oplus  \ker \Gamma(\zeta)    \quad \mathrm{ or } 
\quad 
\CC^{N + N'} =   \EE^- (G(\underline u, \zeta) ) \oplus  \ker \tilde  \Gamma(\zeta) . 
\end{equation}

The \emph{Evans function} is defined as 
\begin{equation}
\label{Evf}
D( \zeta) = \big|  \det{}_{N+N'}  \big( \EE^- (\zeta) ,  \ker \Gamma(\zeta)\big) \big| 
\end{equation}
where, for subspaces $\EE$ and $\FF$ of $\CC^n$, $\det_n (\EE, \FF) $ is equal 
to $0$ if $\dim \EE + \dim \FF \ne n$ and is the  
$n \times n$ determinant formed by  
orthonormal bases in $\EE$ and $\FF$ if $\dim \EE + \dim \FF = n$. 

\begin{rem}
\textup{The definition of the determinant above  depends on   choices of bases.  Note that changing bases in $\EE$ and $\FF$
changes the determinant by a complex number of modulus one, thus 
leaves $|\det (\EE, \FF) |$ invariant.  
But it also  depends on the choice of  a scalar  product  on $\CC^n$. Changing the  scalar products  (or changing of bases in $\CC^n)$) 
changes the function $\det (\EE, \FF)$ to a new function $\widetilde \det(\EE, \FF)$ such that 
$c  |  \det (\EE, \FF) | \  \le |  \widetilde    \det (\EE, \FF) |  \le c^{-1} |  \det (\EE, \FF) | $ where $c > 0$ is independent of the spaces $\EE$ and $\FF$. We will denote by}
 \begin{equation}
\label{equivEvf}
\det \approx \widetilde \det \quad  \mathrm{ or } \quad  D \approx \widetilde D
\end{equation} 
\textup{this property. In particular,  the definition of $D$  is independent  of the choice of orthonormal bases in $\EE^-$ and $\ker \Gamma$ and  all the uniform stability conditions stated below are 
independent of the choice of the scalar product. } 
\end{rem}

\begin{rem}
\label{Evparam}
\textup{If the coefficients of the operator and the profile depend smoothly on parameters
$p$, then the Evans function is also a smooth function of the parameters.}
\end{rem}

These notations being settled,  the \emph{weak stability} condition, which is a necessary condition for well posedness in Sobolev spaces of   \eqref{linq}, 
 reads:

\begin{defi} Given a profile $w$, the linearized equation $\eqref{linq}$ 
satisfies 
the weak spectral stability condition  if  
$D (\zeta) \ne 0$  for all $\zeta \in \overline \RR^{d+1}_+\backslash\{0\}$. 

\end{defi}

\bigbreak

The next lemma  is useful and elementary.

\begin{lem}
\label{lemustab1}   Suppose that $\EE \subset \CC^n $ and  $\Gamma : \CC^n \mapsto \CC^m $, with   $\mathrm{rank} \Gamma = \dim\EE = m$. 
If $ | \det (\EE , \ker \Gamma) |   \ge c > 0 $, then there is $C$, which depends only on
$c$ and $| \Gamma^* (\Gamma \Gamma^*)^{-1} |$ such that 
\begin{equation*}
\forall U \in \EE  \quad |  U | \le C  | \Gamma   U | . 
\end{equation*}
 Conversely, if this estimate is satisfied  then $ | \det (\EE , \ker \Gamma) |   \ge c$ 
where $c > 0  $  depends only on
$C$ and  $| \Gamma  |$. 
\end{lem}
\begin{proof} Let $\pi =  \Gamma^* (\Gamma \Gamma^*)^{-1} \Gamma$ denote the
orthogonal projector  on $(\ker \Gamma)^\perp$. Diagonalizing the hermitian form 
$(\pi e, \pi e)$, yields orthonormal bases $\{ e_j\} $  and $\{ f_j \}$ in 
$\EE$ and $(\ker \Gamma)^\perp$ respectively, such that $\pi e_j = \lambda_j f_j$ 
with $0 < \lambda_j \le 1$.  Take any basis  $\{ g_k \}$ of $\ker \Gamma$.  
Expressing the $e_j$ in the base $\{ f_k , g_l\}$, implies that 
$| \det (\EE , \ker \Gamma) |  = \prod \lambda_j$. 
Since 
  $\lambda_j  \le 1 $ for all $j$, if this determinant  is larger than or equal to $   c > 0 $,  
  then $\min \lambda_j \ge c$ and for all $e \in \EE$
    \begin{equation*}
 c | e | \le | \pi e |  \le   | \Gamma^* (\Gamma \Gamma^*)^{-1}| \, | \Gamma e |. 
\end{equation*}  
Conversely, if the estimate is satisfied, then 
$|e | \le C | \Gamma | \, | \pi e |$ since $\Gamma e = \Gamma \pi e$ for all 
$e \in \EE$. 
Therefore $\lambda_j C | \Gamma | \ge 1$  and the determinant is at least equal 
to $(C | \Gamma |)^{ - m}$. 
\end{proof}

\bigbreak 

\bigbreak
There are analogous definitions for the linearized hyperbolic problem
\eqref{linhypeq} with boundary conditions \eqref{linhypbc}. 
For $\gamma > 0$, $H_0(\underline u, \zeta)$ has no eigenvalues 
 on the imaginary axis, as a consequence of the hyperbolicity 
 assumption (see Remark~\ref{remhyp}). 
The \emph{Lopatinski determinant}  is defined for 
 $\zeta \in \RR^{d+1}_+ :=  \{ \gamma  > 0 \}$  by   
 \begin{equation}
\label{Lopdet}
D_{Lop} (\zeta ) =   \big| \det \big( \EE^- ( H_0(\underline u , \zeta) , 
\ker \underline \Gamma_{  red}  \big) \big|. 
 \end{equation}
 By homogeneity of  $H_0$, this determinant is homogeneous of 
 degree zero in $\zeta$ and  one can restrict attention to 
 $\zeta \in S^d = \{ | \zeta | = 1 \}$. 
 
\begin{defi}
\label{defstab2}
 The linearized equation 
 \eqref{linhypeq} \eqref{linhypbc} 
satisfies 
the weak spectral stability condition  if 
$D_{Lop} (\zeta) \ne 0$ for all $\zeta \in   \RR^{d+1}_+$. 
\end{defi}


\subsection{Uniform spectral stability and maximal estimates}

The weak stability conditions and the reduction to constant coefficients of Lemma~\ref{lemconj} guarantee  the well posedness of \eqref{linq3} for fixed 
$\zeta \in \overline \RR^{d+1}_+ \backslash\{0\}$ and in particular 
estimates of the form 
\begin{equation}
\label{roughest}
\|   u \|_{L^2} + \| \D_z u^2 \|_{L^2}  + | u(0) | + | \D_z u^2(0) |  \le  C(\zeta) 
(  \| f \|_{L^2}  + | g | ) . 
\end{equation}
The next step in the study of 
\eqref{linq}, is  to perform an inverse Fourier-Laplace transform and thus requires
suitable    estimates for the solutions of \eqref{linq3}, with a precise description of the 
constants in the estimate above. 

By continuity in $\zeta$,  
 the weak stability condition implies that the estimate \eqref{roughest} is
satisfied with a uniform constant $C$ when $\zeta $ remains in a compact subset 
of $ \overline \RR^{d+1}_+ \backslash\{0\}$. 
Thus the true question is   to get a detailed behavior of the estimate  when 
$\zeta \to 0$ and when $| \zeta | \to \infty$. 
 
\medbreak

\subsubsection{Low  and medium frequencies}
Consider first the \emph{low-frequency} case. Following \cite{MZ1}, 
the uniform stability condition reads: 

\begin{defi} Given a profile $w$, 
the uniform spectral stability condition for low frequencies  is satisfied 
 when  
there are $c> 0$ and $\rho_0 > 0$  such that 
$D   (\zeta) \ge c $ for all $\zeta \in \overline \RR^{d+1}_+  $ with 
$0 < | \zeta | \le \rho_0$.

\end{defi}

By Assumption~\ref{assbc}, the rank of $\Gamma(\zeta)$ is always $N_b$, and the norms  of $\Gamma(\zeta)$  and $(\Gamma \Gamma^*)^{-1}  $
are  uniformly bounded for  $\zeta$ bounded.  Thus, by Lemma~\ref{lemustab1}, 
the 
low-frequency uniform stability   condition holds if and only if there are $C$ and 
$\rho_0 > 0$ such that 
\begin{equation}
\label{stabcond}
\forall \zeta \in \overline \RR^{d+1}_+ , \ \ 0 < | \zeta | \le \rho_0 , \ \ 
\forall U \in \EE^- (\zeta): \quad |  U | \le C  | \Gamma (\zeta) U | . 
\end{equation}
  
Following \cite{MZ1}, the expected  \emph{maximal estimates}   for low and medium frequencies   for the solutions of 
\eqref{linq3} read
  \begin{equation}
\label{maxest2}
\begin{aligned}
\varphi    \| u  \|_{L^2(\RR_+)}     + 
 \| \D_z u^2    \|_{L^2(\RR_+)}      +  | u  (0& |  + |  \D_z u^2 (0) |    
 \le        
\\  &  C  \big(   \frac{1 }{\varphi }  \| f  \|_{L^2(\RR_+)}    + | g | \big)
      \end{aligned}
\end{equation}
where $\vp = (\gamma + | \zeta |^2) ^{\frac{1}{2}}$ with $C$ independent of 
$\zeta \in \overline \RR^{d+1}_+ \backslash\{0\}$, $| \zeta | \le \rho_0$. 
Note that for fixed $ | \zeta |  > 0$,  
this estimate implies  \eqref{roughest}.

The estimates \eqref{maxest2} correspond to estimates for the solutions 
of the first order system \eqref{linq4}: 
 \begin{equation}
\label{maxest2s}
\varphi    \| U^1  \|_{L^2(\RR_+)}     + 
 \| U^2    \|_{L^2(\RR_+)}      +  | U  (0  |    
 \le        
     C  \big(   \frac{1 }{\varphi }  \| F  \|_{L^2(\RR_+)}    + | g | \big)
\end{equation}
where $U = (U^1, U^2) \in \CC^N \times \CC^{N'}$. 
For the constant-coefficient system
\eqref{linq5} the expected estimates read : 
\begin{equation}
\label{maxest2ss}
\varphi    \| \tilde U^1  \|_{L^2(\RR_+)}     + 
 \| \tilde U^2    \|_{L^2(\RR_+)}      +  | \tilde U  (0  |    
 \le        
     C  \big(   \frac{1 }{\varphi }  \| \tilde F  \|_{L^2(\RR_+)}    + | g | \big). 
\end{equation}

\begin{lem}
\label{lemequivest}
The estimates $\eqref{maxest2ss} $ imply  $\eqref{maxest2ss}$ which imply 
$\eqref{maxest2}$.
\end{lem}

\begin{proof} (See \cite{MZ1}).  Clearly, \eqref{maxest2} is a particular case of 
\eqref{maxest2s} applied to source terms $F$ of the special form \eqref{defF}. 
Moreover, using the conjugation Lemma~\ref{lemconj},  
there holds $U = O(1) \widetilde U$ and $\widetilde U = O(1) U$ and similar estimates
for $F$ and $\widetilde F$. Moreover,  
$$
U^1 =  O(1) \widetilde U,    \quad  
U^2 = O(e^{ - \theta z}) \widetilde U^1 + O(1) \widetilde U^2
$$ 
with $\theta > 0$. 
We use the inequality
$$
\| e^{ - \theta z} \widetilde U^1 \|_{L^2} \ls |   \widetilde U^1(0) | 
+ \|\D_{z} \widetilde U^1  \|_{L^2}. 
$$
Moreover, the form of $G(\underline u, \zeta) $ at $\zeta = 0$ shows that 
$$
\D_z \widetilde U^1  = O(| \zeta | ) \widetilde U^1 + O(1) \widetilde U^2 + \widetilde F^1. 
$$
Therefore, 
$$
\| U^2 \|_{L^2} \ls \| \widetilde U^2 \|_{L^2}  +  | \widetilde U^1(0) | + 
| \zeta | \| \widetilde U^1 \|_{L^2}  + \| \widetilde F ^1 \|_{L^2}. 
$$
Since $|\zeta | \le \vp$, this shows that \eqref{maxest2ss}  implies \eqref{maxest2s}. 
\end{proof}

For $  \zeta  $ in a compact subset of $\overline \RR^{d+1}_+ \backslash\{0\}$, 
all these estimates are true under the weak stability condition (see e.g. \cite{MZ1}). 
Note also (taking $f = 0$ in \eqref{maxest2}) that the uniform stability condition \eqref{stabcond} is necessary 
for the validity of the maximal estimate. 
\emph{The main subject of this paper is to prove 
that the uniform stability condition implies the maximal estimate \eqref{maxest2} for low frequencies, 
under  structural assumptions on the system  weaker than in} \cite{MZ1, GMWZ3, GMWZ4}, allowing for 
instance to consider MHD.

\subsubsection{High frequencies} 
For the \emph{high-frequency} analysis, the maximal extimates that are proven in 
\cite{GMWZ4} concern homogeneous boundary conditions ($ g = 0$) and read
\begin{equation}
\label{maxesthfhom}
\begin{aligned}
(1+ \gamma)    \| u^1  \|_{L^2(\RR_+)}     + 
&\Lambda   \| u^2  \|_{L^2(\RR_+)} +  \| \D_z u^2  \|_{L^2(\RR_+)} 
\\
+ 
(1+ \gamma)^{\frac{1}{2}}     |u^1(0)  |   + 
&\Lambda^\mez   | u^2(0)| + \Lambda^{-\mez} | \D_z u^2 (0)|  \le         
\\   
  & C  \big(     \| f ^1  \|_{L^2(\RR_+)}    +  \Lambda^{- 1}  \| f ^2  \|_{L^2(\RR_+)}  \big). 
      \end{aligned}
\end{equation}
where $\Lambda$ is the natural parabolic weight 
 \begin{equation}
\label{defLambda}
\Lambda(\zeta)  = \big( 1 +  \tau^2 + \gamma^2 + | \eta |^4 \big) ^{1/4}. 
\end{equation}
The balance between the weights for $u^1$ and for $u^2$ is subtle: these components are 
decoupled in  the high-frequency 
principal system \eqref{princpart} and the weights   depend   on their actual coupling through
the  nondiagonal terms and the boundary conditions. 
Here we see the importance of the form \eqref{viscbcd} of the boundary conditions. 
Their linearized version, $\Upsilon' (u, i\eta u^2, \D_z u^2) = g$ reads
  \begin{equation}
\label{linviscbcd}
\left\{\begin{array}{l}
\Gamma_1 u^1(0) := \Upsilon'_1   (w^1(0) )\cdot u^1{} (0)   = g^1 , \\  
\Gamma_2 u^2(0):= \Upsilon'_2   (w^2(0) )\cdot u^2{} (0)   = g^2 , \\  
\Gamma_3 (\zeta) (u^2(0), \D_z u^2(0))      :=   K_d \D_z u^2 (0)  +   K_{\mathrm{tg}}  (\eta) u^2  (0)
    = g^3 .
\end{array}\right.
\end{equation}
with 
   \begin{equation}
   \label{Ktg}
 K_{\mathrm{tg}}    =  \sum_{j=1}^{d-1} i \eta_j K_j (w(0))  
\end{equation}
The complete  maximal estimate with nonvanishing  boundary source terms $g$, 
reads 
 \begin{equation}
\label{maxesthf}
\begin{aligned}
(1+ \gamma)    \| u^1  \|_{L^2(\RR_+)}     + 
&\Lambda   \| u^2  \|_{L^2(\RR_+)} +  \| \D_z u^2  \|_{L^2(\RR_+)} 
\\
+ 
(1+ \gamma)^{\frac{1}{2}}     |u^1(0)  |   + 
&\Lambda^\mez   | u^2(0)| + \Lambda^{-\mez} | \D_z u^2 (0)|  \le         
\\   
  & C  \big(     \| f ^1  \|_{L^2(\RR_+)}    +  \Lambda^{- 1}  \| f ^2  \|_{L^2(\RR_+)}  \big) 
  \\
   & +   C  \big(    (1+ \gamma)^{\frac{1}{2}}     |g^1  |   + 
 \Lambda^\mez   | g^2 | + \Lambda^{-\mez} | g^3|  \big)    
      \end{aligned}
\end{equation}
 with $C$ independent of $\zeta \in \overline \RR^{d+1}_+$   large. 
 Taking $f = 0$, this implies the  following necessary condition~:   there are 
 $C$ and $\rho_1 > 0$ such that 
\begin{equation}
\label{stabcondhf}
\begin{aligned}
 \forall \zeta & \in \overline \RR^{d+1}_+, \ \ | \zeta | \ge \rho_1, \quad 
\forall U = (u^1, u^2, u^3)  \in \EE^- (\zeta): 
\\
  (1& + \gamma)^{\frac{1}{2}}     |u^1   |   + 
 \Lambda^\mez   | u^2 )| + \Lambda^{-\mez} |u^3 |  \le    
 \\
& C \big(  (1+ \gamma)^{\frac{1}{2}}     |\Gamma_1 u^1   |   + 
 \Lambda^\mez   | \Gamma_2 u^2 )| + \Lambda^{-\mez} |\Gamma_3 (\zeta) (u^2, u^3) | \big) 
\end{aligned}
\end{equation}
This can be reformulated in terms of a  \emph{rescaled
Evans function}  (see \cite{MZ1}~:   
In $\CC^{N+ N'}$  and $ \CC^{N_b} $ introduce  the mappings 
\begin{equation}
\label{reschf}
\begin{aligned}
J_\zeta (u^1, u^2, u^3) &  := \big((1+ \gamma)^{\frac{1}{2}}   u^1,  
 \Lambda^\mez  u^2 ,   \Lambda^{-\mez} u^3 \big)  \\
 J_\zeta (g^1, g^2, g^3) & := \big((1+ \gamma)^{\frac{1}{2}}   g^1,  
 \Lambda^\mez  g^2 ,   \Lambda^{-\mez} g^3 \big)  . 
 \end{aligned}
\end{equation}
Note that $ J_\zeta  \Gamma (\zeta)  U =  \Gamma^{sc}(\zeta) J_\zeta U $ 
with 
 \begin{equation}
\label{splitGamma}
 \Gamma^{sc}  U =  \big(   \Gamma_1  u^1  ,  \Gamma_2 u^2 
         ,    K_d u^3  + \Lambda^{-1} K_{\mathrm{tg}} (\eta)  u^2  \big) . 
\end{equation}
Thus \eqref{stabcondhf} reads
\begin{equation}
\label{stabconhfsc}
\forall U \in J_\zeta \EE^- (\zeta) :\quad  
| U | \le C     | J_\zeta \Gamma (\zeta) J_\zeta^{-1} U | 
\end{equation}
Introducing  the \emph{rescaled Evans function} 
\begin{equation}
\label{rescEvf}
D^{sc}( \zeta) = \big|  \det   \big( J_{\zeta}\EE^- (\zeta) ,  J_\zeta  \ker \Gamma(\zeta)\big) \big| .  
\end{equation}
we see that this stability  condition is equivalent to the following definition: 

\begin{defi} Given a profile $w$, the linearized equation $\eqref{linq}$ 
satisfies  the uniform spectral stability condition for high frequencies  when 
there are $c> 0$ and $\rho_1 > 0$ such that 
$D^{sc}  (\zeta) \ge c $ for all $\zeta \in \overline \RR^{d+1}_+  $ with   
$  | \zeta | \ge  \rho_1$.

\end{defi}

Note that for   $\zeta$ in bounded sets, $J_{\zeta} $ and $J_\zeta^{-1}$ are uniformly bounded 
and $D (\zeta) \approx D^{sc}(\zeta)$,  thus the condition 
$D^{sc} (\zeta ) \ne 0$ is nothing but a reformulation of the weak stability condition. 

By   Lemma~\ref{lemustab1}, the high-frequency uniform stability is equivalent to \eqref{stabcondhf}.
In section \ref{S4}, we will recall from \cite{GMWZ4} that 
\emph{the uniform spectral stability implies the high-frequency maximal estimates \eqref{maxesthf}, 
under structural assumptions on the system that are satisfied in many examples}, 
including Navier-Stokes and MHD.

\begin{rem}\label{hfcrit}
\textup{
The structural assumptions we refer to are connected with well-posedness
of the initial-value problem for the viscous equations.
For shock waves, they by themselves guarantee spectral stability
and maximal estimates \cite{GMWZ4}.
For boundary-value problems, they reduce spectral stability to
well-posedness of the frozen-coefficient boundary-value problem
at the boundary; see \cite{MZ1, GMWZ5} for further discussion.
}
\end{rem}

\subsubsection{The inviscid case}
There are analogous definitions for the linearized hyperbolic problem
\eqref{linhypeq} with boundary conditions \eqref{linhypbc}. 
 Recall that the Lopatinski determinant is defined at 
\eqref{Lopdet}.  Definition~\ref{defstab2} of weak stability is strengthened as follows. 
 
\begin{defi}
\label{defstab2u}
 The linearized equation 
 \eqref{linhypeq} \eqref{linhypbc} 
 satisfies the 
  uniform spectral stability condition holds when 
there are $c> 0$  such that 
$D (\zeta) \ge c $ for all $\zeta \in S^d_+ := S^d \cap\{ \gamma > 0 \} $.  
\end{defi}

This uniform stability condition is equivalent 
to a uniform estimate for  all $\zeta \in S^d_+$:  
\begin{equation}
\label{hypsondstab}
\forall u \in \EE^- (H_0(\underline u, \zeta)): \quad 
\big| u \big| \le C \big| \underline \Gamma_{  red}  u \big|
\end{equation}
The expected maximal estimates for solutions of 
 \eqref{linhypeq} \eqref{linhypbc} are 
\begin{equation}
\label{maxestinv}
\gamma^{\mez}  \| u \|_{L^2}  + | u (0 ) |  \le  
C \big( \gamma^{- \mez} \| f \|_{L^2} + | h | \big)
\end{equation}
with $C$ independent of $\zeta \in \RR^{d+1}_+ $.


\subsection{The Zumbrun-Serre-Rousset Theorem and the reduced low-frequency problem }

In this section, we extend the previous results of \cite{ZS} and \cite{Rou} 
which link  the low-frequency uniform stability of the viscous regularizations and the 
uniform stability of the limiting inviscid problem.  
First, we recall that the transversality of the profile is a necessary condition.

\begin{prop}
\label{propstab0}
Given a profile $w$, if the low-frequency uniform spectral stability condition is satisfied, 
then $w$ is transversal. 
\end{prop}

\begin{proof}  Lemma~\ref{lemspec} implies that 
for $\zeta \ne 0$ small enough, $\widetilde U$ is a solution of \eqref{linq5} if and only if 
${}^t( u_H, u_P)= V^{-1} (\zeta)  \widetilde U$ satisfies 
 \begin{eqnarray}
\label{eqH} \D_z u_H   & = & H(\underline u , \zeta)  u_H  + f_H ,  \\
\label{eqP} \D_z u_P   & = & P(\underline u, \zeta)  u_P + f_P ,\\
\label{eqGamma}  \Gamma_H(\zeta) u_H (0)+ 
\Gamma_P (\zeta) u_P {} (0) & :=& \widetilde \Gamma (\zeta) \widetilde U(0)  \ = \  g,    
\end{eqnarray}
where  ${}^t( f_H, f_P)= V^{-1} (\zeta)  \widetilde F$ and 
$\Gamma_H$ [resp $\Gamma_P$] denotes the restriction of 
$\widetilde \Gamma V $ to $\CC^N \times \{ 0 \}$ 
[resp. $\{0\} \times \CC^{N'}$].  
In particular, 
 \begin{equation*}
\EE^- ( G(\underline u, \zeta))  = 
 V(\zeta) \Big( \EE^-(H(\underline u , \zeta)) \oplus \EE^-(P(\underline u,  \zeta)) \Big).   
\end{equation*}
With \eqref{stabcond}, this shows that the low-frequency uniform stability 
condition holds if and only if there are $C$ and $\rho_0 > 0$ such that 
for all $\zeta \in \overline \RR^{d+1}_+ $ with 
$0 < | \zeta | \le \rho_0$
\begin{equation}
\label{stabcond2}
\begin{aligned}
\forall u_H \in & \EE^-(H(\underline u , \zeta)) , \  \  \forall   u_P \in  \EE^-(P(\underline u,  \zeta))  : \\
&\big| u_H \big| + \big| u_P \big| \le  C \big|  \Gamma_H(\zeta) u_H + 
\Gamma_P (\zeta) u_P   \big|. 
\end{aligned}
\end{equation}
In particular, 
\begin{equation}
\label{stabcondP}
 \forall   u_P \in  \EE^-(P(\underline u,  \zeta))  : \quad 
  \big| u_P \big| \le  C \big|    \Gamma_P (\zeta) u_P   \big|  .
  \end{equation}
   By Lemma~\ref{lemspec}, 
$\EE^-(P(\underline u, \zeta))$ is a smooth bundle for $\zeta$ in a neighborhood of $0$. Moreover, 
$\widetilde \Gamma(\zeta)$ and $\Gamma_P(\zeta)$ are smooth around the origin. 
This implies that $  \big| u_P \big| \le  C \big|    \Gamma_P (0) u_P   \big| $ on 
$\EE^-(P(\underline u,0))$, implying that condition i) of Definition~\ref{deftrans} is satisfied. 

Since $\dim (\EE^-(G(\zeta))  = \mathrm{rank} \widetilde \Gamma(\zeta)= N_b$, 
\eqref{stabcond2} implies that for all $h \in \CC^{N_b}$ and all
$\zeta \in \overline \RR^{d+1}_+ $ with 
$0 < | \zeta | \le \rho_0$, there is   $ \tilde U(\zeta) = V(\zeta) \big( u_H(\zeta), u_P(\zeta) \big) $  
in $\EE^-(\zeta) \subset V(\zeta) \big(\CC^N \oplus \EE^-(P(\zeta)) \big)$ 
such that $\widetilde \Gamma(\zeta) \tilde U(\zeta) = h $ and 
$| \tilde U(\zeta) )| \le c | h|$.  By compactness and continuity, letting $\zeta$ tend to zero, implies that there is   $ \tilde U  = V(0) \big( u_H , u_P \big) $  
in $  V(0) \big(\CC^N \oplus \EE^-(P(0)) \big)$ 
such that $\widetilde \Gamma(0) \tilde  U  = h $, showing that condition ii)  
of Definition~\ref{deftrans} is also satisfied. 
\end{proof}

Suppose   that the profile $w$ is transversal. Then, by i) of Definition~\ref{deftrans}
and Remark~\ref{rem0}, 
$\Gamma_P (\zeta )$ is an isomorphism 
 from $\EE^-(P(\underline u, \zeta))$  to its image  $\FF_{0, P} $ when 
 $\zeta = 0$; by continuity
 this extends  to a neighborhood  of the origin  and 
  the decomposition \eqref{dedbc} valid at $\zeta = 0$, extends smoothly 
  on a neighborhood of the origin: 
 \begin{equation}
\label{dedbc2}
\CC^{N_b}  = \FF_H  \oplus \FF_P(\zeta)   , \quad  
  \FF_P (\zeta)  := \Gamma_P (\zeta)  \EE^-(P(\underline u, \zeta)). 
\end{equation}
 Denote by $\pi_H(\zeta) $ and $\pi_P(\zeta) $ the projections associated to this splitting and 
 define the \emph{reduced boundary operator} as
  \begin{equation}
\label{redbc}
\Gamma_{ red} (\zeta) := \pi_H (\zeta) \Gamma_{H} (\zeta),  
\end{equation}
 as well as the \emph{reduced boundary value problem}
  \begin{equation}
\label{redeqH}
\D_z u_H -  H(\underline u, \zeta) u_H = f_H , \quad    \Gamma_{ red} (\zeta) u_H(0) = h.    
\end{equation}

 The \emph{reduced Evans function} is 
 \begin{equation}
\label{redEvf}
D_{red} (\zeta)  = \big|  \det \big( \EE^-(H(\underline u, \zeta) ), \ker \Gamma_{red}(\zeta) \big) \big|.
\end{equation}
\begin{defi}
 The  {reduced   uniform stability condition} is satisfied  if 
$D_{red} (\zeta) \ge c > 0 $ for all $\zeta \in \overline \RR^{d+1} \backslash\{0\}$ with 
$  | \zeta |  $ small enough. 
\end{defi}
This is equivalent to the condition 
\begin{equation}
\label{redstabcond}
\forall u \in \EE^- (H(\underline u, \zeta))
: \quad |  u  | \le C  | \Gamma_{red} (\zeta) u |,  
\end{equation}
for   $\zeta \in \overline \RR^{d+1} \backslash\{0\}$  small. 
 
 \begin{theo}
 \label{ZS1} 
 Given a profile $w$,  the linearized equation $\eqref{linq3}$ satisfies  the 
low-frequency uniform spectral stability condition if and only if 

\quad i) $w$ is transversal, 

\quad ii) the reduced problem $\eqref{redeqH}$ satisfies the 
reduced   uniform stability condition. 
 \end{theo}
 
 \begin{proof}
 We have already shown that the low-frequency uniform stability requires that 
 $w$ is transversal. Moreover, using the splitting \eqref{dedbc2}, we see that 
the uniform stability conditions  \eqref{stabcond} or \eqref{stabcond2} are  equivalent to 
 \begin{equation}
 \label{stabcond3}
 \big| u_H \big| +  \big| u_P  \big| \le C \Big( 
  \big| \Gamma_{red} u_H \big| +  \big| \Gamma_P u_P + \pi_P \Gamma_H u_H  \big| \Big) 
 \end{equation}
 for all $u_H \in \EE^-(H) $ and $u_P \in \EE^-(P)$ (to lighten notations we have 
 omitted the $\zeta$ dependance).  Since $\Gamma_P $ is surjective from 
 $\EE^-(P)$ onto $\FF_P$, for all $u_H \in \EE^-(H)$ there is 
 $u_P \in \EE^-(P)$ such that $\Gamma_P u_P  = - \pi_P \Gamma_H u_H   $ 
 and \eqref{stabcond3} implies \eqref{redstabcond}. 
 
 Conversely, if the profile is transverse,  the estimate \eqref{stabcondP} is valid 
 at $\zeta = 0$ and extend by continuity to $\zeta$ in a neighborhood of $0$. 
 With \eqref{redstabcond}, this clearly implies \eqref{stabcond3}. 
 \end{proof}

 It remains to link the reduced uniform stability condition to the 
 uniform (Lopatinski) stability condition for the hyperbolic 
 boundary value problem, that is for the problem \eqref{linhypeq} with 
 boundary conditions \eqref{linhypbc}. Note that these boundary conditions 
 are given by 
 $\underline \Gamma_{  red} = \Gamma_{red} (0)$ (see Remark~\ref{remhypbc}). 
 
 Because $H$ vanishes at $\zeta = 0$, it is natural to use 
 polar coordinates:  
 \begin{equation}
\label{polcoo}
\zeta = \rho \cz , \quad  \rho = | \zeta |, \ \ \cz \in S^d.
\end{equation}
In these coordinates 
\begin{equation}
\label{polH}
H(\underline u, \zeta) = \rho \chH (\underline u, \cz, \rho), 
\quad  \chH (\underline u, \cz, \rho) = H_0 (\underline u, \cz) + O(\rho). 
\end{equation}
Changing $z$ to $\check z =  \rho z  $,
$u(z)$ to $\check u(\check z)$ and 
$f (z)  $ to $ \rho  \check f ( \check z)$  the reduced problem \eqref{redeqH} is equivalent to 
  \begin{equation}
\label{redeqHe}
\D_{\check z}  \check u_H -  H(\underline u, \cz, \rho ) \check u_H = \check  f_H , \quad   
 \Gamma_{ red} (\zeta) \check u_H(0) = h,    
\end{equation}
which, for $\rho = 0$, is exactly the inviscid problem \eqref{linhypeq1} \eqref{linhypbc}.  
We are thus led to a \emph{nonsingular} perturbation problem. 

 Clearly, 
 for 
 $\zeta \in \overline S^d_+ := S^d \cap \{ \cg \ge 0 \}$, 
 there holds   $\EE^- (H(\underline u, \zeta) ) = \EE^-(\chH(\underline u, \cz, \rho)$
 and  $D_{red} (\zeta) =  \ccD (\cz, \rho) $  with 
 \begin{equation}
\label{red}
\ccD (\cz, \rho) = \big| \det \big( \EE^- (\chH(\underline u , \cz, \rho) , 
\ker \Gamma_{red} (\rho \cz) \big) \big| 
\end{equation}
 
 \begin{rem}
 \label{remdetlop}
 \textup{For $\cg > 0$, $H_0(\underline u, \cz)$ has no eigenvalues 
 on the imaginary axis, as a consequence of   hyperbolicity  (see Remark~\ref{remhyp}). By perturbation, this property holds true for 
 $\chH(\underline u, \cz, \rho) $ for $\rho$ small enough (depending on $\cg > 0$).
 This shows that the vector bundle $\EE^- (\chH(\underline u , \cz, \rho) $
 which was defined on  $ \overline S^d_+ \times ]0, \rho_0]$  has a  smooth
 extension to  $  \in S^+ \times [0, \rho_0]$, as well as $\ccD$. Comparing with 
 the definition of the Lopatinski determinant \eqref{Lopdet}, we see that   }
 \begin{equation}
\label{Lopdet1}
D_{Lop} (\cz)  = \ccD (\cz, 0) , \quad  \mathrm{for} \  \cg > 0. 
\end{equation}
 \end{rem}

 The next theorem, combined with Theorem~\ref{ZS1}, extends 
 Rousset's theorem \cite{Rou} (see also \cite{ZS} for shocks).
 
 \begin{theo}
 \label{ZS2}
 Given a transverse profile $w$,   if the reduced   uniform spectral stability condition
 is satisfied,  then the linearized hyperbolic problem $\eqref{linhypeq} \ \eqref{linhypbc}$ satisfies 
reduced   uniform stability condition. 

Conversely, if the linearized hyperbolic problem  is uniformly stable
and the vector bundle $\EE^- (\chH(\underline u , \cz, \rho) $ has a continous 
extension to $\overline S^d_+ \times [0, \rho_0]$, then  
 the reduced   uniform spectral stability condition
 is satisfied and the linearized problem $\eqref{linq}$ satisfies the uniform
 low-frequency stability condition. 
 \end{theo}  
 
 \begin{proof}
 The uniform estimate \eqref{redstabcond} implies that 
 \begin{equation*}
 | u | \le C | \Gamma_{red}(\zeta) u | 
\end{equation*}
for $u \in \EE^-(\chH(\underline u, \cz, \rho) $, $\cz \in \overline S^d_+$
and $\rho > 0$ small. If $\cg > 0$,  every term is continuous up to $\rho = 0$ and 
the estimate above  implies \eqref{hypsondstab},  that is
\begin{equation*}
 | u | \le C | \Gamma_{red}(0) u | 
\end{equation*}
for $u \in \EE^-(H_0(\underline u, \cz) $, $\cz \in  S^d_+$. This implies that 
the hyperbolic problem in uniformly stable. 

If $\EE^- (\chH(\underline u , \cz, \rho) $ has a continous 
extension to $\overline S^d_+ \times [0, \rho_0]$, the reduced Evans function 
is has a continuous extension to $\overline S^d_+ \times [0, \rho_0]$. 
The hyperbolic   uniform stability  and \eqref{Lopdet1} imply that 
\begin{equation*}
\ccD (\cz, \rho ) \ge c > 0
\end{equation*}
 for $\cz \in S^d_+$ and $\rho = 0$. By continuity, this extends first to $\cz \in \overline S^d_+$ and next to $\rho \in [0, \rho_1]$ for some 
 $\rho_1 > 0$.
 \end{proof}
 
 \begin{rem}
 \textup{It is proved in \cite{MZ3} that when the eigenvalues of the hyperbolic symbol  $\oA(u, \xi)$ have constant multiplicity, and more generally when
 there is a smooth K family of symmetrizers (see the definition below),  then the vector bundle  $\EE^- (\chH(\underline u , \cz, \rho) $ has a continuous extension to $\rho =0$.   The main concern of this paper is to construct K-families 
 for systems with variable multiplicity.  This is possible under suitable assumptions, 
 and therefore the two theorems above extend a result of F.Rousset \cite{Rou}.  
 However, we will also show that  the  bundle $\EE$ does not always 
 admit a continuous extension, with the result that the 
hyperbolic  problem can be uniformly stable while  the 
viscous problem is strongly unstable in the low-frequency regime.  
 This seems to be a new phenomenon.  } 
 \end{rem}
 
 Assuming transversality of $w$, Theorem~\ref{ZS1} implies that the uniform spectral stability for low frequency is equivalent to the
 spectral stability for the reduced problem. There is an analogue for maximal estimates. 
 The \emph{maximal estimates for the reduced problem} \eqref{redeqHe} read  
\begin{equation}
\label{redmaxest}
(\cg  + \rho )^{\mez}  \| \check u_H \|_{L^2}  + | \check u_H (0 ) |  \le  
C \big(( \cg + \rho )^{- \mez} \| \check f_H   \|_{L^2} + | h  | \big)
\end{equation}
with $C$ independent of $\cz  \in \overline S^{d}_+ $ and $\rho \in ]0, \rho_0]$. 
Note that for $\rho = 0$ and $\cg > 0$,  this is the maximal estimate 
for the inviscid problem. 
Scaling back to the original variables, this estimate is equivalent to 
\begin{equation}
\label{redmaxest2}
(\gamma  + | \zeta |^2  )^{\mez}  \|   u_H \|_{L^2}  + |   u_H (0 ) |  \le  
C \big(( \gamma + | \zeta |^2 )^{- \mez} \|   f_H   \|_{L^2} + | h  | \big)
\end{equation}
for the solutions of \eqref{redeqH}. 

 \begin{theo}
 \label{thmaxestred}
 Suppose that the profile $w$ is transversal.  
 Then the maximal estimates $\eqref{maxest2ss}$ are valid for low frequencies
 if and only if the maximal estimates $\eqref{redmaxest}$ for the reduced problem
 hold true. 
 \end{theo}
 
 \begin{proof}
 By Lemma~\ref{lemspec} $P(\underline u, \zeta)$ has no purely imaginary eigenvalues. 
 Thus, using symmetrizers (see e.g.  \cite{MZ1} and Section~\ref{S3} below), there holds 
 \begin{eqnarray}
 \label{estuP+}
\|  u_P^+ \|_{L^2 }   + | u_P^+(0) | \ls   \| f_P^+ \|_{L^2} ,
  \\
  \label{estuP-}
\|  u_P^- \|_{L^2 }   \ls   \| f_P^-  \|_{L^2}    + | u_P^-(0) |, 
\end{eqnarray}
 where $\pm$ denotes the smooth projections on the spaces $\EE^\pm(P(\underline u, \zeta)$. 
 
 The splitting \eqref{dedbc2} implies that the boundary condition \eqref{eqGamma} reads
 \begin{eqnarray*}
&& \pi_H g =   \Gamma_{red} u_H(0)  +   \pi_H \Gamma_P u_P^+ (0),   \\
&&\pi_P g  =    \Gamma_P  u_P^-(0) +  \pi_P \Gamma_{H} u_H(0)  +   \pi_H \Gamma_P u_P^+ (0). 
\end{eqnarray*}
Moreover  $ \Gamma_P $  is invertible on $\EE^-(P)$, hence
$ |   \Gamma_P  u_P^-(0)  | \approx | u_P^-(0)|$ and 
 \begin{eqnarray*}
&&    |  \Gamma_{red} u_H(0) |  \ls  | \pi_H g |  +  |  u_P^+ (0)|,   \\
&& |   u_P^-(0) | \ls  | \pi_P g |    +  |  u_H(0)  |  +   |  u_P^+ (0) |. 
\end{eqnarray*}
Suppose that the estimate \eqref{redmaxest2} is satisfied. Then, 
\begin{equation*}
 \vp  \| u_H \|_{L^2} + | u_H(0)| \ls  
  \vp^{-1}   \| f_H \|_{L^2} + | \pi_H g | +  | u_P^+(0)|.  
 \end{equation*}
 With \eqref{estuP+}, this implies that 
\begin{equation*}
\begin{aligned}
 \vp  \| u_H \|_{L^2} +   \| u_P^-  \|_{L^2} + & | u_H(0)| + | u_P^-(0) |  \\
 &  \ls  
  \vp^{-1}   \| f_H \|_{L^2}  +   \| f_P^- \|_{L^2} + |   g | +  | u_P^+(0)|.  
  \end{aligned}
 \end{equation*}
 Thus, with \eqref{estuP+}, we obtain that 
\begin{equation*}
 \vp  \| u_H \|_{L^2} +   \| u_P \|_{L^2} +  | u_H(0)| + | u_P (0) |  \\
   \ls  
  \vp^{-1}   \| f_H \|_{L^2}  +   \| f_P  \|_{L^2} + |   g | .  
 \end{equation*}
Because $V(\underline u, 0)$ has the special form \eqref{Vtri}, 
$\widetilde U =  V (u_H, u_P) = (\widetilde U^1, \widetilde U^2) $  satisfies 
\begin{equation*}
 \widetilde U^1 = O(1) u_H + O(1) u_P, \quad 
 \widetilde U^1 = O(| \zeta | ) u_H + O(1) u_P
\end{equation*}
 Therefore, the solutions of \eqref{linq5} satisfy 
 \begin{equation*}
 \vp  \| \widetilde U^1  \|_{L^2} +   \| \widetilde U^2 \|_{L^2} +  | \widetilde U  (0) |  \\
   \ls  
  \vp^{-1}   \| \widetilde F    \|_{L^2}   + |   g | .  
 \end{equation*}
 that  is the maximal estimate \eqref{maxest2ss}. 
 
 \medbreak
 
 Conversely, assume that the maximal estimate \eqref{maxest2ss} is satisfied. 
Suppose that  $u_H$ is a solution of \eqref{eqH}. 
By transversality, $\Gamma_P$ is surjective from $\EE^-(P, \zeta)$ to its image 
$\FF_P (\zeta)$ and there exists there is $u_P(0) $ in 
$\EE^-(P, \zeta)$  such that 
\begin{equation}
\label{machin}
\Gamma_P u_P(0)  = - \pi_P \Gamma_H u_H (0) \in \FF_P(\zeta). 
\end{equation}
Consider $u_P = e^{ z P}  u_P(0)$ which is well defined and rapidly decaying 
at infinity since $u_P(0) \in \EE^-(P, \zeta)$. It is a solution of 
\eqref{eqP} with $f_P = 0$. 
Then   
 $\widetilde U = V (u_H, u_P) $  is a solution of \eqref{linq5}  with 
 $\widetilde F = V (f_H, 0)$.  
 Thus  $(u_H, u_P) = V^{-1} \widetilde U$ and there holds
 $$
 \| u_H \|_{L^2} \ls \| \widetilde U \|_{L^2} , \quad   | u_H(0)| \ls | \widetilde U(0) |, 
 \quad  \| \widetilde F \|_{L^2} \ls \| f_H \|_{L^2}  
 $$
  and, by \eqref{machin}, $\widetilde \Gamma \widetilde U(0)  = 
  \Gamma_H u_H(0) + \Gamma_P u_P(0) = \Gamma_{red} u_H(0)$.  
  Thus the estimate \eqref{maxest2ss} immediately implies \eqref{redmaxest2}. 
 \end{proof}

%
%
%

\section{Low frequency analysis:  the main results}
\label{S3}

This section is mainly devoted to the study of the reduced equation \eqref{redeqHe}, 
which is a nonsingular perturbation of the inviscid problem \eqref{linhypeq}. 
 Our goal is to perform an analysis without assuming constant multiplicity of 
 eigenvalues, thus allowing for examples such as MHD. 
 The inviscid case is considered in \cite{MZ2}, and we want to extend the results 
 to  small viscous perturbations.

\subsection{Symmetrizers}

   Consider the constant-coefficient linear first order system \eqref{linq5}. 
 For clarity, we drop the tildes and reserve the notation $u, U... $ for the unknowns
 and call $p \in \cU$ the parameter called $\underline u$ in this equation, which now 
 reads
 \begin{equation}
 \label{nlinq5}
 \D_z U = G(p, \zeta) U + F, \quad   \Gamma (p, \zeta) U (0) = g. 
 \end{equation} 
To prove energy estimates  for the solutions of this equation,   the main step  
is to construct \emph{symmetrizers}. 
They are   self adjoint 
matrices $\Sigma (p, \zeta)$ such that 
\begin{equation}
\label{sym1}
\re \big(\Sigma (p,  \zeta) G(p, \zeta) \big) > 0. 
\end{equation}
The symmetrizer is adapted to the boundary conditions 
and provides maximal estimates for the traces when 
\begin{equation}
\label{sym2}
\Sigma (p, \zeta) > 0  \quad \mathrm{on} \ \ker \Gamma(p, \zeta). 
\end{equation}

The construction of such symmetrizers is in two steps:  
 first, one  constructs
a family of symmetrizers $\Sigma^\kappa$, which is independent 
of the boundary conditions; second one uses the uniform Lopatinski or Evans 
condition, to prove that if $\kappa $ is large enough then the simmetrizer 
is adapted to the boundary condition. 

More precisely, one considers  a splitting
 \begin{equation}
\label{split}
\CC^{N + N'} = \EE^- (p, \zeta)  \oplus \EE^{+}(p, \zeta) 
\end{equation}
where $\EE^- (p, \zeta) $ is the negative invariant space of $G(p, \zeta)$ as above
while $\EE^+ (p, \zeta) $ can be chosen arbitrarily so that the 
splitting \eqref{split} holds .
Denoting  by $ \Pi^\pm(p, \zeta) $ the   projectors
associate to this splitting, the family of symmetrizers $\Sigma^\kappa $ is 
searched so that 
\begin{equation}
\label{sym3}
\Sigma ^\kappa \ge  m(\kappa ) (\Pi^+ )^*  \Pi ^+  -  (\Pi^- )^*  \Pi ^- 
\end{equation}
where $m(\kappa) \to + \infty$ as $\kappa \to + \infty$. 

Since $\EE^- = \ker \Pi^+$, the stability condition \eqref{weakstab} which reads
\begin{equation}
\label{nweakstab}
\ker \Gamma(p, \zeta) \cap \EE^{-}(p, \zeta) = \{ 0 \} 
\end{equation}
 is also 
equivalent to an estimate
\begin{equation}
\label{st1}
| \Pi^- u |^2 \le C (  | \Gamma u |^2 + | \Pi^+ u |^2 ) . 
\end{equation}
Therefore,   if  the family  $\Sigma^\kappa$ satisfies \eqref{sym3},  
then  for  $\kappa$   large enough, 
there holds  
\begin{equation}
\label{sym4}
\Sigma^\kappa \ge c \Id  -  C' \Gamma^* \Gamma , \quad  c > 0 
\end{equation}
and therefore $\Sigma^\kappa$  is adapted to the boundary condition $\Gamma$. 

If $\re \Sigma^\kappa G \ge \delta_\kappa \Id$, then multiplying the equation by 
$\Sigma^\kappa$ and integrating by parts yields the estimate
\begin{equation}
\label{est0}
 \delta_\kappa \| U \|^2 _{L^2} + c | U(0 )|^2   \le \frac{1}{\delta_\kappa} \| F \|^2_{L^2} 
 +  C' | g |^2. 
\end{equation}

This is the sketch of the general argument. To obtain usable 
  estimates,   uniform  versions of \eqref{sym3} \eqref{st1} are needed  as well as 
more precise versions of   \eqref{sym1} (see below). 
Note that in this approach, the construction os symmetrizers is completely independent
of the boundary conditions, and in particular of the validity of the stability conditions.
In this paper, we concentrate on the construction of families of symmetrizers
which satisfy \eqref{sym3}. They are called K-families in \cite{MZ2}.  


\subsection{Main results} 
  The construction of symmetrizers for middle frequencies,
   is performed in \cite{MZ1}.    
 By Lemma~\ref{lemspec}, the matrix $G(p, \zeta) $ has no eigenvalues on the 
 imaginary axis when $\zeta \in \overline \RR^{d+1}\backslash\{0\}$. 
 Therefore,

 \begin{lem}
 \label{lemsymMF} 
 For all $\underline \zeta  \in \overline \RR^{d+1}\backslash\{0\}$, 
 there is a neighborhood of $(\up, \underline \zeta)$ in 
 $\cU \times \RR^{d+1}$ such that for $(p, \zeta) $ in this neighborhood  
 there is a smooth 
 splitting 
 \begin{equation}
\label{splitG}
\CC^{N'} = \EE ^- (p, \zeta) \oplus \EE^{+}(p, \zeta). 
\end{equation}
where $\EE^\pm (p, \zeta)$ denote the invariant space of $G(p, \zeta)$ 
associated to the spectrum in $\{ \pm \re \mu > 0 \}$. 
Denoting  by $ \Pi^\pm(p, \zeta) $ the smooth spectral projectors
associate to this splitting, there is a  smooth family $\Sigma^\kappa (p, \zeta)$ 
of self adjoint matrices such that for all $(p, \zeta)$ in the given neighborhood 
and all $\kappa \ge 1$: 
 \begin{equation}
\label{symMF}
\begin{aligned}
i ) &  \quad  \re  \Sigma^\kappa   G  \ge 0, 
\\
ii)&  \quad  \re \Sigma ^\kappa \ge \kappa (\Pi ^+ )^*  \Pi ^+  - 
( \Pi ^-)^*  \Pi ^-. 
\end{aligned} 
\end{equation}

 \end{lem}

 \begin{cor}
 \label{lemestMF}
 If the weak spectral stability condition is satisfied, then for 
  all $\underline \zeta  \in \overline \RR^{d+1}\backslash\{0\}$, 
   there are a constant $C$ and a neighborhood of $(\up, \underline \zeta)$ in 
 $\cU \times \RR^{d+1}$ such that for $(p, \zeta) $ in this neighborhood  
 the solutions of $\eqref{nlinq5}$ satisfy
 \begin{equation}
 \label{nestMF}
 \| U \|_{L^2} + | U(0 ) |  \le  C \big( \| F \|_{L^2}  + | g | \big). 
 \end{equation}
 
 \end{cor} 
 
 \bigbreak 
We now  concentrate on   low frequencies.  
 By Lemma~\ref{lemspec}, the matrix $G(p, \zeta)$ is locally smoothly conjugated 
 to a block diagonal matrix \eqref{block1} with diagonal blocks 
with $H(p, \zeta)$ of dimension $N \times N$ and 
 $P(p, \zeta)$ of dimension $N' \times N'$. The system \eqref{nlinq5} is 
 therefore equivalent to the equations \eqref{eqH} \eqref{eqP} coupled by the 
 boundary conditions \eqref{eqGamma}.

In the block diagonal reduction \eqref{block1}, we construct symmetrizers 
\begin{equation}
\label{blockS}
\Sigma^\kappa =  
\left(\begin{array}{cc }
\Sigma^\kappa_H   & 0  \\ 0 & \Sigma^\kappa_P    
\end{array}\right)  
\end{equation}
such that the properties \eqref{sym1} ande \eqref{sym3} are satisfied  for 
each block independently.

 The construction of symmetrizers for the elliptic block $P$ is standard
and identical to the construction for middle 
frequencies, since $P(\up, 0)$ has no eigenvalues on the imaginary axis.  
 Denote by $\EE_P^{\pm} (p, \zeta) $ the subspaces of $\CC^{N'}$, invariant for 
 $P(p, \zeta)$, associated to the spectrum in $\{ \pm \re \mu > 0\}$. 
 Thus, for $(p, \zeta)$ in a neighborhood of $(\up, 0)$, there is a smooth 
 splitting 
 \begin{equation}
\label{splitP}
\CC^{N'} = \EE_P^- \oplus \EE_P^{+}. 
\end{equation}
Denote by $ \Pi_P^\pm(p, \zeta) $ the smooth spectral projectors
associate to this splitting.

\begin{prop}
There is a smooth family of self adjoint matrices $\Sigma_P^\kappa $ on a neighborhood 
of $(\up, 0)$ such that 
 \begin{equation}
\label{symP}
\begin{aligned}
i ) &  \quad  \re  \Sigma_P^\kappa   P  > 0 , 
\\
ii)&  \quad  \re \Sigma_P^\kappa \ge \kappa (\Pi_P^+ )^*  \Pi_P^+  - 
( \Pi_P^-)^*  \Pi_P^-
\end{aligned} 
\end{equation}
 
 \end{prop}

 This implies the estimates \eqref{estuP+} \eqref{estuP-} which where used in the 
 previous section.

 \bigbreak 
 To analyze $H$, we 
 use polar coordinates  for $\zeta = \rho \cz $ as in \eqref{polcoo} 
so that  
\begin{equation}
\label{defcH}
H(p, \zeta) = \rho \chH (p, \cz, \rho) , \quad 
\chH (p, \cz, \rho ) = H_0(p, \cz)   + O(\rho). 
\end{equation}
By Lemma \ref{lemspec},  for $\zeta \in \RR^{d+1}_+ \backslash\{0\}$, 
 $\chH$ has no eigenvalue on the imaginary 
axis, hence the number $N^-$  of eigenvalues of 
 $\chH $ in $\{ \re \mu < 0 \} $ is constant.

We fix a point $\underline \cz \in \overline S^d_+ $, that is 
$\underline \cz = (\underline \ct, \underline \ch, \underline \cg) $ in the unit sphere
with $\underline \cg \ge 0$. The goal is to construct smooth symmetrizers 
for $\chH$, for $(p,  \cz, \rho)$ close to $(\up, \underline \cz, 0)$. 
For convenience we introduce the following terminology. 

\begin{defi} 
\label{def215}
A   smooth  symmetrizer  for  $\chH $ on 
a neighborhood  $\omega$ of $ (\up, \underline \cz, 0)$
is a smooth self adjoint  matrix $\ccSi^H  (p, \cz, \rho)$ 
such that 
 \begin{equation}
\label{global41}
  \re  \ccSi^H   \chH   =  \sum V_k^*  \Sigma_k V_k ,  
\end{equation}
where the $V_k$ and $\Sigma_k$ are smooth matrices on $\omega$ 
of appropriate dimension so that the products make sense, 
satisfying 

\quad i)  \ $\sum V_k^* V_k$ is definite positive, 

\quad ii) either $\Sigma_k$ is definite positive or 
$\Sigma_k = \gamma \Sigma_{k, 1} + \rho \Sigma_{k, 2}$ with 
$\Sigma_{k, 1}$ and $\Sigma_{k, 2}$ definite positive. 
 
 \end{defi}
 
 \begin{defi} 
 \label{def216}
 A  family of 
   smooth symmetrizers  $\Sigma^\kappa $ on neigborhoods $\omega^\kappa$ 
 of $ (\up, \underline \cz, 0)$ is called a K-family of symmetrizers for $\chH$ if  there   are  a decomposition 
 \begin{equation}
 \label{dec}   \CC^N  = \underline \EE_H^-  \oplus \underline \EE_H^+
 \end{equation} 
 with $\dim \underline \EE^- = N_{ -}$  and 
 $m(\kappa) \to + \infty $ as $\kappa \to + \infty$  such that  for all 
 $\kappa$ 
\begin{equation}
\label{232}
\Sigma^\kappa  (\up, \underline \cz, 0)  \ge  m(\kappa)  \Pi_+^*   \Pi_ +  - 
 \Pi_-^*   \Pi_-. 
\end{equation}
where $  \Pi_\pm   $ 
  are the projectors associated to the splitting $\eqref{dec}$.

\end{defi}

\begin{rem}
\textup{ Recall from \cite{MZ3}  that if there is K-family of symmetrizers, then 
$\underline \EE^-$ is the limit of the negative spaces 
$\EE^-(p, \cz, \rho)$ as $(p, \cz, \rho)$ tends to $(\up, \underline \cz, 0) $
with  $\rho > 0$. Thus $\underline \EE^-$ is uniquely 
determined. On the other hand, $\underline \EE^+$ is arbitrary, provided that 
the the splitting \eqref{dec} holds:  if \eqref{232} holds for some choice of 
$\underline \EE^+$, then it also holds for another choice for a multiple 
of $\Sigma^\kappa$ with some other function $m(\kappa)$. }
\label{rem312}
\end{rem}

We can now state the main result of this paper, which  extends  \cite{MZ2}.

 \begin{theo} 
 \label{Kfam}
  Suppose that the Assumptions of Section $2.1$   are satisfied. Assume further that one of the following 
  two condition is satisfied:
  
  \quad i)  
all the   real characteristic 
  roots  $(\up, \tau, \xi)$ with $| \xi | = 1$   satisfy the 
  condition (BS) of Definition~$\ref{defBS}$. 
  
  \quad ii) the system is symmetric dissipative in the sense of Definition~$\ref{defsymm}$   and the real characteristic 
  roots  $(\up, \tau, \xi)$ with $| \xi | = 1$   are either totally nonglancing in the sense of Definition~$\ref{def23}$ or satisfy the 
  condition (BS) of Definition~$\ref{defBS}$. 
   
    Then, for all   $\underline \cz \in \overline S^d_+$, there exists   K-families of smooth symmetrizers for $\chH(p, \zeta, \rho)$ near $(\up, \underline \cz, 0)$.

    \end{theo}

 The  condition $(BS)$ ensures that a suitable 
 generalized block structure condition, which makes the 
 construction of symmetrizers given in   \cite{MZ1} work, is satisfied. For hyperbolic problems, it is shown in 
  \cite{MZ2} that the block structure condition is satisfied if the system is smoothly 
  diagonalizable. 
   In the  viscous case, things are 
  more subtle and the generalized block structure condition is discussed in details in Section 3. 
We just point out  here the  following example. 
  
  \begin{theo} 
  \label{GBS}  
  If $(\up, \underline \ct , \underline \cx)$ is a semi-simple 
    characteristic root of constant multiplicity, then 
  the   condition (BS)  of Definition~$\ref{defBS}$ is satisfied at that point.   
     \end{theo} 
 
 Together with Theorem \ref{Kfam}, this implies Theorem~\ref{th11}.  
 Finally, we quote that the existence of K-families implies the validity 
 of the maximal estimates when the boundary conditions satisfy 
 the uniform spectral stability conditions.

 \begin{theo}
 Suppose that there exists a  K-families  of symmetrizers for 
 $\chH$ near $(\up, \underline \cz, 0)$ and suppose that the boundary conditions 
 are such that the uniform spectral 
 stability condition  is satisfied for low frequencies. 
 Then the uniform stability estimates $\eqref{maxest2ss}$ are satisfied. 
 
 Similarly, if the reduced boundary conditions satisfy the reduced uniform stability condition
 then the uniform estimates $\eqref{redmaxest}$ and $\eqref{redmaxest2}$ hold true. 
 \end{theo} 

   
    \subsection{Block reductions }
 
  The advantage of the notion of K-families is that it is independent 
  of the boundary conditions. Therefore, their construction 
  depend only on an analysis of $\chH$. In particular, 
  we can use spectral block decompositions of $\chH$. 
   
 Fix $\underline \cz \in \overline S^d_+$. 
Consider the \emph{distinct} eigenvalues $\underline \mu_k$ of $H_0(\up, \underline \cz)$. 
For $(p, \cz, \rho)$ in a neighborhood of  $(\up, \underline \cz, 0)$, there is a smooth
block reduction 
\begin{equation}
\label{block2}
V^{-1} \chH   V = \diag  ( \chH_k) 
\end{equation}
where the $H_k$ have their spectrum in small discs centered at $\underline \mu_k$
that are pairwise disjoints. Equivalently, there is a smooth decomposition 
\begin{equation}
\label{b2b}
\CC^N = \bigoplus_k \EE_k( p, \cz, \rho) 
\end{equation} 
in invariant spaces for $\chH (p, \cz, \rho)$ and $\chH_k$ is  the restriction 
of $\chH$ to $\EE_k$. 
We denote by $N_k$ the dimension of $\EE_k$, that is the size of $\chH_k$.

The K-families of symmetrizers are constructed   for each block 
$\chH_k$ separately.  
If $\Sigma_k^\kappa$ is a K-family for $\chH_k$, it is clear that 
$\Sigma^\kappa = V^* \diag (\Sigma^\kappa_k) V $ has the form 
\eqref{global41} and is a K-family for $\chH$.

  When the mode is \emph{elliptic}, that is when $\re \underline \mu_k \ne 0$, 
 the construction of symmetrizers is easy (see e.g. \cite{Kr, CP, MZ1}). 
 
 \begin{prop}
 \label{prop219}
 Suppose that $\underline \mu_k$ is an eigenvalue of $H_0(\up, \underline \cz)$ 
 with $\re \underline \mu_k  \ne 0$. 
Then is a smooth family of self adjoint matrices $\Sigma^\kappa_k $ on a  neighborhood 
of $(\up, \underline \cz, 0)$ such that 
 \begin{equation}
\label{symE}
\begin{aligned}
i ) &  \quad  \re  ( \Sigma^\kappa_k   \chH_k )   > 0 , 
\\
ii)&  \quad  \re  \Sigma^\kappa_k \ge   \kappa  \Id  \quad if  \ \ 
\re \underline \mu_k > 0, 
\\ & \quad 
\re  \Sigma^\kappa_k \ge    - 
\Id  \quad if \ \ 
\re \underline \mu_k  < 0. 
\end{aligned} 
\end{equation}
 
 \end{prop}

 Therefore we now  restrict our attention to  a nonelliptic mode: 
 \begin{equation}
\label{mur}
\underline \mu_k = i \underline \cx_{d} , \quad   \underline \cx_{d} \in \RR. 
\end{equation}
By definition of $H_0$,  this implies that 
$ - \underline \ct  + i \underline \cg $  is an eigenvalue $\underline \lambda$
of $A(\up, \underline \cx)$ with $\underline \cx = (\underline \ch, \underline \cx_d)$. 
In particular,  by hyperbolicity, this  can only happen when  $\underline \cg = 0$. 
By Lemma~\ref{lemspec}, $\chH_k$ has no eigenvalues on the imaginary axis
when $\rho > 0$, thus the number of eigenvalues in 
$\{ \re \mu < 0 \}$ is constant. We call it $N^-_k$. 
The next definition reformulates Definitions \ref{def215} and \ref{def216} for nonelliptic blocks
$\chH_k$.

\begin{defi} 
\label{def316}
A   smooth  symmetrizer  for  a nonelliptic block $\chH_k $ on 
a neighborhood  $\omega$ of $ (\up, \underline \cz, 0)$
is a smooth self adjoint matrix $\Sigma (p, \cz, \rho)$ 
such that, for some  $C$, $c>0$, there holds for 
all $(p, \cz, \rho) \in \omega $, 
 \begin{equation}
\label{global41b}
 \re  \Sigma  \chH_k   = \cg \Sigma_1 + \rho \Sigma_2 ,  
\end{equation}
with $\Sigma_1(\up, \underline \cz, 0)$ and $\Sigma_2  (\up, \underline \cz, 0)$
definite positive. 

A  family of 
   smooth symmetrizers  $\Sigma^\kappa_k $ on neigborhoods $\omega^\kappa$ 
 of $ (\up, \underline \cz, 0)$ is called a K-family of symmetrizers for $\chH_k$ if  there   are  a decomposition 
 \begin{equation}
 \label{deck}   \EE_k (\up, \underline \cz, 0) = \underline \EE^-_k  \oplus \underline \EE_k^+
 \end{equation} 
 with $\dim \underline \EE_k^-$ equal to $N_{k}^-$  and 
 $m(\kappa) \to + \infty $ as $\kappa \to + \infty$  such that  for all 
 $\kappa$ 
\begin{equation}
\label{243y}
\Sigma^\kappa_k  (\up, \underline \cz, 0)  \ge  m(\kappa)  (\Pi_k^+)^*   \Pi_k^+  - 
 (\Pi_k^-)^*   \Pi_k^-. 
\end{equation}
where $  \Pi_k^\pm   $ 
  are the projectors associated to the splitting $\eqref{deck}$.

\end{defi}

Given $\underline \cz = (\underline \ct, \underline \ch, 0) \in \overline S^{d}_+$ and a nonelliptic mode 
$\underline \mu_k  = i \underline \cx_d$, $-\underline \ct,$ is an eigenvalue of $A(\up, \underline \cx)$
with $\underline \cx = (\underline \ch, \underline \cx_d)$.  
Therefore,  Theorem~\ref{Kfam}, is an immediate corollary of Proposition~\ref{prop219} and the 
following two theorems.  

\begin{theo}
\label{th221}
Suppose that the system is symmetric dispersive in the 
sense of Definition~$\ref{defsymm}$; suppose in addition that $ (\up, \underline \ct, \underline \cx) $ is a totally 
incoming or outgoing characteristic root in the sense of Definition~$\ref{def23}$. Then there are K-families of symmetrizers
for the associated block $\chH_k$, with 
$\underline \EE^-_k = \{ 0 \}$ in the outgoing case and 
$\underline \EE^-_k = \CC^{N_k}$ in the incoming case. 
\end{theo}

\begin{theo}
\label{th222}
If  $ (\up, \underline \ct, \underline \cx) $ is a 
characteristic root 
 which satisfies the generalized block structure condition of Definition~$\ref{BS2}$.   Then there are K-families of symmetrizers
for the associated block $\chH_k$. 
\end{theo}

%
%

\section{The generalized block structure condition} 

\subsection{Hyperbolic multiple roots} 
 
 We first 
recall from \cite{MZ2} several notations and definitions concerning the 
characteristic roots of the hyperbolic part $L$. For simplicity, 
we suppose, as we may, that the coefficient of $\D_t$ is 
$A_0 =\Id$, so that, with notations \eqref{y2},  $L = \overline L$. 
The characteristic determinant is denoted by  
\begin{equation}
\label{Delta}
\Delta (p, \tau, \xi ) :=  \det ( \tau \Id + A(p, \xi))   . 
\end{equation}

\begin{defi}
\label{def22}
Consider a root $(\up, \utau, \uxi)$ of $\Delta (\up,  \utau,  \uxi)) = 0$, of 
algebraic multiplicity $m$ in $\tau$. 

i)  $(\up, \utau, \uxi)$ is algebraically regular, if on a neighborhood $\omega$ of 
$(\up, \uxi)$ there are $m$ smooth real functions  
$\lambda_j(p, \xi)$, analytic in $\xi$,  such that
$\lambda_j(\up, \uxi) = - \utau$ and  for 
$(p, \xi) \in \omega$: 
\begin{equation}
\label{22}
\Delta (p, \tau,  \xi)  = e(p, \tau, \xi) \prod_{j=1}^m \big(\tau + \lambda_j(p, \xi) \big)
\end{equation}
where $e$ is a polynomial in $\tau$ with smooth coefficients such that 
$e(\up, \utau, \uxi) \ne 0$. 

ii) $(\up, \utau, \uxi)$ is geometrically regular if in addition 
there are $m$ smooth functions $e_j(p, \xi)$ on $\omega$ with values in 
$\CC^N$,  analytic in $\xi$, such that 
\begin{equation}
\label{23}
A(p, \xi) e_j(p, \xi) = \lambda_j(p, \xi) e_j(p, \xi), 
\end{equation}
and the $e_1, \ldots, e_m$ are linearly independent. 

iii)   $(\up, \utau, \uxi)$ is semi-simple with constant multiplicity if 
all the $\lambda_j$'s are equal. 
\end{defi}
 
Case $iii)$     occurs when $\lambda(p, \xi)$ is a continuous semi-simple eigenvalue of $A(p, \xi)$ 
with constant multiplicity near $(\up, \uxi)$,  such $\utau + \lambda(\up, \uxi) = 0$. This implies that 
$\lambda$ is smooth and analytic in $\xi$ as well as the eigenspace $\ker (A  - \lambda )$. 
In this case,  one can choose  for  $\{ e_j \}$ any smooth  basis of 
of this  eigenspace. 

If all the roots at  $(\up, \uxi)$ are geometrically regular, then, locally near $(\up, \uxi)$, $A(p, \xi)$ 
is smoothly diagonalizable, meaning that it has  a smooth basis of eigenvectors.

\begin{exam}
\textup{For the inviscid MHD, the multiple eigenvalues are algebraically regular, 
but some are not geometrically regular (see \cite{MZ2} and Section~\ref{SecMHD} below).}
\end{exam}

\medbreak
The second notion which plays an important role 
in the analysis of hyperbolic boundary value problems is the notion of 
\emph{ glancing modes }.  
Recall from \cite{MZ2} the following definition. 
If  $\utau $ is a root of multiplicity $m$ of the polynomial   
 $\Delta(\up, \cdot,   \uxi) $, then by hyperbolicity, the  Taylor expansion of 
 $\Delta$ at $(\up, \utau, \uxi) $ at the order $m-1$ vanishes so that 
 \begin{equation}
 \label{Deltatan}
 \Delta (\up, \utau+ \tau , \uxi+ \xi) = \underline \Delta_m(\tau, \xi) + O(|\tau, \xi|^{m+1})
 \end{equation}
and $\underline \Delta_m $ is homogeneous of degree $m$. Moreover, 
$\underline \Delta_m $ is 
hyperbolic in the time direction. Indeed, any direction of hyperbolicity 
for $\Delta (\up, \cdot)$ is a direction of hyperbolicity for $\underline \Delta_m $. 
Denote by $\underline \Gamma_+$ the open convex cone of hyperbolic 
directions fot $\underline \Delta_m $ which contains $dt$.

\begin{defi}
\label{def23}
The  root $(\up, \utau, \uxi)$ of $\Delta $, of 
 multiplicity $m$, is said nonglancing when the boundary is noncharacteristic for $\underline \Delta$. 
 
 It is totally incoming [resp. outgoing] when  the inward [resp. outward] conormal to the boundary belongs to $\underline \Gamma_+$. 
 It is totally nonglancing if is either totally incoming or totally outgoing. 
 \end{defi}

\begin{exam} \textup{This definition agrees with the usual one for simple roots, given by 
$\tau + \lambda(p, \xi ) = 0$. In this case $\D_t + \na_\xi \lambda\cdot \D_x $ is 
is the Hamiltonian transport   field for the propagation of singularities or oscillations and 
the glancing condition  $\D_{\xi_d} \lambda = 0$ precisely means that the field is tangent 
to the boundary. More generally, 
if the  root $(\up, \utau, \uxi)$ of $\Delta $ is algebraically regular, then, with notations 
as in \eqref{22}
  }
\begin{equation}
\label{gl1}
\underline \Delta_m (\tau, \xi) = e(\up, \utau, \uxi) 
 \prod_{j=1}^m \big(\tau + \xi \cdot \na_\xi \lambda_j(\up, \uxi) \big)
\end{equation}
\textup{The mode is nonglancing if none of the tangential speed 
$\D_{\xi_d} \lambda_j(\up, \uxi)$ vanish. It is totally incoming [resp.  outgoing]
if they all are positive [resp. negative]. In particular, in the  constant multiplicity case, 
all the $ \lambda_j$ are equal and they are all glancing, incoming or outgoing at the same time.  }
\end{exam}

In the study of boundary value problems, the dichotomy incoming vs outgoing plays a crucial 
role: for instance, for transport equations one boundary condition is needed in the 
first case and none in the second.  Using symmetrizers to prove energy estimates,  
 they are constructed in opposite ways.  The  general  Kreiss construction  
 also reflects this dichotomy.  Introduce the following definition: 
 
 \begin{defi} 
 \label{defiIO}
 Suppose that $(\up, \utau, \uxi)$  is an algebraically regular root of $\Delta$. With notations
 as in $\eqref{22}$, denote by $\nu_j$  the order of 
 $\underline \xi_d$ is a root of order of the equation 
$\underline \tau + \lambda_j( \up, \underline  \xi_1, \ldots, \underline \xi_{d-1}, \cdot ) = 0$, that is  the positive integer 
such that 
 \begin{equation}
\label{multig}
\D^a_{\xi_d} \lambda_j (\up, \underline \xi)  = 0  \quad \mathrm{for} \ a < \nu_j \quad 
\mathrm{and} \quad   \beta_j := \frac{1}{\nu_j !} \D^{\nu_j}_{\xi_d} \lambda_j (\up, \underline \xi)  \ne 0.    
\end{equation}
 
 We say that $\lambda_j$ is of type
 $I$ when either $\nu_j$ is even or $\nu_j$ is odd and $\beta_j > 0$. 
 It is of type $O$ when $\nu_j$ is odd and $\beta_j < 0$. 
 
 We denote by $J_O$ [resp. $J_I$] the set of indices $j$ of the corresponding type.

 \end{defi}

 \begin{rem}
 \label{rem38}
 \textup{When   $(\up, \underline \ct, \underline \cx  ) $ is non glancing, then 
 the  all the $\nu_j$ are equal to $1$, and being of type $I$ [resp. type $O$] means  to be incoming 
 [resp. outgoing]. They are all of the same type exactly when the mode is totally nonglancing. 
 }
 \end{rem}

\begin{rem}
\textup{The details of the construction of Kreiss' symmetrizers depend strongly 
on being of type $I$ or $O$, see \cite{Kr, CP, Met3} and Section~5. There are no reason 
other than technical why even roots are of type $I$ rather than $O$. }
\end{rem}

\subsection{The decoupling condition}

The spectral properties of $A(\xi)$ are modified by the perturbation $B$. 
 In particular, since the construction of symmetrizers depends deeply on the property 
 of being incoming/outgoing, it is very important that the perturbation 
 respects the decoupling between the different type of modes. 
 
 \begin{defi} 
 \label{defdec}
Suppose that $(\up, \underline \tau, \underline \xi  ) $ is a geometrically regular root of $\Delta$
of order $m$. 
Consider a basis  $\{ e_j\}$ as in  \eqref{23} and  dual  left eigenvectors 
$\underline \ell_j$  such that
\begin{equation}
\underline \ell_j \big( \utau \Id + A(\up, \uxi) \big) = 0, \quad  
\underline \ell_j \cdot e_{j'} (\up,  \uxi) = \delta_{j, j'}. 
\label{dual1}
\end{equation}
Consider the 
 and the $m \times m$ matrix with entries 
 \begin{equation}
\label{bsharp}
B^\sharp_{j, j'}  = \underline \ell_j B (\up, \uxi) e_{j'} (\up, \uxi) 
\end{equation}

i) We say that the decoupling condition is satisfied if 
\begin{equation}
\label{331x}
 B^\sharp_{j, j'} = 0  \quad when \quad  (  j ,  j') \in (J_O \times J_I ) \cup  ( J_I \times J_O)
\end{equation}
where $J_0$ and $J_I$ are introduced in Definition~$\ref{defiIO}$. 

ii) We say that the basis $\{ e_j\}$ is adapted to $B$ if
\begin{equation}
\label{334x}
\re B^\sharp > 0.
\end{equation} 

\end{defi}

\begin{defi}
\label{defBS}
We say that the root  $(\up, \underline \tau, \underline \xi  ) $ of $\Delta$ satisfies the condition 
(BS) if is 
is   geometrically regular root, satisfies the decoupling condition $\eqref{331x}$ and there is 
an eigenbasis  basis $\{e_j\}$  adapted to $B$. 
\end{defi}

We give several examples and counterexamples. 
The next result  rephrases Theorem~\ref{GBS}. 

\begin{prop} {(Constant multiplicity) }
\label{prop316}
  Suppose that    $(\up, \underline \tau, \underline \xi)$ is a semi-simple 
characteristic root with constant multiplicity of $\Delta$.   Then the   condition (BS) 
is satisfied. 
  
\end{prop}

\begin{proof}
  For semi-simple 
characteristic root $\lambda$ with constant multiplicity either $J_O$ or $J_I$ is empty so that  the decoupling condition  \eqref{331x} is trivially satisfeid. Moreover, it is proved in \cite{MZ1} that 
(H1) implies that the spectrum of $B^\sharp$ is located in $\{\re z  > 0 \}$. 
Thus there is a basis $ \{\underline e_j \}$ in  $\ker(A(\up, \uxi)  + \utau \Id ) $
such that $\re B^\sharp$ is definite positive.
Next, since any smooth basis  $\{e_j \}$  in $\ker(A - \lambda) $ satisfies 
\eqref{23}, one can choose it 
such that $e_j (\up, \uxi) = \underline e_j$.   
\end{proof}

\begin{prop}
\label{prop317} (Artificial viscosity) 
   Suppose that  $(\up, \underline \tau, \underline \xi)$  is  geometrically regular  for $iA + B$ in the sense 
   that there  
 are $m$ smooth functions $\lambda_j( p, \xi, \rho)$ 
 and $m$ linearly independent smooth vectors 
 $e_j(p, \xi, \rho)$ on a neighborhood of  $(\up, \uxi, \rho)$,   analytic in $\xi$, such that 
$\lambda_j( \up, \uxi, 0) = - \underline \tau $ for all $j$ and 
\begin{equation}
\label{23b}
\big ( i A(p, \xi) + \rho B(p, \xi) \big)  e_j(p, \xi, \rho )   =  i \lambda_j(p, \xi, \rho ) e_j(p, \xi, \rho). 
\end{equation}
Then, the  decoupling condition is satisfied and the basis
$\{ e_j {}_{| \rho =0 } \}$ is adapted to $B$. 

 \end{prop} 

\begin{proof}
Alternately, differentiating \eqref{23} with respect to $\rho$ 
and multiplying on the left by $\underline \ell_{j'}$, implies that $B^\sharp_{j', j} = 0$ when 
$j \ne j'$. Moreover, (H1) implies that $B^\sharp_{j,j } > 0$. 
\end{proof}

For example, if $(\up, \underline \tau, \underline \xi)$  is   geometrically regular  for $A$ in the sense of 
Definition \ref{def22} and if $B = \Delta_x \Id $ is an artificial viscosity, then 
$(\up, \utau, \uxi)$  is geometrically regular for $iA + B$. 
However, this condition is too restrictive for applications, in particular when 
$A$ and $B$ do not commute. 

 \begin{exam}
 \textup{If the root is totally nonglancing, then the decoupling condition is trivially satisfied
 since either $J_I$ or $J_O$ is empty. This applies to fast shocks in MHD.   } 
 \end{exam}
 
 \begin{cexam}
 \textup{Slow shocks in MHD do not satisfy the decoupling condition, see Section~\ref{SecMHD}.}  
 \end{cexam}
 
      The decoupling condition  is crucial in the construction of symmetrizers. 
      The second condition \eqref{334x} is more technical. One could expect that 
    with   the positivity Assumption (H1), one could always find an adapted basis. 
    This is not clear, except for mutliplicity 2 or symmetric systems.

   \begin{prop}  \label{prop414} 
    Suppose that  $(\up, \underline \tau, \underline \xi)$  is  geometrically regular of multiplicity $m$. 
  Assume that either $m = 2$  or that the symmetry assumption (H1') is satisfied. 
 There is a basis $\{ e_j \}$ adapted to $B$. 
 
 If in addition all the eigenvalues $\lambda_j$ are of the same type O or I, then 
 the condition (BS) is satisfied. 
\end{prop}

The proof is given Section~\ref{Sec6}. 

\subsection{The hyperbolic block structure condition}

We turn back to the construction of symmetrizers for nonelliptic blocks $\chH_k$ in 
the splitting \eqref{block2}. 
The construction of K-families is performed in \cite{MZ1} provided that 
$\chH_k$ can be put in a suitable normal form. This is the so called 
\emph{ block structure condition}. We first review this condition in the 
hyperbolic case, and next extend it  to the   hyperbolic-parabolic case.

Consider $\up$ and a frequency $\underline \cz = (\underline \ct, \underline \ch, 0) \ne 0$ and a purely imaginary eigenvalue 
\eqref{mur} $\underline \mu_k = i \underline \cx_d$ of  $H_0(\up, \underline \cz)$. 
Let $\underline \cx = ( \underline \ch, \underline \cx_d)$. Then 
$(\up, \underline \ct, \underline \cx) $ is a root of $\Delta$.  
We consider the  
block $\chH_k$ associated to $\underline \mu_k$ and denote by $\EE_k$  the corresponding
invariant space of $\chH$. We use the notations 
$\chH_{k,0} (p, \cz) = \chH_k(p, \cz, 0)$  and $\EE_{k,0} (p, \cz) = \EE_k(p, \cz, 0)$. 

\begin{defi}
\label{Ddefbl2}
$\chH_{k,0}  $ has the block structure property near $(\up, \underline \cz  )$ if there exists a smooth invertible matrix 
$V_{k, 0}  $ on a neighborhood of that point such that 
$V_{k,0}^{-1} \chH_{k, 0}  V_{k,0}  $ is block diagonal,
\begin{equation}
\label{block2b}
V_{k, 0} ^{-1} \chH_{k, 0}    V_{k,0}  =  
\left[\begin{array}{cccc}
Q_1  & 0 &  &   
\\
  0    &      \ddots &0\\
 &     0 & Q_{m'}
\end{array}\right]\, , 
\end{equation}
 with  diagonal blocks $ Q_j $
of size $\nu_j \times \nu_j$ such that :

  $Q_j(p, \cz)$  has purely imaginary  coefficients   
  when $\cg = 0$,   
 \begin{equation}
\label{Dbl1}
 Q_j (\up, \cz)    =  \underline \mu_k \Id +   i 
\left[\begin{array}{cccc}
0  & 1 & 0&   
\\
0  &0  & \ddots  &  0
\\
  & \ddots &      \ddots & 1 \\
 &  &  \cdots & 0
\end{array}\right]\, , 
\end{equation}
 and  the real part of the lower left hand corner of 
${\partial_{\cg}  Q_j} (\up, \underline \zeta)$, denoted by $  q^\flat_j$,  
does not vanish. 

\end{defi}
 
  When  $\nu_j = 1$,  $Q_j(p, \cz)$  is a scalar.  In this case, \eqref{Dbl1}   has to be understood 
  as $ Q_j (\up, \cz)    =  \underline \mu_k$, with no Jordan's block.  
  The lower left hand corner of the matrix is $Q_j$ itself and the condition reads 
$  q^\flat_j := \partial_{\cg} Q_j (\up, \underline \cz)  \ne 0$.

\begin{prop} [\cite{MZ2}] 
\label{prop33}
If the root $(\up, \underline \ct, \underline \cx)$
 is geometrically regular in the sense of Definition~$\ref{def22}$,  
 the corresponding block $\chH_{k, 0}$ satisfies 
 the block structure condition.
 
  Conversely, if  $\chH_{k, 0}$ satisfies 
 the block structure  condition  with 
 matrices $V$  that are real analytic in $\cz$, then the  root 
 $(\up, \underline \ct, \underline \cx)$
 is   geometrically regular. 
 
 \end{prop}
 
\begin{rem}
\textup{There is a slight discrepancy here between the necessary and the sufficient condition, due to  analyticity conditions.   Definition~$\ref{def22}$ requires analyticity in $\cx$. 
This is used in the proof of sufficiency. In addition, it  implies that 
the block structure  condition  holds with   matrices $V$  that are real analytic in $\cz$. 
Thus, there is an  ``if and only if'' theorem. However, for the construction of symmetrizers, 
analyticity of $V_k$ is not needed, this is why we do not insist on it in the definition above. 
In addition, note that \emph{for fixed $p$}, the existence of $C^\infty$ eigenvalues and eigenvectors
for $A$, implies that these eigenvalues are real analytic in $\xi$ and that one can choose analytic 
eigenvectors (see e.g \cite{Shi, Mal}. The question is to control the domain of analyticity 
as $p$ varies. In applications, for this problem, proving  analyticity is not harder than proving the  $C^\infty$ smoothness}.  
\end{rem}

 To prepare the hyperbolic-parabolic analysis, we have to review the proof of Proposition~\ref{prop33}. In particular, we 
  reformulate the conditions of Definition~\ref{Ddefbl2} in a more intrinsic way.  
 The choice of  a smooth matrix $V_{k, 0}$ is equivalent to the choice of  a smooth
  basis of $\EE_{k, 0} $, denoted by $\{ \vp_{j, a} (p, \cz)\}_{ 1 \le j \le m', 1 \le a \le \nu_j} $. 
  The property \eqref{Dbl1} reads 
  \begin{eqnarray}
\label{Dbl2}
&& (H_0 (\up, \underline \cz) - \underline \mu_k )  \vp_{j, 1}  (\up, \underline \cz) = 0, 
\\
\label{Dbl3}
&& (H_0 (\up, \underline \cz) - \underline \mu_k )  \vp_{j, a}  (\up,  \underline \cz) = 
i  \vp_{j, a -1}  (\up, , \underline \cz), \quad  2 \le a \le \nu_j. 
\end{eqnarray}
  With \eqref{b2b}, there is a unique smooth 
  dual basis $\psi_{j, a} (p, \cz)$ such that 
  \begin{equation}
\label{dualb}
\begin{aligned}
& \psi_{j, a} \cdot \EE'_{k,0} = 0 , 
\\
& \psi_{j, a} \cdot \vp_{j', a'} = \delta_{j, j'} \delta_{a, a'}. 
\end{aligned}
\end{equation}  
Here,  $\EE'_{k,0}$ denotes the invariant space of 
$H_0(p, \cz)$  such that $\CC^N = \EE_{k, 0} \oplus \EE'_{k, 0}$. 
It is the sum of invariant subspaces associated to 
eigenvalues $\underline \mu_{k'} \ne \underline \mu_k$.

In the basis $\vp_{j, a}$, the entries of the matrix $V^{-1}_{k, 0} \chH_{k, 0} V_{k, 0}$ are 
$\psi_{j, a} H_0 \vp_{j', a'}$. 
The diagonal block structure means that 
\begin{equation}
\label{Dbl4}
\psi_{j, a} H_0 \vp_{j', a'} = 0  \quad  \mathrm{when} \quad j \ne j'. 
\end{equation}
 The other conditions read: 
 \begin{eqnarray}
\label{Dbl5}&&  \re (  \psi_{j, a} H_0 \vp_{j, a'} )  = 0  \quad  \mathrm{when} \quad  \cg = 0 ,  
  \\
\label{Dbl6}
& &  \re  \D_{\cg} ( \psi_{j, \nu_j } H_0 \vp_{j, 1} )  (\up, \underline \cz) \ne 0.  
\end{eqnarray}
 
We  first show how to compute this quantity in terms of $A$ only.

\begin{lem}
\label{lem31}
Suppose that  $\chH_{k, 0}$ has a block diagonal
decomposition $\eqref{block2b}$ in a smooth basis $\vp_{j, a}$  of 
$\EE_k(p, \cz, 0)$ which satisfies $\eqref{Dbl2}$ $ \eqref{Dbl3}$. 
Let   $\psi_{j, a}$ denote a dual basis satisfying $\eqref{dualb}$. 
The lower left hand corner entry 
of  $\D_{\cg} Q_j (\up, \underline \cz)$ is  equal to the lower left hand corner entry of 
$ - i  \D_{\ct} Q_j (\up, \underline \cz) $  and equal to 
\begin{equation}
\label{qj}
\underline  q_j  =   -  \psi_{j, \nu_j} (\up, \underline \cz) \, A_d^{-1} (\up) \,   \vp_{j, 1} (\up, \underline \cz). 
\end{equation}
\end{lem} 
\begin{proof}
Let $\underline H_0 = H_0 (\up, \underline \cz)$. 
 Then 
$\underline H_0   - \underline \mu_k$ is invertible 
on $\EE'_{k, 0} (\up, \underline \cz) $. With \eqref{Dbl2} \eqref{Dbl3}, this implies that 
\begin{eqnarray}
\label{orth}
&&\mathrm{range} \big( \underline H_0 - \underline \mu_k \Id \big)  
= \{   \psi _{1, \nu_1} (\up, \underline \cz), \ldots ,   \psi _{m', \nu_{m'}}(\up, \underline \cz) \}^{\perp}, 
\\
\label{ker}
&&\ker  \big( \underline H_0 - \underline \mu_k \Id \big)  
= \{  \vp _{1, 1} (\up, \underline \cz), \ldots ,   \vp _{m', 1}(\up, \underline \cz) \}. 
\end{eqnarray}
In particular,  
\begin{equation}
\label{new3}
 \big( \underline H_0   - \underline \mu_k \Id \big) \vp_{j, 1} = 0 \quad \mathrm{
and} \quad  \psi_{j, \nu_j } \big( \underline H_0  - \underline \mu_k \Id \big) = 0. 
\end{equation}

  The entry in consideration is 
  $$
  q_j( p, \cz) =  \psi_{j, \nu_j } H_0 \vp_{j, 1} =  
  \psi_{j, \nu_j } \big( H_0  - \underline \mu_k \Id \big) \vp_{j, 1}  + \underline \mu_k \delta_{\nu_j, 1}. 
  $$
 Therefore, 
differentiating in $\cg$ and $\ct$ and using \eqref{defH0}, implies that 
\begin{equation}
\label{dtau=idgamma}
\D_{\cg} q_j (\up, \underline \cz) =  - i \D_{\ct} q_j (\up, \underline \cz)  = \underline q _j 
\end{equation}
is given by \eqref{qj}. 
\end{proof}

We now discuss how much flexibility there is in the choice of the basis 
$\vp_{j, a}$. 
Recall that we are considering   a purely imaginary eigenvalue 
$\underline \mu_k =  i \underline \xi_{d}  $ 
of $H_0(\up, \underline \cz)$,  so that $ - \underline \ct $  is an eigenvalue 
$\underline \lambda$
of $A(\up, \underline \cx)$ with $\underline \cx = (\underline \ch, \underline \cx_d)$.

\begin{lem} 
\label{lem32}
Suppose that  $\chH_{k,0}  $ has the block structure property near 
$(\up, \underline \cz  )$ in a smooth  basis $\vp_{j, a}$ and denote by  $\psi_{j, a}$ the 
dual  basis \eqref{dualb}.  Then, 

i) $\underline \lambda$ is a semi-simple eigenvalue of $A(\up, \underline \cx)$
with multiplicity $m$ equal to the number $m'$ of blocks $Q_j$, 

ii)  on a neighborhood of $(\up, \underline \cx)$, there are $m$ smooth 
eigenvalues $\lambda_j (p, \cx)$ of $A(p,  \cx) $ and $m$ smooth linearly independent eigenvectors $e_j(p, \cx)$, 
such that  
 \begin{eqnarray}
 \label{311}
&& \lambda_j(\up, \underline \cx) = \underline \lambda, 
\\
\label{312}
&&   A(p, \cx) e_j (p, \cx) = \lambda_j(p, \cx) e_j(p, \cx), 
\\
\label{313}
&& e_j(\up, \underline \cx) = \vp_{j, 1} (\up , \underline \cz), 
\end{eqnarray}

iii) the order of $\underline \cx_d$ as  a root  of 
$\underline \ct + \lambda_j( \up, \underline \ch, \cdot ) = 0$ is equal to $\nu_j$,

iv) denoting by $\{\underline \ell_j\} $ the left eigenvector dual basis of $\{ e_j \}$ as in 
$\eqref{dual1}$, there holds 
\begin{equation}
\label{dual2}
   \underline \ell_j  A_d(\up) = \beta_j \psi_{j, \nu_j} (\up, \underline \cz)   . 
\end{equation}
 with  $  \beta_j := \frac{1}{\nu_j !} \D^{\nu_j}_{\xi_d} \lambda_j (\up, \underline \cx) $
 as in $\eqref{multig}$, 
 
v) the lower left hand corner entry 
of  $\D_{\cg} Q_j (\up, \underline \cz)$ is  
\begin{equation}
\label{qj2}
\underline q_j = - 1 / \beta_j \in \RR.
\end{equation} 

\end{lem} 

\begin{proof} {\bf a) } Define $\widetilde \vp_{j, \nu_j} = \vp_{j, \nu_j}$
and for $a < \nu_j$ 
\begin{equation}
\label{nphi}
\widetilde \vp_{j, a} (p, \zeta) =   - i \big( H_0(p, \zeta) - \underline \mu_k \big) 
\vp_{j, \nu_j} . 
\end{equation}
By \eqref{block2b}\eqref{Dbl1}, there holds
\begin{equation}
\label{nphip}
\widetilde \vp_{j, a} (\up, \underline \cz ) =    \vp_{j, a} (\up, \underline \cz ). 
\end{equation}
Moreover, in the new  basis   $ \tilde \vp_{j, a} $, the  matrix of 
$  Q_j$ has the form 
\begin{equation}
\label{Dbl1n} 
 Q_j  =   i \underline \cx_d \Id  + i  \begin{pmatrix}  *  & 1 & \ldots &0
\\
\vdots & 0 & \ddots &0
\\ * & 0 & \ldots &1
\\  *  & 0 & \ldots &0
\end{pmatrix} 
\end{equation} 
Thanks to \eqref{nphip}, the dual basis $\{ \tilde \psi_{j, a} \} $ associated to 
$ \{ \tilde \vp_{j, a} \} $ also satisfies 
$\widetilde \psi_{j, a} (\up, \underline \cz ) =    \psi_{j, a} (\up, \underline \cz )$. 
This implies that the lower left hand corner of $\D_{\cg}  Q_j (\up, \underline \cz)$ is unchanged
in the new basis.

 \medbreak
 {\bf b) }  Consider the determinant
\begin{equation*}
\Delta_j ( p, \cz, \cx_d) = \det \big( \xi_d \Id + i Q_j (p, \cz) \big). 
\end{equation*}
It is independent of the basis  $\{ \psi_{j, a} \} $ or  $\{ \tilde \psi_{j, a} \} $. 
Thus, it  is real when $\cg = 0$ and vanishes at $(\up, \underline \cz, \underline \cx_d)$. 
Moreover,  \eqref{Dbl1}  implies that
\begin{equation*}
 \D_{\ct} \Delta_j( \up, \underline \cz, \underline \xi_d) =  - \underline q_j . 
\end{equation*}
As a byproduct,  using also \eqref{dtau=idgamma}  this shows that 
\begin{equation}
\label{317}
\underline q_j \in \RR \quad  \mathrm{thus} \quad \underline q_j = \re \underline q_j = q^\flat_j \ne 0 . 
\end{equation}

In particular, the implicit function theorem implies that there is a smooth function 
$\lambda_j(p, \cx)$, in a real neighborhood of $(\up, \underline \cx)$, such that 
$\lambda_j(\up, \underline \cx) = -  \underline \ct $ and for 
$\cz = (\ct, \ch, 0)$: 
\begin{equation}
\label{318}
\Delta_j( p, \cz, \cx_d) =  \alpha_j ( p, \cz, \cx_d) \big( \ct + \lambda_j (p, \cx) \big)  
\end{equation}
with $\alpha_j (\up, \underline \cz, \underline \cx_d) \ne 0$. 

\medbreak
{\bf c) } Consider next the eigenvector equation 
\begin{equation}
\label{eigveceq}
\big( \cx_d \Id + i Q_j (p, \cz) \big) e_j   = 0  . 
\end{equation}
 By \eqref{Dbl1n}, in   the basis  $\{ \tilde \psi_{j, a} \} $,   the $\nu_j - 1$ first equation   determine the last $\nu_j - 1$ components of $e_j$ 
 \begin{equation}
 \label{hpxx}
 (e_{j})_{ a}  =  (\cx_d - \underline \cx_d)^{a -1} (e_{j})_{ 1} , \quad   a \ge 2 . 
 \end{equation}
  Substituting these values,  the last equation
 is a scalar equation equivalent to $\Delta_j = 0$. 
Introduce  
\begin{equation*}
\zeta_j (p, \ch, \cx)  =   \big( - \lambda_j(p, \cx), \ch, 0 \big), 
\end{equation*}
and 
\begin{equation}
\label{defejnnn}
e_j (p, \cx) = \tilde \vp_{j, 1} (p, \cz) + \sum_{a=2}^{\nu_j}   (\cx_d - \underline \cx_d)^{j-1} 
\tilde \vp_{j, a} (p, \cz ) . 
\end{equation}
This vector  is smooth and satisfies  \eqref{eigveceq}, thus 
 \begin{equation*}
\big(A(p, \cx) - \lambda_j (p, \cx) \Id \big) e_j (p, \cx) =  
A_d(p)  \big(  i H_0 (  p, \cz_j) + \cx_d \Id \big) e_j (p, \cx)  = 0. 
\end{equation*}
Moreover, the $e_j( \up, \underline \cx) = \vp_{j, 1} (\up, \underline \cz)$   are linearly independent. 

\medbreak

{\bf d) } By \eqref{318}, for $\cz = (\ct, \ch, 0)$, there holds 
\begin{equation*}
\begin{aligned}
\det \big(\ct \Id + A(p, \cx) \big)  & = \det (A_d)\,  \det \big( i H_0 (p, \cz) + \cx_d \Id  \big) 
\\ & = \alpha (p, \ct, \cx) \prod_{j= 1}^{m'}   \big( \ct + \lambda_j(p, \cx) \big)
\end{aligned}
\end{equation*}
where $\alpha (\up, \underline \ct, \underline  \cx) \ne 0$ and $m'$ 
is the number of blocks $Q_j$. This shows that 
$  -  \underline \ct $ is an eigenvalue of algebraic   order  $m' $ of 
$  A(\up, \underline \cx) \big)  $. By step c), the geometric multiplicity is at least $m'$, 
implying that $  -  \underline \ct $ is semi-simple of   order $ m'$. 

Moreover, by \eqref{Dbl2}, there holds
\begin{equation*}
\Delta_j (\up, \underline \cz, \cx_d) =  ( \cx_d - \underline \cx_d) ^{\nu_j}, 
\end{equation*}
showing that $\underline \cx_d$ is a root of multiplicity $\nu_j$ of $\Delta_j $, 
thus of $\underline \ct + \lambda_j(\up, \underline \ch, \cx)  = 0$. 

\medbreak
{\bf e) } Let $\underline \ell_j$  satisfy \eqref{dual1}. 
Thus  
\begin{equation*}
\begin{aligned}
\mathrm{Range} \big( \chH_0(\up, \underline \cz) - \underline \mu_k \Id \big)  
&=  A_d^{-1} (\up)\mathrm{Range} \big(\underline \ct  \Id + A(\up, \underline \cx) \big)   
\\
& 
=  A_d^{-1} (\up) \{ \underline \ell _{1}, \ldots , \underline \ell_{m} \}^{\perp}.  
\end{aligned}. 
\end{equation*}
Comparing with \eqref{orth}, this implies that 
\begin{equation}
\label{319}
\mathrm{span} \big\{ \psi_{j, \nu_j} (\up, \underline \cz) , 1 \le j \le m \big\} = 
\mathrm{span} \big\{   \underline\ell_j , 1 \le j \le m \big\} .  
\end{equation}

For $a \in \{ 1, \ldots, \nu_j \}$, introduce 
\begin{equation}
\label{base}
\underline e_{j, a}   = \frac{1}{ (a-1)! }  \D_{\xi_d} ^{a-1} e_j (\up, \underline \cx). 
\end{equation} 
Because $\underline \cx_d$ is a root of order 
$\nu_j $ of $\underline \ct  + \lambda_j( \up, \underline \ch, \cx)= 0$, the definition
\eqref{defejnnn} implies that  
\begin{equation*}
\underline e_{j, a } =  \tilde \vp_{j, a}  (\up, \underline \cz) 
= \vp_{j, a} (\up, \underline \cz)   \quad   \mathrm {for}   \  1 \le  a \le  \nu_j .
\end{equation*}
In particular,   \eqref{dualb} implies that 
\begin{equation}
\label{322}
\psi_{j', \nu_{j'}} (\up, \underline \cz) \cdot \underline e_{j, \nu_j} = 
\psi_{j', \nu_{j'}} (\up, \underline \cz) \cdot \vp_{j, \nu_j}(\up, \underline \cz) = \delta_{j, j'}. 
\end{equation}

 Differentiating    the equation 
 \begin{equation}
 \label{320}
\Big( A(\up, \cx) - \lambda_j (\up,    \cx) \Big) 
e_j ( \up, \cx) = 0 
\end{equation} 
with respect to $\cx_d$  and at order $\nu_j $ yields 
\begin{equation*}
 \big( \underline \ct \Id + A(\up, \underline \cx)\big)  \D_{\xi_j}^{\nu_j} e_j (\up, \underline \cx) 
 = - \nu_j  A_d(\up)  \D^{\nu_j-1}_{\xi_d} e_j (\up, \underline \cx)  + 
  \D_{\xi_j}^{\nu_j } \lambda_j (\up, \underline \cx)  e_j (\up, \underline \cx). 
\end{equation*}
Multiplying on the left by $ \underline \ell_{j'}  $ annihilates the left hand side, implying 
\begin{equation*}
 \underline \ell_{j'} A_d(\up) e_{j, \nu_j} (\up, \underline \cz) = 
\beta_j  \underline \ell_{j'} \cdot e_j (\up, \underline \cx) = \beta_j \delta_{j', j}. 
\end{equation*}
By \eqref{319},  the $\underline \ell_j A_d$ and $\psi_{j, \nu_j}$ span the same space. , 
 Therefore, comparing with \eqref{322} implies that 
 $ \underline \ell_{j'} A_d(\up) = \beta_j \psi_{j', \nu_{j'}} (\up, \underline \cz)$. 
 
 \medbreak
 {\bf f) } By \eqref{qj} and \eqref{dual2}, we have
 $$
- \beta_j  \underline q_j =    \underline \ell_j \vp_{j, 1} (\up, \underline \cz) = 
\underline \ell_j e_j(\up, \underline \cx) = 1. 
 $$
The proof of the lemma is complete.  \end{proof}

  \begin{rem}
  \label{rem33}
  \textup{This lemma  is a variation on the  necessary part in Proposition~\ref{prop33}
  (see \cite{MZ2}), 
  with useful additional remarks. It 
  shows that the block structure condition is 
  closely related to a smooth diagonalisation of $A$.   Conversely, if one starts from a 
  smooth basis $e_j$ and a  root of 
  $\underline \ct + \lambda_j (\up, \underline \cx)$ with 
  \eqref{multig}, one constructs a basis $\vp_{j, a}$ such that $\vp_{j, a} (\up, \underline \cz)$
  is given by \eqref{base}, using an holomorphic extension of $e_j$  to complex values of 
  $\cx_d$ (see \cite{MZ2}). Lemma \ref{lem32} implies that the change of bases which preserve
  the block structure form are linked to change of bases which preserve the smooth 
  diagonalization of $A$. }
  \end{rem} 
 
 The construction of K-families of symmetrizers for the blocks $Q_j $ is performed in 
 \cite{Kr, Ma1, Met3}. The sign of $\beta_j$ and the parity of $\nu_j$ play an important 
 role. Hyperbolicity implies that   $H_0$ and thus the $\chH_k$ and $Q_j$ have no 
 purely imaginary eigenvalues when $\cg > 0$. 
 Denote by $\EE^-_{Q_j}$ the invariant space of $Q_j$ associated to 
 the spectrum in $ \{ \re \mu < 0 \}$
  since the definition of the limiting space $\underline \EE_{Q_j} ^-$. Recall 
  that the limit space at $(\up, \underline \cz)$ is 
 \begin{equation}
 \label{325y}
  \underline \EE_{Q_j}^- =  \CC^{ \nu'_j } \times \{ 0 \}^{\nu_j - \nu'_j} 
  \end{equation}
  with 
  \begin{equation} 
  \label{326y}
    \nu'_j = \left\{ \begin{array}{llllll}
 \nu_j / 2   & & \mathrm{when } & \nu_j \ \mathrm{ is \ even}, 
 \\
  (\nu_j + 1)  / 2   & & \mathrm{when } & \nu_j \ \mathrm{ is \ odd} & \mathrm{and} 
  & \beta_j > 0, 
  \\
  (\nu_j - 1)  / 2   & & \mathrm{when } & \nu_j \ \mathrm{ is \ odd} & \mathrm{and} 
  & \beta_j < 0. 
 \end{array}\right.
 \end{equation} 
 
 \begin{rem}
 \textup{As a corollary, we have the following characterization of  the sets $J_0$ and $J_I$:   }
 \begin{equation}
\label{4J}
\left\{ \begin{aligned} &  j \in J_I  \quad  if \ \nu_j \ is \ even \ or \ \nu_j \  is \ odd \ and  \ q_j^\flat < 0, 
\\
& j \in J_0  \quad  if \  \nu_j \  is \ odd \ and \  q_j^\flat > 0. 
\end{aligned}\right. 
\end{equation}
 \end{rem}

 
 \subsection{The hyperbolic-parabolic case}

We still consider a block $\chH_k$ associated to a purely imaginary eigenvalue 
\eqref{mur}.  In the next section, we show that the following technical conditions are  
the natural one for the 
 construction of  Kreiss symmetrizers. 

 \begin{defi}
 \label{BS2}
$\chH_{k}  $ has the generalized block structure property near $(\up, \underline \cz , 0 )$ if there exists a smooth invertible matrix 
$V_{k}  $ on a neighborhood of that point such that 
\begin{equation}
\label{bdiag2}
V_{k}^{-1} \chH_{k}  V_{k}  = \left(\begin{array}{ccc}
  Q_1  & \cdots  & 0
\\
\vdots  &\ddots  & \vdots
\\
0 & \cdots  & Q_m
\end{array}\right)  +   \rho \left(\begin{array}{ccc}
  \tilde B_{1, 1}   & \cdots  & \tilde B_{1, m} 
\\
\vdots  &\ddots  & \vdots
\\
\tilde B_{m, 1}  & \cdots  & \tilde B_{m,m}
\end{array}\right)
\end{equation}
 where the $Q_j (p, \cz)$  satisfy the properties of Definition~$\ref{Ddefbl2}$. 
 Moreover, the   $m \times m$ matrix $B^\flat$  with entries $B^\flat_{j, j'}$ equal to 
the lower left hand corner of  $\tilde B_{j, j'}(\up, \underline \cz, 0)$ satisfies  
\begin{equation}
\label{328x}
 B^\flat_{j, j'} = 0  \quad when \quad  (  j ,  j') \in (J_O \times J_I ) \cup  ( J_I \times J_O)
\end{equation}
where $J_0$ and $J_I$ are defined by $\eqref{4J}$ and there is a real diagonal matrix $D^\flat$, with entries 
$d_j^\flat$  such that 
\begin{equation}
\label{bls2}
d^\flat_j  q^\flat_j > 0 , \quad  \re D^\flat B^\flat > 0.  
\end{equation}
 \end{defi}

We   show that  these conditions are related to the  condition $(BS)$ of Definition~\ref{defBS}
formulated on the original system. 
We need  first a more detailed form of the block reduction  $H$ in \eqref{block1}. 
 Introduce  the following notations:  
\begin{eqnarray}
   \label{eq643b}
&  B_{* *} (p, \zeta) & := \sum_{j=1, k }^{d-1}  \eta_j  \eta_k  B_{j, k} (p) , 
\\ \label{eq643c}
& B_{* d} (p, \zeta) & := \sum_{j=1 }^{d-1}  \eta_j   (   B_{j, d} (p)  + B_{d,j}(p) ) 
\end{eqnarray}

\begin{lem}
\label{lem642b}
 
 One can choose the matrix $V$ in \eqref{block1} such that 
there holds 
 \begin{equation}
 \label{neqH}
 H(p, \zeta) =  H_0(p, \zeta)  -  H_1 (p, \zeta) + O(\vert \zeta \vert^3)
 \end{equation}
where 
\begin{equation}
  \label{eq644}
 H_1  =   A_{d}^{-1} \Big(  B_{*,*} -  i    B_{*,d} H_0  -    B_{d,d} H^2_0\Big),   
  \end{equation}
   
\end{lem}

 \begin{proof}
   Direct computations show that  
 the kernel of $G(p, 0) $ is 
 $\CC^N \times \{0\}$  and, using that $A_d$ is invertible, that  $\ker G(p, 0) \cap \range G(p,0) = \{0 \}$
  This shows that   
 $0$ is a semi-simple eigenvalue of $G(p, 0)$. 
  
 If $\mu  $ is a purely imaginary eigenvalue of 
 $G(p, 0)$, then $0$ is an eigenvalue of 
 $ i A(p, \xi) +  B(p, \xi) $ with $\xi = (0, - i \mu)$. By Assumption (H1) 
 this requires that $\xi = 0$, thus $\mu = 0$. 
 This shows that the nonvanishing eigenvalues of $G(p, 0)$ are not 
 on the imaginary axis.

 This implies that there is a smooth matrix $V(p, \zeta)$ on 
a neighborhood of $(\up, 0)$ such that 
\eqref{block1} holds with $H(p, 0) = 0$ and 
$P(p, 0) $ invertible  with no eigenvalue  on the imaginary axis. 

 The image of the  first $N$ columns  of $V$  is the invariant space 
 of $G$, and $H$ is the restriction of $G$ to that space. 
 At $\zeta= 0$ this space is $\ker G$, and 
 performing a smooth change of basis in $\CC^N$, we can always assume that 
 the first $N$ columns  of $V$ are of the form 
 \begin{equation}
 \label{colV} 
V_I (p, \zeta) = \begin{pmatrix}
     \Id_{N \times N}     \\
      W (p, \zeta)    
\end{pmatrix}
\end{equation}
with $W$ of size $N' \times N$ vanishing at $\zeta = 0$. 
This implies   \eqref{Vtri}.

 By \eqref{block1}
$ G V_I = V_I H $, hence $M V_I  = G_d V_I H$ and 
\begin{equation*}
 \cM    =  - \cA H + \oB_d W H , 
\quad 
 W   =  J H . 
\end{equation*}
Therefore, 
\begin{equation}
\label{conjH}
 \cM    =  - \cA H + \oB_d J H^2 = - \cA H + B_{d,d} H^2. 
\end{equation}
Taking the first order term at $\zeta = 0$ shows that the first order term 
in $H_0$  in $H$ satisfies
$$
 ( i \tau  + \gamma) \Id  +
\sum_{j=1}^{ d-1}   i \eta_j A_j = -A_d(p) H_0
$$
 and hence is given by \eqref{defH0}. The second order term $H_1$ in $H$  
satisfies
 $$
 B_{*,*} = - A_d H_1 +  i B_{*, d} H_0  +  B_{d, d} J H_0^2 
 $$
   implying \eqref{neqH} and \eqref{eq644}. 
   \end{proof}

Parallel to Lemma~\ref{lem31}, we can now state:  
 
 \begin{lem}
\label{lem37}
Suppose that  the matrix of $\chH_{k}$ is given by the right hand side of 
$\eqref{bdiag2}$ in a smooth basis $\vp_{j, a}$  of 
$\EE_k(p, \cz, \rho)$ which satisfies $\eqref{Dbl2}$ and $ \eqref{Dbl3}$ for $\rho = 0$. 
Let  $\{ \underline \ell_j\} $  denote the  dual basis of 
$\{ e_j = \vp_{j, 1} \}$ satisfying $\eqref{dual1}$. 
 The entries of $B^\flat$ are 
\begin{equation}
\label{bjj}
B^\flat_{j, j'} 
 =  -  \frac{1}{\beta_j}  \underline \ell_j    \, B(\up, \underline \cx)  \, 
   \vp_{j', 1} (\up, \underline \cz, 0). 
\end{equation}
\end{lem}

\begin{proof} In the block reduction 
\eqref{bdiag2},  the lower left hand corner entry of the $(j,j')$-block  is 
\begin{equation*}
h_{j, j'} = \psi_{j, \nu_j} \chH \vp_{j', 1}  =  \psi_{j, \nu_j} \big( \chH - \underline \mu_k) \vp_{j', 1} 
+ \underline \mu_k \delta_{j, j'}. 
\end{equation*}
 Differentiating   in $\rho$   
 and using the relations \eqref{new3}  yields 
\begin{equation*}
- B^\flat_{j, j'} = \D_{\rho } h_{j, j'}  (\up, \underline \cz , 0) =  -  \underline \psi_{j, \nu_j} 
  \,  \tilde B (\up, \underline \cz) \, 
\underline \vp_{j, 1}   , 
\end{equation*}
where  $\underline \psi_{j, \nu_j} $ and $\underline \vp_{j, 1}$ stand for the 
evaluation at  $(\up, \underline \cz, 0)$ of the corresponding function. 
Using the explicit form of $\tilde B$ and the relations 
$$
\underline  H_0  \underline \vp_{j, 1}  = i \underline \cx_d \underline \vp_{j, 1}, 
\quad
\underline \psi_{j, \nu_j} \underline  H_0   = i \underline \cx_d \underline \psi_{j, \nu_j} 
$$
 we obtain
 $$
 \begin{aligned}
 \underline \psi_{j, \nu_j} \tilde B (\up, \underline \cz) \underline \vp_{j, 1}  & = 
  \underline \psi_{j, \nu_j}  A_d^{-1}  \Big(  B_{*,*}  (\up, \underline \ch  ) 
  +   \underline \cx_d B_{*; d} (\up, \underline \ch)  + 
  \underline \cx_d^2 B_{d, d} (\up) \Big) \underline \vp_{j, 1} 
  \\ 
  & =   \underline \psi_{j, \nu_j}   B (\up, \underline \cx) \underline \vp_{j, 1} 
 \end{aligned}
 $$
 With \eqref{dual2},  this implies \eqref{bjj}. 
\end{proof}

\begin{theo}
\label{th39}
   If   $(\up, \underline \ct, \underline \cx)$ is a geometrically regular 
characteristic root of $\Delta$ which satisfies the condition (BS) of 
Definition~$\ref{defBS}$. Then the associated block 
$\chH_k$ satisfies the generalized block structure condition.

\end{theo} 

\begin{proof}
Since  $(\up, \underline \ct, \underline \cx)$ is  geometrically regular, the hyperbolic part
$\chH_{k, 0}$ satisfies the block structure condition. 
Moreover, if $e_j$ is a basis analytic in $\xi$,  there is a basis $\varphi_{j, a}$ such that 
$\varphi_{j, a} (\up, \underline \cz) = e_j(\up, \underline \cx)$ (see Remark~\ref{rem33} or  
\cite{MZ2}). By Lemma \ref{lem37}, \eqref{331x} is equivalent to  \eqref{328x}. 

If once can choose the base $\{e_j\} $ such that 
\eqref{334x} holds, then choose   $d^\flat_j = - \beta_j$ and by \eqref{qj2} and \eqref{bjj}
there holds $d_j^\flat q_j^\flat = 1$ so that 
  $ D B^\flat =  B^\sharp$ satisfies \eqref{bls2}. 
  \end{proof} 

\begin{rem}
\textup{Conversely,  if the generalized block structure condition
holds with matrices $V_k$ which are real analytic in $\cz$, then, by Proposition~\ref{prop33}
$(\up, \underline \ct, \underline \cx)$ is  geometrically regular.  
By \eqref{bjj}, \eqref{328x} is equivalent to 
the decoupling condition \eqref{331x}. Moreover,  \eqref{bls2} implies that there is a diagonal matrix 
with positive entries $d^\sharp_j = d^\flat_j / q^\flat_j$ such that 
$\re D^\sharp B^\sharp > 0$.  Consider the diagonal matrix 
$C = (D^\sharp)^{-1/2} = \mathrm{diag} (c_j)$ and  the 
new basis   $\tilde e_j = c_j e_j$. The new dual basis 
is $\tilde{\underline \ell_j} = c_j^{-1} \underline c_j$ and the new matrix 
$\tilde B^\sharp$ is  $C^{-1} B^\sharp C = C D^\sharp B^\sharp C$ and therefore 
$\re \tilde B^\sharp  = C \re (D^\sharp B^\sharp) C$ is definite positive.} 
\end{rem}


\section{Symmetrizers}

In this section, we prove Theorems~\ref{th221} and \ref{th222}. 
We are given a frequencies  $\underline \cz = (\underline \ct,  \underline \ch,   0)$
and  a purely imaginary eigenvalue $\mu_k = - i \underline \cx_d $ of 
$H_0(\up, \underline \cz)$, so that 
$(\up, \underline \ct, \underline \cx$, with $\underline \cx = (\underline \ch, \underline \cx_d)$
is a root of the characteristic determinant $\Delta$, of multiplicity $m$. 
Our goal is to construct K-families of symmetrizers for the block
$\chH_k(p, \cz, \rho)$ associated to $\mu_k$.

\subsection{Proof of Theorem~\ref{th222}}

We assume here that $(\up, \underline \ct, \underline \cx$ is geometrically regular
and satisfies the condition $(BS)$. 
 
We follow closely  \cite{MZ1} (Lemma~4.11 and Appendix A therein. See also
\cite{Met3}) where the constant multiplicity case is studied. In this case,  
all the blocks $Q_j$ are  equal and thus have  the same dimensions $\nu$, 
but more importantly, all the eigenvalues are of the same type $O$ or $I$. 
 So we review the main steps of the construction
and indicate where  the proof of \cite{MZ1, Met3} has to be  modified.

In the block reduction  
\eqref{bdiag2} of $\chH_k$, we choose the symmetrizers $\Sigma^\kappa_k  $  to be block diagonal:
\begin{equation}
\label{defcsk}
\begin{aligned}
\Sigma^\kappa_k  & =
 \left(\begin{array}{ccc}
  S^\kappa_1  & \cdots  & 0
\\
\vdots  &\ddots  & \vdots
\\
0 & \cdots  & S^\kappa_m
\end{array}\right)\\
S^\kappa_j(p, \cz, \rho) & =  E^\kappa_j + \widetilde E^\kappa_j(p, \cz)  -
 i \gamma  F^\kappa_j - i \rho \widetilde F^\kappa_j
\end{aligned}
\end{equation}
where  $E^\kappa_j$ and $\widetilde E^\kappa_j $ are real symmetric matrices,
and $F^\kappa_j$ and $\widetilde F^\kappa_j$ are real and skew symmetric.
Moreover, $E^\kappa_j$ ,  $F^\kappa_j$ and $\widetilde F^\kappa_j$
 are constant, $\widetilde E^\kappa_j $ depends only on $(p, \ct , \ch)$ 
and 
the  $E^\kappa _j$ have the special form
$$
E^\kappa _j =
\left[\begin{array}{ccccc}
0 & \cdots & \cdots  &  0 & e^\kappa _{j, 1}
\\
\vdots  &    &  &\adots & e^\kappa_{j, 2}
 \\ \vdots & & \adots{}  & \adots
\\
0 & \adots  & \adots \\
 e^\kappa _{j, 1} &  e^\kappa _{j, 2}   &  & &  e^\kappa_{j, \nu_j}
\end{array}\right],
$$
and
$\widetilde E^\kappa_j(\underline p, \underline \cz) = 0$.

The block structure condition implies that 
\begin{equation} 
  \chH_k   = \diag(Q_j{}_{\vert \gamma = 0} )   +  
   \gamma  \diag (\D_\gamma Q_j {}_{\vert \gamma = 0} )  + \rho \widetilde B  _{\vert \rho = 0}  
    O(\gamma^2 + \rho^2) 
\end{equation} 
 $\Sigma^\kappa_k$ is  a symmetrizer   for $\chH_k$, on a neighborhood (depending on $\kappa$)
of $(\up, \underline \cz, 0)$,  if 
\begin{eqnarray}
\label{43y} &   &  \re \big( ( E_j^\kappa + \widetilde E_j^\kappa ) Q_j {}_{ \vert \gamma = 0} \big) = 0,   \\
 \label{44y} &   & \re \big( E_j^\kappa  \D_\gamma Q_j (\up, \underline \cz, 0)  -  i F^\kappa_j  Q_j(\up, \underline \cz) \big)   > 0, 
  \\
   \label{45y} &   & \re \big( \diag( E_j^\kappa  ) \widetilde B (\up, \underline \cz, 0) 
    -   i \diag (\widetilde F^\kappa_j  Q_j(\up, \underline \cz) ) \big)   > 0. 
\end{eqnarray}
Moreover, the condition \eqref{243y} reads 
\begin{equation}
\label{Cond1}
(E^\kappa_j w, w )  \ge
C_1 \Big( \kappa \vert \Pi_j^+ w \vert^2 - \vert \Pi_j^- w \vert^2\Big) ,
\end{equation}
where $\Pi_j^\pm$ is the projection onto
$\underline \EE_j^{\pm} $  in the  decomposition
$\CC^{\nu_j} = \underline \EE_j^{-} \oplus \underline \EE_j^{+}$, where 
\begin{equation}
\label{47y}
 \underline \EE_j^{-} = \CC^{\nu_j'} \times \{0\}^{\nu_j - \nu_j'}  , \quad  
 \underline \EE_j^{+} = \{0\}^{\nu_j'} \times \CC^{\nu_j - \nu_j'} , 
\end{equation} 
with $\nu'_j$ given by \eqref{326y}.

Before starting the construction,  we note that Ralston's  lemma \cite{Ra} (see also \cite{MZ1, Met3}) implies that  one can perform
 an additional  change of basis $\Id + \rho \tilde V$ such that the matrices $\widetilde B_{j,j'}$ in 
 \eqref{bdiag2} are of the form 
 \begin{equation}
\label{45}
\widetilde B_{j, j'}(\up, \underline \cz)  = \begin{pmatrix}
     * &  0 & \ldots & 0  \\
     \vdots  &  0 & \ldots & 0 \\
    B^\flat_{j, j'}  & 0 & \ldots & 0 
\end{pmatrix}. 
\end{equation}
 This does not affect the previous choices, made at $\rho =0$. 
 Next, we introduce nome notations. 
 A vector $w \in \CC^{N_k} = \oplus \CC^{\nu_j} $,
  is broken into $m $ blocks
$w_j \in \CC^{\nu_j}$,   with
components denoted by 
$w_{j, a}$. 
 We now proceed to the construction of the symmetrizers. 
 
 \medbreak

{\bf a) }  One first choose the $E_j^\kappa$ such that \eqref{Cond1} holds and 
\begin{eqnarray}
\label{minek}
&&\re \big(E^\kappa _j \partial _\gamma Q_j(\underline p, \underline \cz)  w_j , w_j
\big) \ge 2
\vert w_{j,  1}
\vert^2 - C_\kappa  \vert w'_j \vert^2
\\
\label{Cond4}
&&\re \big( \diag(E_j^\kappa) \widetilde B (\underline p, \underline \cz, 0)    w_j
, w_j \big)
\ge 2 \vert w_{*,1} \vert^2 - C'_\kappa  \vert w'_* \vert^2.
\end{eqnarray}
with $w_{j,  1}$ denoting the first component of $w_j  \in \CC^{\nu_j} $ and
$w'_j\in \CC^{\nu_j-1} $ denotes the other components and  
  $w_{*, 1} \in \CC^{m} $ is the collection of the first
components $w_{j,1}$ while 
$w_*'$ denotes the remaining  components. 

 Note that 
 \begin{equation*} 
\begin{aligned}
&\re \big(E^\kappa _j \partial _\gamma Q_j(\underline p, \underline \cz)  w , w
\big)  =  e_{j,1}^\kappa  q_j^\flat | w_{j, 1} |^2   +  O( | w_j | | w'_j | ) , 
\\
&\re \big( \diag(E_j^\kappa) \widetilde B (\underline p, \underline \cz, 0)    w_j
, w_j  \big)
 =  \re \big( E^\flat B^\flat   w_{*,1} , w_{*, 1} \big)  + O ( \vert w'_* |  | w | ) ,  
\end{aligned}
 \end{equation*}
 where $E^\flat$ is the $m \times m $ diagonal matrix with entries
 $e^\kappa_{j,1}$.  Moreover, the decoupling condition 
 \eqref{331x} implies that $B^\flat$ has a block diagonal structure: 
 ordering the base $\{e_j\}$ according to the type $I$ or $O$, with obvious
 notations there holds :  
 \begin{equation}
\label{411y}
B^\flat = \begin{pmatrix}
    B^\flat_I   &  0  \\
     0 &  B^\flat_O
\end{pmatrix}
\end{equation}
 Similarly we note $E^\flat = \diag (E^\flat_I, E^\flat_O)$ and 
  \eqref{minek}  and \eqref{Cond4} are satisfied
 if 
 \begin{equation}
\label{Cond5}
e_{j,1}^\kappa  q_j^\flat \ge 3, \quad   \re  E^\flat_I B^\flat_I    \ge 3 \Id, 
\quad      \re  E^\flat_O B^\flat_O   \ge 3 \Id. 
\end{equation}

 On the other hand, to satisfy \eqref{Cond1}, one chooses the $e^\kappa_{j, a}$ inductively, 
 starting from $a = 1$, but this choice depends on the type of the eigenvalue.  Remember also from \eqref{qj2} that 
$ \beta_j = - 1 / q_j^\flat$. According to  
 \cite{Kr, CP} or Lemma 8.4.2 in  \cite{Met3}, the $e_{j, a}^\kappa$ are chosen as follows. 
 
 \quad  1)   If   $  \lambda_j  $ is  of type $I$ ,  then $e^\kappa_{j, 1}= e_{j, 1}$  is taken $O(1)$, independent 
 of $\kappa$,   and  the $e^\kappa_{j, a}$ for 
 $a \ge 2$ are chosen successively and depend on $\kappa$.  In particular,  when 
 $\nu_j $ is even,   $e^\kappa_{j, 2} \ge c \kappa$.   When $\nu_j$ is odd and 
 $\beta_j > 0$, then $q_j^\flat < 0$ and 
  $e_{j, 1}  < 0 $;   when $\nu_j \ge 3 $, then 
 $e^\kappa_{j, 3} \ge c \kappa$.

 \quad 2)  If   $  \lambda_j  $ is  of type $O$, that is    $\nu_j$   odd and 
 $q^\flat_j > 0$, one   chooses  $e^\kappa_{j, 1} \ge c \kappa$
 and the other $e^\kappa_{j, a}$ are chosen inductively.

 By assumption, there is a diagonal real matrix $D^\flat = \diag (D^\flat_I, D^\flat_O)$ such that 
    \begin{equation*}
\label{Cond5*}
d^\flat_j   q_j^\flat > 0 , \quad   \re  D^\flat_I B^\flat_I   > 0, 
\quad      \re  D^\flat_O B^\flat_O   > 0  . 
\end{equation*}
Therefore, there is a positive constant $c$ such that  if we choose $e^\kappa_{j, 1}   =  c  d^\flat_j$ 
when $\lambda_j$ is of type $I$ and 
$e^\kappa_{j, 1} = c \kappa d^\flat_j $ when $\lambda_j$ is of type $I$, the condition \eqref{Cond5}  
 is satisfied.    Next, according to   \cite{Kr, CP, Met3}, we can choose  the $e^\kappa_{j, a}$ for $a \ge 2$ such that 
the inequality \eqref{Cond1} is also satisfied. 
 
 \begin{rem}
 \label{rem41}
 \textup{The construction above shows that the conditions of Definition~\ref{BS2} 
 are more or less necessary for the construction of K-families of symmetrizers. 
First,  the different  magnitude in $\kappa$ of 
 $e^\kappa_{j, 1}$ for different types forces the decoupling \eqref{411y}, that it condition \eqref{331x}. Second,  a spectral condition on $B^\sharp$ is not sufficient in general 
 to insure the existence of a diagonal matrix $E^\flat$ such that \eqref{Cond5} holds. 
 This indicates that condition \eqref{bls2} is also necessary for the construction above. 
 } 
 \end{rem}
 
 \medbreak
 {\bf b) }  Once the matrices $E_j$ are chosen, the construction goes on 
 as in \cite{MZ1, Met3}. We omit the details. 
By \eqref{Dbl1}, $ \re ( E_j Q_j(\up, \underline \cz)) = 0 $.
Next, using the implicit function theorem and the property that
$  \frac{1}{i} Q_k $ is real when $\cg = 0$,  the real symmetric matrix 
 $\widetilde E_k(p, \ct, \ch )$   is chosen
so that  such that 
$ \re (E_j + \widetilde E_j) (  Q_j{}_{\vert \gamma = 0}   ) = 0  $. 

 Since $F_j$ is real and skew symmetric, there holds 
   $ \re  -  i F_j Q_j (\up, \underline \cz) =  \re F_j J_j $ 
where $J_j$ is the Jordan matrix in \eqref{Dbl1}. One can choose $F_j$  such that
$$
\re (F_j J_j    w_j  , w_j    ) \ge
 -  \vert w_{j, 1} \vert^2 +  (C+1)   \vert w'_j \vert^2.
$$
where $C$ is the constant  in \eqref{minek}.
Adding to \eqref{minek}  implies \eqref{44y}. 

 Similarly,  $ \re  -  i \widetilde F_j Q_j (\up, \underline \cz) =  \re \widetilde F_j J_j $ 
   and  one can choose $\widetilde F_j$ such that 
 $$
\re (F_j J_j    w_j   , w_j    ) \ge
 -  \vert w_{j, 1} \vert^2 +  (C'+1)   \vert w'_j  \vert^2.
$$
where $C'$ is the constant  in \eqref{Cond4}, 
  implying  \eqref{45y}.


\subsection{Proof of Theorem~\ref{th221}}

We now assume that the symmetry hypothesis (H1') is satisfied and that 
the root $(\up, \underline \ct, \underline \cx$ is totally non glancing. 
 In \cite{MZ2}, symmetrizers for$ \chH_k (p, \cz, 0) $ 
 are constructed. We show that they also symmetrize 
 $\chH_k (p, \cz, \rho) $ when $\rho > 0$. 

 In \cite{MZ2}, it is proved that the nonglancing condition implies that 
 the multiplicity of $\mu_k$ as an eigenvalue of $H_0 (\up, \underline \cz)
 = \chH (\up, \underline \cz, 0)$ is equal to $m$.   
Denote by $V_k$ the $N \times m $ sub-matrix of $V$ which corresponds to the block
$\chH_k$. Therefore, for $(p, \cz, \rho)$ close to
$(\up, \uzeta, 0)$, the corresponding invariant space of 
$\chH_h $ is $\EE_k (p, \cz, \rho) = V_k(p, \cz, \rho) \CC^m$ and  
\begin{equation}
\label{Vk}
V_k \chH_k = \chH  V_k 
\end{equation}
Recall that $\EE_{k}^- (p, \cz, \rho)$ is the negative space of 
$\chH_k$ for $\cz \in \overline S^d_+$, $\rho \ge 0$ with $\cg > 0 + \rho > 0$. 

\begin{lem}
\label{lem5b1}
 
ii) If    $( \up, \utau, \uxi)$  is totally incoming, then,  for 
$(p, \zeta)$ in a neighborhood of $(\up, \uzeta)$, $\EE_k^-  (p, \zeta) = \CC^m $. 

iii) If    $( \up, \utau, \uxi)$  is totally outgoing, then,  for 
$(p, \zeta)$ in a neighborhood of $(\up, \uzeta)$, $\EE_k^- (p, \zeta) = \{0\}$.
\end{lem}
 
 \begin{proof}
 The dimension is constant for $\cg > 0 + \rho > 0$, and the result is proved in \cite{MZ2}
 when $\rho = 0$. 
 \end{proof}
      
      By assumption, there is a definite positive matrix $S(p)$ such that 
the $S A_j$ are symmetric.

\begin{lem}
\label{lem5b2}
The symmetric matrix  
\begin{equation}
\label{Sig}
\Sigma_{k,0}  (p, \zeta) =  - V_{k}^*(p, \zeta, 0)  S(p)  A_d(p)  V_{k}(p, \zeta,0)\,.  
\end{equation}
is a  symmetrizer for $\chH_{k}$ on a neighborhood of $(\up, \underline \cz, 0)$. 
 More precisely, there holds 
\begin{equation}
\label{eq59x}
\re \Sigma_{k} \chH_{k} = \gamma   R_1 + \rho R_2
\end{equation}  
with  $\Sigma_1(\up, \underline \cz, 0)$ and $\Sigma_2(\up, \underline \cz, 0)$ definite positive.

In addition, 
$\Sigma_k (\up, \uzeta, 0)  $  is definite positive [resp. negative]
 when the mode is totally incoming [resp. outgoing].  
 \end{lem}

\begin{proof}
According to \eqref{defcH}, there holds 
\begin{equation*}
\chH (p, \cz, \rho) = H_0(p, \cz) + \rho H' (p, \cz, rho). 
\end{equation*}
Using \eqref{Vk} and the definition \eqref{defH0} of $H_0$, one obtains 
the identity \eqref{eq59x} with 
\begin{eqnarray}
& & R_1  =   V_k^* S V_k ,    \\
&  & R_2   =  V_k^* (\re S A_d H' ) V_k. 
\end{eqnarray}
Because $S$ is definite positive, $R_1$ also has this property. 
Next,  Lemma~\ref{lem642b} implies that $H'(p, \cz, 0) = - H_1(p, \cz)$ 
with $H_1$ given by \eqref{eq644}. Since $H_0 (\up, \underline \cz) = \mu_k \Id = - i \xi_k \Id$ on 
$\EE_k (\up, \underline \cz, 0)$, there holds
\begin{equation*}
H'(\up, \underline \cz, 0) V_k (\up, \underline \cz, 0) =  - 
A_d^{-1}(\up) B(\up, \underline \cx) . 
\end{equation*}
Therefore, at the base point $(\up, \underline cz, 0)$, there holds 
\begin{equation*}
R_2 (\up, \underline cz, 0) = V_k^* ( \re S B ) V_k.  
\end{equation*}
The assumption (H1') implies that $SB$ is definite positive 
on the space $\EE_k (\underline p, \underline \cz, 0) = \ker (A(\up, \underline \cx) + \underline \ct \Id)$, implying that $R_2$ is definite positive at $(\up, \underline cz, 0)$, hence on a neighborhood
of that point. 

That $\Sigma_k (\up, \uzeta, 0)  $  is definite positive [resp. negative]
 when the mode is totally incoming [resp. outgoing] is proved in 
\cite{MZ2}.
\end{proof}

 With Lemma~\ref{lem5b1}, this implies that 
 \begin{equation}
\Sigma^\kappa _k = \left\{\begin{array}{rl}
   \Sigma_k  & \quad \mathrm{in \ the \ incoming \ case, } 
\\
\kappa \Sigma_k  & \quad \mathrm{in \ the \ outgoing \ case. } 
\end{array}
\right. 
\end{equation}
are K-familes of symmetrizers for $\chH_k$.


%
%

\section{Further remarks and examples}

\label{Sec6} 
\subsection{Adapted basis. Proof of Proposition~\ref{prop414}}
In this section, we always assume that  Assumptions (H1) is  satisfied. 
Consider   a geometrically regular root 
 $(\up, \underline \ct, \underline \cx)$  of $\Delta$. We show that there 
 are eigenbasis $\{e_j\}$ satisfying \eqref{23} which are adapted to $B$, in the sense 
 of Definition~\ref{defdec}, either when the multiplicity is 2 or when 
 the system is symmetric. 

\medbreak
  \noindent 1. {\sl  The case of multiplicity two}
  
Projecting on the 2-dimensional invariant space of 
$i A (p, \cx) +  \rho B (p, \cx) $ associated to the eigenvalues close the 
$i \lambda_1$ and $i \lambda_2$, we are reduced to consider 
$2 \times 2$ matrices 
\begin{equation}
\label{330}
i \tilde A (p, \cx) + \rho \tilde B(p, \cx, \rho)  \quad \mathrm{with} \quad 
\tilde A = \begin{pmatrix}
     \lambda_1 &    0 \\
    0  &  \lambda_2
\end{pmatrix} . 
\end{equation}
 The Assumption (H1) implies that the spectrum
 $ i \tilde A + \rho \tilde B$ is contained in $\re \lambda \ge c \rho$, 
 for $(p, \cx)$ close to $(\up, \underline \cx)$ and $\rho \in [0, \rho_0]$, 
 for some $\rho_0 > 0$. 
 Changing $\cx$ to $- \cx$ and using (H1) near $- \underline \cx$, we see that 
 the spectrum
 $ \pm  i \tilde A + \rho \tilde B$ is contained in $\re \lambda \ge c \rho$.  

 We show that, changing the 
base $\{ e_1, e_2\}$ if necessary,  one always 
meet condition \eqref{334x}. 

\begin{lem}
With assumptions as above, there is a smooth change of bases preserving 
$\eqref{330}$, such that $\re \tilde B(\up, \underline \cx)$ is definite positive. 
\end{lem}

\begin{proof} 
 
 The constant multiplicity case $\lambda_1 =\lambda_2$ being already 
 treated, we  assume that $\lambda_1 \ne \lambda_2$ on any neighborhood 
 of $(\up, \underline \cx)$. 
  In this case we are limited to consider 
 diagonal change of basis  and we prove  that there exists a  diagonal matrix 
 $D$, such that 
 \begin{equation}
\label{331}
\re \big ( D \tilde B(\up, \cx, 0) D^{-1}  \big) > 0 . 
\end{equation}

{\bf a) }  Recall that  there is $c > 0$ such that the spectrum
 $ \pm  i \tilde A + \rho \tilde B$ is contained in $\re \lambda \ge c \rho$.
 We first show that for all $t \in \RR $  the spectrum of 
\begin{equation}
\label{332}
\begin{pmatrix}
      0 &0    \\
      0 &i t   
\end{pmatrix}  +    \tilde B(\up, \underline \cx, 0)  - \frac{c}{4} \Id 
\end{equation}
is contained in $\re \lambda  > 0 $.    
If not, there are  $t $, $\rho_1 > 0$   and a neigborhood  $\omega$ 
of $(\up, \underline \cx)$ such that  
\begin{equation}
\label{333}
\begin{pmatrix}
      0 &0    \\
      0 &i t   
\end{pmatrix}  +    \tilde B(\up, \cx, \rho) 
\end{equation}
has an eigenvalue  in $\re \lambda <  c/  2 $ when $(p, \cx) \in \omega$
and $\rho \in [0, \rho_1]$.  
There is  $(p', \cx') \in \omega$ such that 
$ \lambda_2 (p', \cx') - \lambda_1(p', \cx') = t_1 \ne 0$. 
Choose $\rho \in [0, \rho_1[$ such that 
$ \rho |t | \le  |t_1|$. By continuity, since $\lambda_2 - \lambda_1$ vanishes 
at $(\up, \underline \cx)$, 
 there is  $(p, \cx) \in \omega$ such that 
$ \lambda_2 (p, \cx) - \lambda_1(p, \cx) = \pm t \rho$.  
Therefore the matrix 
$\pm i \tilde A(p, \cx) + \rho \tilde B(p, \cx, \rho)$ has an eigenvalue 
in $\{ \re \lambda \le \rho c/2 \}$. 

\medbreak 
{\bf b) } Consider a matrix
 \begin{equation*}
 \begin{pmatrix}
    a  &  b  \\
    c  &  d
\end{pmatrix}. 
\end{equation*}
Its spectrum is contained in $ \{ \re \lambda > 0 \}$, if and only if 
$$
\re (a + d ) > 0  \quad \mathrm{and} \quad  | \re (\sqrt f) |^2  \le \re (a+ d)^2 
$$
where $ f= (a-d)^2 + 4 bc$. Since $ | \re (\sqrt f) |^2  = \frac{1}{2} ( \vert f \vert + \re f ) $, 
the second condition reads
$$
| \im  f  |^2  <   4 (\re a+ \re d)^2  \big(  (\re a+ \re d)) ^2 - \re f \big) , 
$$ 
or 
$$
\begin{aligned}
\big( \re(a-d)&  \im (a-d) +  2 \im ( bc)  \big)^2  < 
\\
 & (\re a+ \re d)^2  \big( (\im a - \im d) ^2 + 4  \re a \re d -  4 \re (bc)  \big)
  \end{aligned}
$$
or 
\begin{equation}
\label{335}
\begin{aligned}
 \re a \re d   \big( (\im a - \im d) ^2 
 &- (\re a - \re d) \im (bc) (\im a - \im d) 
 \\
  +     (\re a + \re d) ^2 &( \re a \re d -    \re bc ) - | \im (bc)|^2 \big) > 0. 
  \end{aligned}
\end{equation}
We apply this criterion to the matrices   \eqref{333}. In this case, when $t$ varies
in $\RR$ the coefficient $\im (a - d)$ varies from $- \infty$ to $+ \infty$ while the other coefficients 
are fixed. Therefore, if the corresponding inequality \eqref{335}  is satisfied  for all
$t$,  then $\re a + \re b > 0$, $\re a \re d \ge  0$ and 
$$
\begin{aligned}
(\re a - \re d)^2 & | \im (bc) |^2  \le \\
& 4 \re a \re d  \big( 
 (\re a + \re d) ^2 &( \re a \re d -    \re bc ) - | \im (bc)| ^2 \big)  . 
\end{aligned}
$$
Thus 
$$
| \im (bc) |^2   \le  4 \re a \re d   ( \re a \re d -    \re bc ) 
$$
and 
\begin{equation*}
| bc |  + \re (bc) \le  2 \re a \re d . 
\end{equation*}
Denoting by $b_{j,k} $ the entries of 
$\tilde B (\up, \underline \cx, 0)$, we see that the spectral condition of step a) implies the 
following conditions: 
\begin{equation}
\label{336}
\re b_{1 1} > 0, \quad \re b_{2 2} > 0 , \quad   | b_{1 2} b_{2 1}  |  +
 \re (b_{12} b_{21} )  <   2 \re b_{11} \re b_{22} . 
\end{equation}

\medbreak
{\b c) } Similarly, we note that the condition  $\re \tilde B (\up, \underline \cx, 0) > 0$ 
is equivalent to 
\begin{equation}
\label{337}
\begin{aligned}
&\re b_{1 1} > 0, \quad \re b_{2 2} > 0 , \quad 
\\
&  | b_{1 2} |^2 + | b_{2 1}  |^2   +
2  \re (b_{12} b_{21} )  <   4 \re b_{11} \re b_{22} . 
\end{aligned}
\end{equation}
With $D = \begin{pmatrix}
    1  &  0  \\
     0  &  \delta 
\end{pmatrix}$, the conjugation $D \tilde B D^{-1}$ changes 
$B $ to $B'$ with 
$$
b'_{11}  = b_{11},  \quad b'_{22} = b_{22}, \quad  b'_{21} = \delta b_{21}, \quad 
b'_{12} = \frac{1}{\delta} b_{12}. 
$$
For all $\eps > 0$, one can choose $\delta$ such that 
$$
|b'_{12} |^2 + |b'_{21} |^2 \le | b_{12} b_{21} | + \eps
$$
and therefore, \eqref{336} implies that there  is  $\delta$  such that 
\eqref{337} holds for the $b'_{j k}$. 
\end{proof}

   \noindent 2. 
   \textit{Symmetric systems.} 
     
     \begin{lem}  
      Suppose  that $(\up, \underline \ct, \underline \cx)$  is  geometrically regular  and 
       Assumption (H1') holds.  Then one can choose the eigen-basis $\{ e_j\} $ such that 
       \eqref{331x} holds. 
     \end{lem} 
     
     \begin{proof} Denote by 
      $ S$ the symmetrizer. 
 We show that one can choose the eigen-basis $e_j$ such that 
 \begin{equation}
\label{ejorth}
{}^t e_j (\up, \underline \cx) \, S(\up)  \, e_{j'}  (\up, \underline \cx) = \delta_{j, j'}. 
\end{equation}
In this case, $\underline \ell_j = {}^t e_j( \up, \underline \cx) S(\up)$ 
and 
 \begin{equation}
 \label{bssym}
B^\sharp_{j, j'} = {}^t e_j( \up, \underline \cx) S(\up) B(\up, \underline \cx) e_{j'}( \up, \underline \cx) 
\end{equation}
showing that $\re B^\sharp$ is the restriction of 
$\re \big( S(\up) B(\up, \underline \cx)\big)$ to the space spanned by the 
$e_j( \up, \underline \cx) $, and hence positive. 

To prove \eqref{ejorth}, 
consider the partition of  $\{1, \ldots, m \}$ into susbets $J_a$ 
such that $j$ and $j'$ belong to the same class   $J_a$ if and only if 
$\lambda_j = \lambda_j' $ on a neighborhood of $(\up, \underline \cx)$. 
Denote by $\FF_a(p, \cx) $ the space spanned by the $e_j(p, \cx)$ for 
$j \in J_a$. Then, near $(\up, \underline \cx)$,  $A(p, \cx) = \tilde \lambda_a \Id $ on this space, where 
$\tilde \lambda_a$  is the common value of the $\lambda_j $ for $j \in J_a$. 
Thus, locally, one can find a smooth basis  of $\FF_a$, analytic in $\cx$ and 
orthonormal for the scalar product $S(p)$. 
Collecting these bases of  $\FF_a$, \eqref{ejorth}   holds when $j$ and $j'$ 
belong to the same class $J_a$. 

When $j$ and $j'$ do not belong to the same class $J_a$, 
there is a sequence $(p^n, \cx^n)$ converging to $(\up, \underline \cx)$ such that 
$\lambda_j(p^n, \cx^n) \ne \lambda_{j'} (p^n, \cx^n) $. 
The symmetry  of $ S(p^n) A(p^n, \cx^n)$ implies  that 
\begin{equation*}
 {}^t e_j (p^n, \cx^n) \, S(p^n)  \, e_{j'} (p^n, \cx^n) = 0 . 
\end{equation*}
Therefore, passing to the limit, we see that \eqref{ejorth} is also satisfied
when $j$ and $j'$ do not belong to the same class $J_a$. 
 \end{proof}

\subsection{ Discontinuity of the negative spaces $\EE^-$ } 
    
    We show the decoupling condition \eqref{331x} is necessary for the continuity of 
    $\EE^-(p, \cz, \rho)$ at $\rho = 0$. 
    Before stating  the result, we make the following remark.

     \begin{lem}
      \label{lem63}          
      Suppose that  $(\up, \underline \ct, \underline \cx)$  is  geometrically regular   and  nonglancing.
      With notations as in $\eqref{22}, \eqref{multig}, \eqref{bsharp}$ let 
      $\beta_j = \D_{\xi_d} \lambda_j( \up, \underline \cx) \ne 0$. Then, there is $c > 0$ such that for all 
      $\rho > 0$ and $t \in \RR$, the spectrum of 
      $i  t \mathrm{diag} (\beta_j)  + \rho B^\sharp$ is contained in 
      $\{ \re \mu \ge c \rho \}$. 
      \end{lem}
    
    \begin{proof} We fix $p = \up$ and forget it in the notations. 
    For $\xi$ close to $\underline \cx$, consider the invariant space of $i A( \xi) + \rho B(\xi) $ 
    associated to eigenvalues close to $- \underline \ct$. In the basis $\{e_j \}$, its matrix
    is 
    \begin{equation}
\label{610}
i \mathrm{diag} \big(\lambda_j(\xi ) \big)  + \rho \tilde B (\xi, \rho), 
\end{equation}
with $\tilde B(\underline \cx, 0) = B^\sharp$. 
The Assumption (H1) implies that the spectrum of this matrix lies in $\{ \re \mu \ge c \rho \}$. 

Adding $ \underline \ct \Id$, we can assume, without loss of generality, that 
$\lambda_j (\underline \cx) = 0$. Taking  $t > 0$, 
$\xi = \underline \cx +\pm (0, t )$ and   $\rho = t \sigma > 0 $  the matrix in 
 in \eqref{610} is 
 \begin{equation*}
 t  M(t, \sigma) = t \big( \pm  i  \mathrm{diag} (\beta_j)  + \sigma B^\sharp + O(t) \big) 
\end{equation*}
 and the spectrum of $M(t, \sigma) $ lies in $\{ \re \mu \ge c \sigma \}$. 
Letting $t $ tend to zero, implies that the 
spectrum of $M(0, \sigma)$ is also contained in $\{ \re \mu \ge c \sigma \}$
and the lemma follows by homogeneity.
    \end{proof}
    
    \begin{cor}
    \label{cor64}
 If  $(\up, \underline \ct, \underline \cx)$  is  geometrically regular   and  nonglancing, 
 then for all $\gamma \ge 0$ and $\rho \ge 0$, with $\gamma + \rho > 0$, 
 the matrix $ \diag(\beta_j^{-1}) \big( \gamma \Id + \rho B^\sharp)$  has no eigenvalues
 on the purely imaginary axis. 
    \end{cor}

    Consider  $\underline \cz = (\underline \ct, \underline \ch, 0) \ne 0$ and a purely imaginary eigenvalue 
  $\underline \mu_k = i \underline \cx_d$ of  $H_0(\up, \underline \cz)$. 
Let $\underline \cx = ( \underline \ch, \underline \cx_d)$. Then 
$(\up, \underline \ct, \underline \cx) $ is a root of $\Delta$.  
We denote by   $\chH_k$ the block associated to $\underline \mu_k$ and, 
for $  \rho > 0$, we denote by
 $\EE_k^-(p, \cz, \rho) $   the negative invariant space  of $\chH_k$.        
      
      \begin{prop}
      \label{prop65}    Suppose that  $(\up, \underline \ct, \underline \cx)$  is  geometrically regular   and  nonglancing and suppose that there exist $j \in J_I$ and $j' \in J_0$ such that 
\begin{equation}
\label{coup}
B^\sharp_{j', j} \ne 0.
\end{equation}
      Then the negative space $\EE_k^- (\up, \cz, \rho)$  has no limit as
      $(\cz, \rho) \to (\underline \cz, 0)$. 
     
     In particular, there are no smooth $K$-families of symmetrizers for 
     $\chH_k$ near $(\up, \underline \cz)$. 
       \end{prop}

      \begin{proof}
   
      By Lemma~\ref{lem31} and Lemma~\ref{lem37}, the block decomposition 
      \eqref{bdiag2} implies that in a suitable basis 
      \begin{equation}
      \label{611}
 \chH_k (\up, \underline \ct, \underline \ch, \gamma, \rho) =  - 
  \diag(\beta_j^{-1}) \big( \gamma \Id + \rho B^\sharp) + O( \gamma^2 + \rho^2). 
\end{equation}
Denote by $\EE^-(\gamma, \rho)$ the negative space 
of $ \chH_k (\up, \underline \ct, \underline \ch, \gamma, \rho)$.  We show that 
\begin{equation}
\label{612}
\lim_{\gamma \to 0} \EE^-(\gamma, 0) \ne \lim_{\rho \to 0} \EE^- (0, \rho), 
\end{equation}
  which implies that $\EE^-(\gamma, \rho)$ has no limit as 
  $(\gamma, \rho) \to (0, 0)$.  

Consider first the case where $\rho = 0$. Then, \eqref{611} implies that 
the first limit in \eqref{612} is the space 
$\EE_I$ spanned by the vectors $e_j$ of the basis such that $\beta_j > 0$, that is
such that $j \in J_I$. 

On the other hand, Corollary~\ref{cor64}  implies that 
$ B^\flat = -  \diag(\beta_j^{-1})   B^\sharp $ has no eigenvalue on the imaginary axis.
Therefore, the second limit in \eqref{612} is the negative space $\EE^-_{B^\flat}$ of $B^\flat$. 
If it were equal to $\EE_I$, this would mean that $\EE_I$ is invariant by $B^\flat$, 
thus by $B^\sharp = - \mathrm{diag}(\beta_j) B^\flat$, which contradicts \eqref{coup}. 

By \cite{MZ3}, the existence of smooth $K$-families of symmetrizers implies that 
the limit of $\EE^-_k$ at $(\underline \cz, \rho)$ exists, and is equal to the space 
$\underline \EE_k^-$ of Definition~\ref{def316}. Therefore, \eqref{612} 
implies that there are no smooth $K$-families of symmetrizers. 
\end{proof}

      
 \subsection{Viscous instabilities}

      Consider boundary conditions as in Assumption~\ref{assbc}. 
      When the negative space $\EE^-$ is not continuous in $(\cz, \rho)$, then
      the Evans function is likely not continuous and one can 
      expect that the low-frequency uniform  stability condition for the viscous problem is strictly 
      stronger than the similar condition for the inviscid problem. 
      In particular, the inviscid problem can be strongly stable while 
      the viscous one is strongly unstable. 
    We illustrate  here this phenomenon on 
     an  explicit  example.       
      
      \medbreak
 \noindent
1. \emph{An example.}      
      Consider the system
      \begin{equation}
\label{exeq}
\left\{\begin{aligned}
&(\D_t + \D_y ) u_1 + \D_x u_2 =   \eps \mu \Delta u_1, \\
&(\D_t + \D_y ) u_2 + \D_x u_1 = \eps \nu \Delta u_2. 
\end{aligned}
\right.
\end{equation}
      Taking linear combinations and changing $\eps$, the system is equivalent to 
         \begin{equation}
\label{exeq2}
 (\D_t + \D_y) \Id +  A \D_x  -  \eps B \Delta , \quad  A = \begin{pmatrix}
    1  &0    \\
      0& -1  
\end{pmatrix}, \quad B = \begin{pmatrix}
1 & a\\
a & 1
\end{pmatrix},
\end{equation}
with $a = | \nu - \mu |/ (\nu + \mu) \in [0, 1 [$. 
This system is symmetric and satisfy the assumptions (H1) and (H2). 

The hyperbolic part is diagonal: the eigenvalues are 
\begin{equation}
\label{ex3}
\lambda_1 = \eta  + \xi,  \quad \lambda_2 = \eta - \xi. 
\end{equation}
     They cross on the line  $\xi = 0$  and are trivially geometrically regular since the system is already
     in diagonal form. One of the eigenvalue is incoming, one is outgoing. 
      The decoupling condition \eqref{331x} is satisfied if and only if $a = 0$. 
      In the sequel, we assume that $a > 0$.
       
\bigbreak
   \noindent 2. \emph{Boundary conditions.} 
   Next,   consider boundary conditions for \eqref{exeq2}: 
    \begin{equation}
\label{bce}
u_{| x = 0}   +  \eps \Gamma \D_x u  _{| x = 0}  = 0 . 
\end{equation}
   We first compute the  limiting inviscid boundary conditions, using boundary layers. 
   The bounded solutions $u = w(x/ \eps)$ of \eqref{exeq2} are 
   \begin{equation}
\label{profex}
w(z) =  u  + e^{ z B^{-1} A} h, \quad   h \in \EE^-_{B^{-1} A},  \ u \in \CC^2. 
\end{equation}
where $\EE^-_{B^{-1} A}$ is the negative space 
of $B^{-1} A$. Therefore, $ u $ is the endpoint of a profile which satisfies the boundary 
condition \eqref{bce}, if and only if 
\begin{equation}
\label{hybxex}
u \in \big(  \Id + \Gamma    B^{-1} A \big) \EE^-_{B^{-1} A}. 
\end{equation}
Note that given any complex number $\underline c$, one can choose  $\Gamma$
such that this boundary condition reads
\begin{equation}
\label{bcex}
u_1 = \underline c u_2 
\end{equation}

\noindent 3. \emph{Low frequency stability.} 
The first order system \eqref{linq4} reads
\begin{equation}
\label{exy1}
\D_z U - G(\zeta) U , \quad  G (\zeta) = \begin{pmatrix}
   0   &   \Id  \\
   \sigma B^{-1} +  \eta^2 \Id &    B^{-1} A    
\end{pmatrix}, 
\end{equation}
with $\zeta = (\tau, \eta, \gamma)$  and $ \sigma = \gamma + i (\tau + \eta)$. 
 Perform the small frequency reduction \eqref{block1}, 
 using the change of unknows
  \begin{equation*}
 \begin{pmatrix}
    u \\  \D_z u  \end{pmatrix} = V(\zeta) 
\begin{pmatrix}
       u_H \\ u_P 
\end{pmatrix}. 
\end{equation*} 
Then, by Lemma~\ref{lem642b}, there holds 
\begin{equation*}
V^{-1} G V = \begin{pmatrix}
     H &0    \\
     0 &  P
\end{pmatrix}
\end{equation*}
with $P(0) = B^{-1} A$ and 
 \begin{equation}
\label{ex4}
H (\zeta) = -  \sigma A  +  ( \sigma^2 - \eta^2) AB  + O( | \zeta |^3) , 
\quad \end{equation}

Since $V(0)$ has the triangular form \eqref{Vtri}, we see that the 
boundary condition reads
\begin{equation}
\label{nbc}
u_H  +  \tilde \Gamma (\zeta) u_P = 0, \quad  \tilde \Gamma(0 )  = \Gamma + A^{-1} B . 
\end{equation}
The Evans condition is violated at $\zeta$ if there is $u_H \in \EE^-_H(\zeta)$
and $u_P \in \EE^-_P (\zeta)$  satisfying this  boundary condition. 
 The negative space of $P(\zeta)$,  $\EE^-_P(\zeta)$ is smooth in $\zeta$ and 
 equal to $\EE^-_{B^{-1} A}$ when $\zeta =0$. Thus, the   Evans condition is violated at $\zeta$ 
 if and only if if there is $u_H \in \EE^-_H(\zeta)$ such that 
 \begin{equation*}
\EE^-_H(\zeta) \cap    \tilde \Gamma (\zeta) \EE^-_P(\zeta)  \ne  \{ 0 \}   
\end{equation*}
Since $A^{-1} B = (B^{-1} A )^{-1}$, there holds
$$
\tilde \Gamma(0) \EE^-_P(0) = ( \Id + \Gamma B^{-1} A) \EE^-_{B^{-1} A}. 
$$
Comparing with \eqref{hybxex} and \eqref{bcex}, we see that for $\zeta$ small, 
the space    $ \tilde \Gamma (\zeta) \EE^-_P(\zeta) $ is generated 
by ${}^t (c(\zeta), 1)$ where 
  $c(\zeta)$ is a smooth function  such that 
$c(0) = \underline c$. Therefore, the  Evans condition is violated at $\zeta$ if and only if 
\begin{equation}
\label{nEC}
\begin{pmatrix}
      c(\zeta)     \\
      1  
\end{pmatrix}  \in \EE_H^-(\zeta) . 
\end{equation}

\begin{rem}
\textup{Using the terminology of \cite{MZ2}, the analysis above 
shows that the \emph{reduced boundary condition}  for the hyperbolic part 
$H(\zeta)$ reads}
\begin{equation}
\label{bcxy}
u_1 = c(\zeta) u_2. 
\end{equation}
Taking $\zeta = 0$ in this equation, we recover that \eqref{bcex} is the natural limiting 
boundary condition for the hyperbolic operator $H_0$.
\end{rem}

\begin{prop}
There are choices of $a$ and $\Gamma$, such that

i)  the inviscid problem
$\eqref{exeq2}$ for $\eps = 0$ with the boundary condition $\eqref{bcex}$
is maximal striclty dissipative thus uniformly stable, 

ii) the viscous problem with boundary conditions $\eqref{bce} $ is strongly unstable
for small frequencies, in the sense that the Evans functions vanishes 
for arbitrarily small frequencies $\zeta$ with $\gamma > 0$. 

\end{prop}
  
 \begin{proof} The matrix 
  \begin{equation}\label{symm}
S=
\begin{pmatrix}
1 & 0\\
0 & s
\end{pmatrix}, \quad  s > 0
\end{equation}
 is a symmetrizer for the inviscid problem. 
If 
 \begin{equation}
 \vert \underline c \vert ^2  < s, 
 \end{equation} 
  the boundary condition 
 is  strictly dissipative for  $S$. This implies that the uniform Lopatinski condition 
 is satisfied. 
  
   We consider frequencies $\zeta = \rho \cz$ with 
   $\cz$ close to $(-1, 1, 0)$  where $H_0(\cz) = 0$ has a double eigenvalue. 
   More precisely we consider frequencies 
   \begin{equation}
\label{choice}
\zeta = ( - \rho + \rho^2  \hat \tau, \rho, \rho^2 \hat \gamma ) . 
\end{equation}
  In this case, we see that $G$ is a function of $\hat \sigma = \hat \gamma + i \hat \tau$
  and $\rho$, holomorphic in $\hat \sigma$, as well as 
  $V$, $P$, $H$ and $c$. Moreover
  \begin{equation}
\label{redsy}
H(\zeta) = - \rho^2 ( \hat \sigma A + A^{-1} B + O(\rho) ) = \rho^2 \hat H( \hat \sigma, \rho). 
\end{equation}
   The model operator is 
   $$
   \hat H(\hat \sigma, 0) = -  \hat \sigma A -  A^{-1} B = 
   \begin{pmatrix}
  -  \hat \sigma -  1    &   - a   \\
   a     &    \hat \sigma + 1
\end{pmatrix}
   $$ 
  $\hat H(1, 0)$ has one eigenvalue with positive real part, with 
  eigenvector ${}^t( \underline b, 1)$ with 
  $\underline b = ( 2 + \sqrt {4 - a^2}) / a$ (Note here the importance of the assumption $a \ne 0$). 
  Therefore, for 
 $\hat \sigma $ close to $1$ and $\rho$ small , the negative space of  $\hat H(\hat \sigma, \rho )$
 is generated by 
 ${}^t ( b  (\hat \sigma, \rho) , 1)$ where $b$ is smooth and holomorphic in $\hat \sigma$
 and $b(1, 0) = \underline b$. Moreover
 \begin{equation}
\label{impl}
\D_{\hat \sigma} b (1, 0) = \frac{1}{a} \big( 1 + \frac{2}{\sqrt{4 - a^2}} \big) \ne 0. 
\end{equation}
   Comparing with \eqref{nEC}, we see that the stability condition is violated 
   at   $\zeta$ given by \eqref{choice}, if and only if 
   \begin{equation}
\label{nECr}
b(\hat \sigma, \rho) = c(\zeta) = \hat c ( \hat \sigma, \rho). 
\end{equation}
 
 Given $ a \in ]0, 1[$, we choose $\underline c = \underline b$ and 
 $\Gamma $ such that the inviscid boundary condition reads \eqref{bcex}. 
 Note that $\hat c ( \hat \sigma, 0) = \underline c$ for all $\hat \sigma$. 
  Thus the equation \eqref{nECr} holds at 
 $\hat \sigma = 1 $ and $\rho = 0$. Moreover, with \eqref{impl}, the 
 implicit function theorem shows that for $\rho > 0 $ small, there is 
 $\hat \sigma (\rho$ close to $1$ solution of \eqref{nECr}, 
 providing frequencies $\zeta (\rho ) = O(\rho)$ 
 with $\gamma(\rho) \sim \rho^2   > 0$, where 
 the stability condition is violated.
    \end{proof}
 
 \bigbreak
 
 \noindent 4. \emph{Smooth symmetrizers. }
 We briefly discuss here the existence of smooth symmetrizers for 
 the hyperbolic operator $\chH$ \eqref{defcH}. In the present case, we deduce from 
 \eqref{ex4} that 
      in polar coordinates $\zeta = \rho \cz$,   there holds 
       \begin{equation}
\label{ex5}
\chH (\cz, \rho) = -  \check \sigma A  +  \rho  ( \cs^2 - \ch^2) AB  + O( \rho^2) , 
\quad  \cs = \cg + i (\ct+ \ch).
\end{equation}
      Fix  $\underline \cz = (1, -1, 0)$, which corresponds to a multiple root
      of the hyperbolic part. Then $\underline \cs = 0$, and near 
    $(  \underline \cz, 0) $ 
           \begin{equation}
\label{ex6}
\chH (\cz, \rho) = -   A   (\cs \Id + \rho  \beta (\cz) B )  + O( \rho^2)  
\end{equation}
with $\beta(\underline \cz) = 1$. 
      Dropping the $\check{} $ , and changing $\rho b $ to $\rho$, the matrix
      $\chH$ is  a perturbation for 
      $(\sigma, \rho)$ close to $(0, 0)$ of the 
      following  \emph{canonical example}
    \begin{equation}\label{eg}
\begin{pmatrix}
1 & 0\\
0 & -1
\end{pmatrix}
\partial_x
+ \sigma 
\begin{pmatrix}
1 & 0\\
0 & 1
\end{pmatrix}
+\rho
\begin{pmatrix}
1 & a\\
a & 1
\end{pmatrix}, \quad \re  \sigma \ge 0, \ \rho \ge 0, 
\end{equation}
  Note that \eqref{redsy} derives from \eqref{ex6} choosing 
  $\cs = \rho \hat \sigma$.

  Denote by $\EE^- $ the negative space of $\chH$ for 
  $\re \sigma + \rho > 0$.  On can check directly on this example that the 
negative spaces  have no limit  as $(\sigma, \rho) \to (0, 0)$:  the limits are  different 
when $\rho = 0$ and $\sigma = 0$, since the positive  spaces of $A$ and $AB$ are different 
when $ a \ne 0$. 

On the other hand, blowing up once more the local coordinates near $\underline \cz$, 
that is taking polar coordinates $(\sigma, \rho) = r  (\hat \sigma, \hat \rho) $, 
is is clear from \eqref{ex6} that $\EE^-$ is a smooth function of 
$ (\hat \sigma, \hat \rho) $. 
   
\bigbreak
  
  If $\Sigma(\cz, \rho)$ is a smooth symmetrizer for 
   $\chH$, then  \eqref{global41} implies that 
   $\underline \Sigma = \Sigma( \underline \cz, 0) $ must be a symmetrizer 
   for $- ( \sigma A + \rho A B) $ for all $\sigma$ and $\rho$, equivalently 
   that $S = \underline \Sigma A $ is a symmetrizer for \eqref{eg}, that is 
   \begin{equation}
\label{symex}
 S = S^*  \gg 0 , \quad  SA = AS  ,  \quad  \re (S B) \gg 0. 
\end{equation}
The first two conditions  are satisfied if and only if $S$ is diagonal  and positive. 
Multiplying it by a positive factor, it must be of the form \eqref{symm}. 
  
 The third condition holds if and only if
\begin{equation*}
s  > a^2 ( 1+ s)^2/ 4 .  
\end{equation*}
Denoting by $ s_{\mathrm{min}} (a) <  1 < s_{\mathrm{max}} (a) < \infty  $  
the roots of the equation $ 4 s = a^2 (1 + s^2)$, the condition reads 
\begin{equation}
\label{CNs} 
s_{\mathrm{min}} (a) <  s < s_{\mathrm{max}} (a) . 
\end{equation} 
 This shows that the choice of symmetrizers  is much more limited
 in the viscous case compared to the inviscid one. 
 In particular, when $a$ is close to 1,   \eqref{CNs}  forces to choose $s$
  in a small interval around $1$. 
  \medbreak

 The boundary condition \eqref{bcxy} is strictly dissipative for $\Sigma$, then 
 \eqref{bcex} is strictly dissipative for $\underline \Sigma$. This holds if and only if 
 $s > | \underline c|^2$. Therefore:  

 \emph{  
  There is a smooth  symmetrizer $\Sigma(\cz, \rho)$ for $\chH$
  on a neighborhood of  $(\underline \cz, 0)$, adapated to the boundary conditions
  \eqref{bcxy} 
 only if 
 \begin{equation}
 \label{maxdiss}
| \underline c |^2  < s_{\mathrm{max}} (a) . 
\end{equation} 
 }

%
%
%

\section{The high-frequency analysis}
\label{S4} 

\subsection{The main high-frequency estimate}

This section is devoted to an analysis of  uniform  maximal estimates
 for high-frequencies.  We still assume that the Assumptions of Section 2 are 
 satisfied and we prove that the anticipated \eqref{maxesthf} are satisfied when the uniform spectral stability conditions are satisfied,  under 
 the following additional  structural  assumptions which strengthens  (H3):  it means first that the block 
 $L^{11}$ is hyperbolic with constant multiplicity with respect to time, and second that it is totally incoming our outgoing.  
 
 \begin{ass}
\textup{(H5) } For all $u \in \cU^*$ and $\xi \in \RR^d\backslash\{0\}$, the eigenvalues of $\overline A^{11}(u, \xi) $  are real,  semi-simple and have  constant multiplicities  . 

\textup{(H6)}  $L^{11}(u, \D) $ is also hyperbolic with respect to the normal direction $dx_d$. 
 
 \end{ass}

For Navier-Stokes and MHD equations and in many examples
$L^{11}$ is a transport field
\begin{equation}
L^{11}   =  \D_t + \sum_{j=1}^d a_j(u) \D_j 
\end{equation} 
and the condition reduces to $ a_d (u) \ne 0$  for $u \in \cU^*$, that is to 
Assumption~\ref{assnoncar2}, which means inflow or outflow boundary conditions.  
The
hyperbolicity condition (H6) in the normal direction  is   important as shown on an example below.
On the other hand the constant multiplicity  condition (H5) is more technical, and could be replaced 
by symmetry conditions: this is   briefly discussed in Remark~ \ref{remH5}.  

We consider the linearized equation \eqref{linq4}: 
\begin{equation}
\label{linq4nn}
\D_z  u  = \cG(z, \zeta) u +  f,   \quad   \Gamma (\zeta) u(0) = g 
\end{equation}
with  $u = {}^t (u^1, u^2, u^3) $,  $f = {}^t (f^1, f^2, f^3) $, $ \Gamma$ as in 
\eqref{linviscbcd} and $g = {}^t (g^1, g^2, g^3) $. 
\begin{theo}
\label{tm7n1} 
With assumptions as indicated above, assume that the 
uniform spectral stability condition is satisfied for high frequencies. 
Then there are $\rho_1 > 0$ and $C$ such that for all 
$\zeta \in \overline \RR^{d+1}_+$ with $| \zeta | \ge \rho_1$, 
the solutions of $\eqref{linq4nn}$ satisfy 
\begin{equation}
\label{mainestHFss} 
\begin{aligned}
(1+ \gamma) \| u^1 \|_{L^2} + \Lambda & \| u^2 \|_{L^2}  + \| u^3 \|_{L^2} 
\\
+ (1+ \gamma)^\mez | u^1(0)| & + \Lambda^\mez |  u^2 (0) |  +
\Lambda^{-\mez}  |  u^3(0) | 
\\ 
\le 
 C   \big( & \| f^1 \|_{L^2} +   \| f^2 \|_{L^2}  + \Lambda^{-1} \| f^3 \|_{L^2} \big) 
  \\
& + C \big(   (1+ \gamma)^\mez | g^1 |  + \Lambda^\mez |  g^2   |  +
\Lambda^{-\mez}  |  g^3 | \big). 
\end{aligned}
\end{equation}
\end{theo}

High frequencies require a particular analysis for two reasons. First, the splitting
hyperbolic vs parabolic is quite different in this regime and second
  the conjugation operator $\Phi$ of Lemma~\ref{lemconj} is not uniform for large 
$\zeta$.  The analysis is made in \cite{MZ1} for full viscosities and Dirichlet boundary conditions. 
For partial viscosities and shocks, that is for transmission condition, the problem is solved  in 
\cite{GMWZ4}.  The presentation below is more systematic  and allows for more general 
boundary conditions of the form \eqref{viscbcd}. 

We now explain the general strategy of the proof. We use the notations
\begin{equation}
\label{nota74}
\begin{aligned}
 & \| u\|_{sc} = (1+ \gamma) \| u^1 \|_{L^2} + \Lambda   \| u^2 \|_{L^2}  + \| u^3 \|_{L^2} , \\
 & \| f \|'_{sc} =  \| f^1 \|_{L^2} +   \| f^2 \|_{L^2}  + \Lambda^{-1} \| f^3 \|_{L^2},
\\
 & | u(0) |_{sc} = (1+ \gamma)^\mez | u^1(0)|   + \Lambda^\mez |  u^2 (0) |  +
\Lambda^{-\mez}  |  u^3(0) | ,
\\
 & | g |_{sc} =  (1+ \gamma)^\mez | g^1 |  + \Lambda^\mez |  g^2   |  +
\Lambda^{-\mez}  |  g^3 |. 
\end{aligned}
\end{equation}

{\bf 1) }  The main step in the proof of the theorem is to separate off  the incoming and outgoing components of $u$. This is done using  a change of variables 
$ \hat u = \cV^{-1} (z, \zeta)  u $ which transforms 
the equation \eqref{linq4nn} to 
\begin{equation}
\label{linq4nnn}
\D_z  \hat u  = \hat \cG(z, \zeta) \hat u +  \hat  f,   \quad   \hat \Gamma (\zeta) \hat u (0) = g . 
\end{equation}
There are norms similar to \eqref{nota74} for $\hat u$ and $\hat f $ as well; with little 
risk of confusion, we use here the same notations. An important property 
is that:  
\begin{equation}
\label{inva75}
\begin{aligned}
 & \| u\|_{sc}  \le C  \| \hat u \|_{sc}, \quad \quad \ \ 
   \| \hat f \|'_{sc}  \le C \| f \|'_{sc} ,
   \\
  &  | u(0) |_{sc} \le C | \hat u(0) |_{sc},     \quad 
  | \hat u(0) |_{sc} \le C |   u(0) |_{sc},   
  \end{aligned}
\end{equation}
with $C$ independent of $\zeta$.  Moreover, $\hat \Gamma (\zeta) =\Gamma(\zeta) 
\cV (0, \zeta)$ satisfies 
\begin{equation}
\label{inva76}
| \hat \Gamma(\zeta)  \hat u (0) |_{sc}  \le   C | \hat u(0) |_{sc}. 
\end{equation}

The new matrix $\hat \cG$ has the important property that 
 \begin{equation}
 \hat \cG = \begin{pmatrix} \hat \cG^+ & 0 \\ 0 & \hat \cG^- \end{pmatrix}  + 
\hat  \cG ' 
 \end{equation}
 with 
 \begin{equation} 
 \label{err79nn}
 \| \hat \cG' \hat u \|'_{sc}  \le  \eps(\zeta) \| \hat u \|_{sc}  
 \end{equation}
 where $\eps(\zeta)$ tends to $0$ as $ |\zeta |$ tends to infinity.  
The block structure corresponds to a splitting 
$
\hat u = ( \hat u^+  ,  \hat u^-)$
   with 
  $\hat u^- \in \CC^{N_b  } $ and $\hat u^+ \in \CC^{N+N'- N_b}$ denoting the incoming and outgoing components respectively.

\medbreak

{\bf 2) }  One proves separate estimates for the incoming and outgoing components: 
\begin{eqnarray}
\label{estu+}
&&\| \hat u^+ \|_{sc}   + | \hat u^+(0) |  \le C \| (\D_z - \hat \cG^+) \hat u^+ \|_{sc} , 
\\
&&\label{estu-}
\| \hat u^- \|_{sc}     \le C \| \D_z - \hat \cG^-) \hat u^- \|_{sc}  +  C  | \hat u^-(0) |,  
\end{eqnarray}
with $C$ independent of $\zeta$.  (The norms are defined, identifying 
$\hat u^- \in \CC^{N_b} $ to $ (0, \hat u^-) \in \CC^N$ etc).
As a result, with \eqref{err79nn}, this implies that if $\hat u$ is a solution of 
\eqref{linq4nnn}, then 
 \begin{eqnarray}
\label{estu+b}
&&\| \hat u^+ \|_{sc}   + | \hat u^+(0) |  \le C \| \hat f  \|_{sc} + \eps (\zeta) \| \hat u \|_{sc} , 
\\
&&\label{estu-b}
\| \hat u^- \|_{sc}     \le C \| \hat f  \|_{sc} + \eps (\zeta) \| \hat u \|_{sc}    +  C  | \hat u^-(0) |,  
\end{eqnarray}

\medbreak 

{\bf 3) }  We show that the estimates above imply that if the uniform spectral 
stability condition is satisfied, then the solutions of 
\eqref{linq4nnn} satisfy for $| \zeta |$ large enough
\begin{equation}
\label{estHFsss1}
\| \hat u \|_{sc} + | \hat u(0) |_{sc} \le C \big( \| \hat f \|_{sc} + | g |_{sc} \big) 
\end{equation}
 implying that the solutions of \eqref{linq4nn} satisfy 
 \begin{equation}
\label{estHFsss2}
\|   u \|_{sc} + |   u(0) |_{sc} \le C \big( \|   f \|_{sc} + | g |_{sc} \big) 
\end{equation}
 that is \eqref{mainestHFss}. 
 
 $\bullet$ Indeed,  by definition,  $h \in \EE^- (\zeta)$ if and only if there is $u$ solution 
 of $ \D_z u = \cG u$ with $u(0) = h$. The corresponding $\hat u = \cV^{-1} u$ 
 satisfies  by \eqref{estu-b} 
 $$
   \| \hat u^-  \|_{sc} \le     C | u^-(0) | +  \eps(\zeta)  \| \hat u^+   \|_{sc}  
 $$
 if $\zeta$ is large enough. Therefore, \eqref{estu+b} implies that 
 for $\zeta$ large and   all  $h \in \EE^- (\zeta)$, 
$  \hat h  = \cV^{-1} (0, \zeta) h = (\hat h^+, \hat h^-)$ satisfies
 \begin{equation}
 \label{est156}
| \hat h^+ |_{sc} \le \eps (\zeta) | \hat h^- |_{sc} . 
 \end{equation} 
  
 $\bullet$ In addition $\hat \EE^-(\zeta) :=  \cV^{-1} (0, \zeta)\EE^- (\zeta) $ has dimension 
 equal to $N_b$, as the space of the $\hat h^-$. Therefore, \eqref{est156} shows that 
 for $\zeta$ large, the projection $h \mapsto h^-$ is bijective from 
 $\hat \EE^-(\zeta) $  to $\CC^{N_b}$, with inverse uniformly bounded in the norm 
 $| \cdot |_{sc}$. 
 
 The uniform spectral stability condition reads 
 \begin{equation}
 \forall h \in   \EE^-(\zeta), \qquad  | h |_{sc}  \le  C  | \Gamma (\zeta) h |_{sc} 
 \end{equation}
 (see \eqref{stabcondhf}). Using \eqref{inva75}, this implies 
 \begin{equation}
 \forall \hat h \in   \hat  \EE^-(\zeta), \qquad  | \hat h |_{sc}  \le  C  | \hat \Gamma (\zeta) \hat h |_{sc} . 
 \end{equation}
 Using the isomorphism between  $\hat \EE^-(\zeta) $  and  $\CC^{N_b}$, we see that for $\zeta$ large enough and  $\hat h^- \in \CC^{N_b}$, there is $\hat h^+$ such that $(\hat h^+, \hat h^-) \in \hat \EE^-(\zeta) $.  Together with
 \eqref{est156} and \eqref{inva76}, there holds
 $$
 | \hat h^- |_{sc} \le  | \hat h |_{sc}  \le  C  | \hat \Gamma (\zeta) \hat h |_{sc} 
 \le  C  | \hat \Gamma (\zeta) (0, \hat h^-)  |_{sc} + \eps(\zeta)  | \hat h^- |_{sc} . 
 $$
 For $\zeta$ large, the last term can be dropped, increasing $C$. Finally, we conclude that for all $\hat h \in \CC^N$
 \begin{equation}
 \label{est167} 
  | \hat h |_{sc}  \le  C  | \hat \Gamma (\zeta) \hat h |_{sc}   +  C | \hat h^+ |_{sc}.  
\end{equation}
 Applying   this estimate to $\hat u (0)$, combining with \eqref{estu+} and \eqref{estu-} 
 and absorbing the error term $\hat \cG' \hat u$ for $\zeta $ large, we immediately obtain
 \eqref{estHFsss1}.

 \bigbreak
 
 The third part of the proof will not be repeated. We will focus on the reduction \eqref{linq4nnn}
 and on the proof of the estimates for $\hat u^\pm$.     
 

\subsection{Spectral analysis of the symbol}

Consider the linearized operator \eqref{linL}
$$
- \cB \D_z^2  + \cA \D_z  + \cM . 
$$
The   coefficients  satisfy 
\begin{equation}
\label{coefflin}
\begin{aligned}
\cB (z) &= B_{dd} (w(z))\\
  \cA (z, \zeta) &= A_d(w(z) ) -  \sum_{j=1}^{d-1} i\eta_j \big( B_{j d} +  B_{d, j} \big)(w(z))   + 
   E_d (z)  
   \\
   \cM(z, \zeta) &=  (i \tau + \gamma) A_0(w(z) ) +  \sum_{j=1}^{d-1} i  
    \eta_j  \big( A_j(w(z)) + E_j(z) \big) 
   \\
   & \qquad   \qquad \qquad \qquad   \quad +  \sum_{j, k = 1}^{d-1} \eta_j \eta_k B_{j, k}(w(z) )   +   E_0 (z) 
   \end{aligned}
\end{equation}
where the $E_k $ are functions, independent of $\zeta$, which involve derivatives 
of $w$ and thus converge to $0$ at an exponential rate when $z$ tends to infinity. 
Moreover, we note that 
\begin{equation}
\label{Einfty}
 E^{11}_k = 0, \quad     E^{12}_k = 0 \quad   \mathrm{ for }\  k > 0. 
\end{equation}
With \eqref{struc1}, we also remark that $\cM^{12}$ does not depend on $\tau$ and $\gamma$ 

\medbreak

We start with a spectral analysis of the matrix $\cG $ in 
\eqref{linq4}.  It is convenient to use here the notations
$u = (u^1, u^2, u^3) \in \CC^{N-N'} \times \CC^{N'} \times \CC^{N'}$.  
In   the corresponding block decomposition of matrices and using  the notations above, there holds
\begin{equation}
\label{deccG}
\cG  = \begin{pmatrix}
     \cG^{11} & \cG^{12}  & \cG^{13}    \\  0 & 0&  \Id \\
      \cG^{31} & \cG^{32}  & \cG^{33}  
\end{pmatrix}
\end{equation}
where    
\begin{eqnarray*}
\cG^{11} = -  (\cA^{11})^{-1} \cM^{11}, && 
\cG^{31} =   (\cB^{22})^{-1} ( \cA^{21}\cG^{11} + \cM^{21})  , 
\\  \cG^{12} = -  (\cA^{11})^{-1} \cM^{12}, & 
 &  
\cG^{32} =   (\cB^{22})^{-1} ( \cA^{21}\cG^{12} + \cM^{22}), 
\\
 \cG^{13} = -  (\cA^{11})^{-1} \cA^{12},  &&
\cG^{33} =   (\cB^{22})^{-1} ( \cA^{21}\cG^{13} + \cA^{22})  . 
\end{eqnarray*}
Note that $\cG^{11}$, $\cG^{12}$, $\cG^{31}$ and $\cG^{33}$ are first order 
(linear or affine  in $\zeta$), that  $\cG^{32}$ is second order (at most quadratic in $\zeta$)
and that $\cG^{13} $ is of order zero (independent of $\zeta$). 
We denote by $\cG^{ab}_{\mathrm{p}} $ their principal part (leading order part as polynomials).  
We note that 
\begin{equation}
\label{princcG}
\cG^{ab}_{\rmp} (z, \zeta) =  G^{ab}_\rmp (w(z), \zeta) 
 \quad  \mathrm{when} \  (a, b) \ne (3, 1) , 
\end{equation}
with 
\begin{equation*}
\begin{array}{lll}
G^{11}_{\mathrm{p}}(u, \zeta)  & = -  (A_d^{11}(u))^{-1}  \big(( \gamma+ i\tau) A_0^{11}(u) +  
\sum_{j=1}^{d-1} i \eta_j A_j^{11}(u) \big),  \\
G^{12}_{\mathrm{p}}(u, \zeta)  & = -  (A_d^{12}(u))^{-1}    
\sum_{j=1}^{d-1} i \eta_j A_j^{12}(u) \\
G^{13}_{\mathrm{p}}(u)  & = -  (A_d^{11}(u))^{-1}    A_d^{12}(u)
\\
G^{32}_{\mathrm{p}}(u, \zeta)  & =   (B^{22}(u))^{-1}  
\sum_{j, k =1}^{d-1}   \eta_j  \eta_k B_{j,k}^{22}(u)   \big), \\
G^{33}_{\mathrm{p}}(u, \zeta)  & =  -  (B^{22}(u))^{-1}  
\sum_{j=1}^{d-1} i \eta_j \big(B_{j, d}^{22}(u)  + B^{22}_{d,j} (u) \big) .  
 \end{array}
\end{equation*}
The principal term of $\cG^{3, 1}$ involves derivatives of the profile $w$. 
Denoting by $p = \lim_{z \to + \infty}  w(z) = w(\infty)$ the end state of 
the profile $w$, we note that the end state of $\cG^{31}_\rmp$ is 
$$
\cG^{31}_{\mathrm{p}}(\infty, \zeta)   =  (B^{22}(p))^{-1} 
\big(( \gamma+ i\tau) A_0^{21}(p) +  
\sum_{j=1}^{d-1} i \eta_j A_j^{21}(p)  + A^{21}_d(p) G^{11}_\rmp(p, \zeta) \big).   
$$
There are similar formulas using  the matrices $\overline A_j$
and $\overline B_{j, k}$ of \eqref{y2}.

The spectral analysis is easier 
when all the terms are reduced to first order. If 
$u = (u^1, u^2, u^3)$ is replaced by 
$\tilde u  =   \rmh _{|\zeta |} u  :=  (  u^1,   u^2,  | \zeta |^{-1}  u^3)$, $\cG$ is replaced by 
\begin{equation}
\label{dectG}
\tilde \cG =  \rmh_{| \zeta |} \cG  \rmh^{-1}_{| \zeta |} = \begin{pmatrix}
     \cG^{11} &   \cG^{12}  &  | \zeta |  \cG^{13}    \\  0 & 0&  | \zeta | \Id \\
       | \zeta |^{-1}  \cG^{31} &  | \zeta |^{-1}  \cG^{32}  & \cG^{33}  
\end{pmatrix} := \begin{pmatrix}
   \cG^{11}    &    \cP^{12}\\
   \cP^{21}   &   \cP^{22}
\end{pmatrix}
\end{equation}
with obvious definitions of $\cP^{a b}$. Note that $\tilde \cG$ is or order one, while $\cP^{21}$ is of order zero. Thus 
\begin{equation}
\label{dectG1}
\tilde \cG (z, \zeta)   = \tilde \cG_{\mathrm{p}} (z, \zeta)   + O(1) , \quad  
\tilde \cG_{\mathrm{p}}= \begin{pmatrix}
  \tilde \cG_{\mathrm{p}}^{11}    &  \cP^{12}_\rmp  \\    0 
      &    \cP^{22}_{\mathrm{p}}
\end{pmatrix}  = O(| \zeta |). 
\end{equation}
Moreover, since the coefficients in $\cG$ converge exponentially at infinity, the 
remainder in \eqref{dectG1} is uniform in $z \in \RR_+$ and $|\zeta | \ge 1$. 
Moreover, the principal part of $\tilde \cP^{22}$ is if the form
$ \tilde \cP^{22}_\rmp(z, \zeta) = P^{22}_\rmp(w(z), \zeta)$.

\begin{lem}
\label{lemspec2}  i)  For all $\zeta \in \overline \RR^{d+1}_+ $ 
with  $\gamma  > 0  $ and $ \eta \ne 0  $   and for all 
   and $z \ge 0$,     
  $ \tilde \cG_\rmp(z,  \zeta)$   has no eigenvalues on the imaginary axis; 
  moreover, the number of eigenvalues in $\{ \re \mu < 0 \}$ 
  is $N_b = N^1_+ + N'$.

ii) for all compact subset of $\cU^*$, there are $c > 0$ and $\delta > 0$ such that  for 
all $u$ in the given compact and all $\zeta \in \overline \RR^{d+1}_+ $ such that 
either  $    \gamma     \le \delta | \zeta |  $  
or  $ | \eta | \le \delta | \zeta |   $, the distance 
between the spectrum of $G^{11}_{\rmp}(u, \zeta) $ and the spectrum of $P^{22}_\rmp(u, \zeta)$
is larger than $c | \zeta |$.

\end{lem}

\begin{proof} The spectrum of $\tilde \cG_{\rmp}$ is the union of the spectra of 
$G^{11}_\rmp $ and $P^{22}_\rmp$. 
By homogeneity, it suffices to consider $\zeta \in \overline S^d_+ $. 

{\bf a) }   $G^{11}_\rmp$  is related to $L^{11}$  since 
$A^{11}_d ( i \xi + G^{11}_\rmp(u, \zeta))   = L^{11} (u, \gamma + i \tau, i \eta, i\xi)$. 
By Assumption (H3),  $L^{11}$ is hyperbolic in the time direction, hence  $G^{11}_\rmp$  has no eigenvalues on the imaginary axis  when  $\gamma > 0$; moreover, the boundary is noncharacteristic for $L^{11}$ by Assumption~\ref{assnoncar2}, implying that 
 the number of eigenvalues of $G^{11}_\rmp$ in $\{ \re \mu < 0 \}$ is equal 
 to the number of positive eigenvalues of $A_d^{11}$, that is 
  is $ N^1_+ $. 

Next, note that 
\begin{equation*}
P^{22}_\rmp  = \begin{pmatrix}
     0  &   | \zeta | \Id  \\
     | \zeta |^{-1} G^{32}_\rmp  &  G^{33}_\rmp 
\end{pmatrix}
\end{equation*}
Thus,  $i \xi $ is an eigenvalue of $P^{22}_\rmp$ if and only if 
$0$ is an eigenvalue of $B^{22} (\eta, \xi)$, which is impossible by (H2)  if  $\eta \ne 0$. 
Thus, the eigenvalues of $P^{22}_\rmp$ are not  purely imaginary when $\eta \ne 0$. 
   Moreover, the number of eigenvalues in $\{ \re \mu < 0 \}$ 
  is $    N'$ (see \cite{MZ1}). This finishes the proof of $i)$. 
  
  \medbreak

  {\bf b) } 
  If $\eta = 0$, $G^{32}_\rmp$  and $ G^{33}_\rmp $
  vanish, hence 
  the spectrum of $P^{22}_\rmp$ is $\{ 0 \}$. On the other hand $0$ is not an eigenvalue of    
  $G^{11}_\rmp = - (\gamma + i \tau) (A_d^{11})^{-1} A_0^{11}$  
     since $A_d^{11}$ and $  A_0^{11}$
  are invertible and $| \gamma + i \tau | = | \zeta | = 1 $.

If $\gamma = 0$ and $\eta \ne 0$, the eigenvalues of $P^{22}_\rmp$ are not  in $i \RR$.  
On the other hand, by Assumption (H6)  the eigenvalues of $G^{11}_\rmp$ are 
purely imaginary, thus  $P^{22}_\rmp$ and  $G^{11}_\rmp$ have no common eigenvalue. 
This finishes the proof of $ii)$. 
\end{proof}

The analysis in a purely  ``elliptic'' zone  
$\{   \gamma \ge \delta  | \zeta |\  \mathrm{and} \ 
| \eta | \ge  \delta | \zeta | \}$ with $\delta > 0$,  is easy, see below. The most 
difficult and important part 
is  to understand the    ``hyperbolic-parabolic'' decoupling in an arbitrarily small  cone 
 \begin{equation}
 \label{defCC}
 {\rC}_\delta   =  
 \{  0 \le  \gamma   \le     \delta  | \zeta | \}  \cup \{ | \eta | \le \delta  | \zeta | \} 
 \end{equation}   with  
 $\delta$  such that property $ii)$ of Lemma~\ref{lemspec2} holds for 
 $u$ in a simply connected neighborhood $\cU^*_0  $ of a compact set which contains the  
 curve $\{ w(z), \ z \in [0, + \infty [   \}$.  
 There,  the usual homogeneity and the parabolic homogeneity are in competition, 
 leading to different classes of symbols. 
 We use the following terminology: let $\zeta=(\tau,\gamma,\eta)$ and for a multi-index
$\alpha=(\alpha_\tau,\alpha_\eta, \alpha_\gamma)  \in \NN \times \NN^{d-1} \times \NN$, set
$$
|\alpha|=\alpha_\tau+|\alpha_\eta| \quad \mathrm{ and}
\quad \la \alpha \ra =2 (\alpha_\tau+\alpha_\gamma) +  |\alpha_\eta|. 
$$ 
Recall that the parabolic weight is $\Lambda = ( 1 + \tau^2 + \gamma^2 + | \eta |^4)^{\frac{1}{4}}$. 
\begin{defi}\label{w42a}

i )  $\Gamma^m(\Omega) $ denotes the space of \emph{homogeneous symbols} of
order $m$, that is of functions $h(z,  \zeta)\in
C^\infty(\overline \RR_+ \times \Omega)  $ such that there is $\theta > 0$ 
such that for all $\alpha \in \NN^{d+1}$ and all 
$k \in \NN$, there are constants $C_{\alpha, k}$  such that  for  $ |\zeta|\geq 1$ : 
\begin{align}
& |   \partial_{\zeta}^\alpha h|   \leq
C_{ \alpha, 0 }     |\zeta|^{m-|\alpha|},  & if \ k = 0,  &
\\
& | \partial_{z }^k \partial_{\zeta}^\alpha h|   \leq
C_{ \alpha, k }  e^{ - \theta  z}  |\zeta|^{m-|\alpha|},  & if \ k > 0, &
\end{align}

ii )   $\PG^m(\Omega) $ denotes 
the space of \emph{parabolic symbols} of order $m$, that is of functions $h(z,  \zeta)\in
C^\infty(\overline \RR_+ \times \Omega)  $ satisfying similar estimates with $| \zeta |^{m - | \alpha |}$
replaced by $\Lambda^{ m - \la \alpha \ra}$. 

 We use the same notation for spaces of homogeneous or parabolic
matrix symbols of any fixed dimension.
\end{defi}

 \begin{lem}
   For all $\hat \zeta \in   S^d \cap \rC_{\delta}$, there is a 
   a conical neighborhood $\Omega$ of $\hat \zeta$ and there are matrices 
  $ \cW^{12}_\rmp  \in \Gamma^0(\Omega) $ and $\cW^{21}_\rmp$,   
   homogenous of degree 0 in $\zeta$ for $u \in \cU^*_0$  such that 
 \begin{eqnarray}
\label{conjW2}
&&  \cW^{21}_\rmp  \cG^{11}_\rmp  -   \cP^{22}_{\rmp} \cW^{21}_\rmp =  | \zeta | \cP^{21}_\rmp. 
\\
\label{conjW2n}
 && \cG^{11}_\rmp  \cW^{12}_\rmp -   \cW^{12}_\rmp \cP^{22}_{\rmp} 
   =  -   \cP^{12}_\rmp. 
\end{eqnarray}
 \end{lem} 
 
 \begin{proof} By homogeneity, it is sufficient to construct $\cW^{21}_\rmp$ 
 for $|\zeta | = 1$. 
 By Lemma~\ref{lemspec2}, for $\zeta \in S^{d+1} \cap \rC_{\delta}$ and $u \in \cU^*_0$, the spectra
 of $G^{11}_\rmp(u, \zeta)$ and $P^{22}_\rmp(u, \zeta)$ do not intersect, so that the 
 linear system of equation 
 $$
X  G^{11}_\rmp(u, \zeta) -  P^{22}_\rmp (u, \zeta) X = Y 
 $$
 has a unique solution  $X  = \cX (u, \zeta) Y$. Therefore   
  $\cW^{21}_\rmp (z, \zeta) =  | \zeta |   \cX(w(z), \zeta) \cP^{21}_\rmp(z, \zeta)$ 
  satisfies \eqref{conjW2}  (Note that $\cP^{21}$ is of degree $0$). 
  
  The construction of $\cW^{12}_\rmp$ is similar, noticing that 
 $\cP^{12}_\rmp$ is of degree $1$. 
 \end{proof}

In the block structure of $\cG$, there holds
\begin{equation}
\label{conjW3}
\cW^{21}_\rmp = \begin{pmatrix}
     \cV^{21}_\rmp   \\   \cV^{31}_\rmp
\end{pmatrix}, \qquad 
\cW^{12}_\rmp = \begin{pmatrix}
     \cV^{12}_\rmp    &    \cV^{13}_\rmp
\end{pmatrix}
\end{equation}
and \eqref{conjW2} reads 
 \begin{eqnarray}
 \label{conjW2b}
\cV^{21}_\rmp \cG^{11}_\rmp   -    | \zeta |    \cV^{31}_\rmp    & = & 0 , 
\\
\label{conjW2c}
\cV^{31}_\rmp \cG^{11}_\rmp    
 - | \zeta |^{-1} \cG^{32}_\rmp  \cV^{21}_\rmp -      \cG^{33}_\rmp \cV^{31}_\rmp  & = &  \cG^{31}_\rmp. 
\end{eqnarray}
Similarly, 
 \begin{eqnarray}
 \label{conjW2nb}
 \cG^{11}_\rmp  \cV^{12}_\rmp   -    | \zeta |^{-1}   \cV^{13 }_\rmp \cG^{32}_\rmp & = &   -  \cG^{12}_\rmp  
\\
\label{conjW2nc}
 \cG^{11}_\rmp    \cV^{13}_\rmp
 - | \zeta |   \cV^{12}_\rmp -       \cV^{13}_\rmp  \cG^{33}_\rmp& = &  | \zeta | \cG^{13}_\rmp. 
\end{eqnarray}
 For further use, we make the following remark : by \eqref{princcG}, we see that 
  $\cG^{12}_\rmp$ and $\cG^{32}_\rmp$ vanish when $\eta = 0$. Therefore, 
  \eqref{conjW2nb} implies that $\cV^{12}$ also vanishes when $\eta = 0$ and hence
  \begin{equation}
  \label{vanV12}
  \cV^{12} (z, \zeta ) =  O ( | \eta |  / | \zeta |).  
  \end{equation}
  
\medbreak 
 With these notations, let 
 \begin{equation*}
\label{conjW5}
\cV_I (z, \zeta) =  
\begin{pmatrix}
  \Id      &   0  &  0    \\
  | \zeta |^{-1}    \cV^{21}_\rmp  &  \Id & 0 
     \\
     \cV^{31}_\rmp &0&\Id 
\end{pmatrix}, 
\quad 
\cV_{II} (z, \zeta) =  
\begin{pmatrix}
  \Id      &       \cV^{12}_\rmp   &   | \zeta |^{-1}    \cV^{13}_\rmp    \\
  0  &  \Id & 0 
     \\ 
    0  &0&\Id 
\end{pmatrix}
\end{equation*}
 and $   \cV = \cV_I \cV_{II}$.  
  Using the conjugation $u =\cV \hat u $, 
  $f =  \cV \hat f$, for $\zeta$  in the  in the cone $\rC_\delta$, the equation \eqref{linq4nn} is transformed 
  to 
  \begin{equation}
  \label{linq8}
  \D_z \hat u   =   \hat \cG  \hat u   + \hat f,  
  \quad \hat \Gamma \hat u(0) = g
  \end{equation}
  with $ \hat \cG =  \cV^{-1}  \cG \cV  -   \cV^{-1} \D_z \cV$ 
  and $ \hat \Gamma (\zeta) = \Gamma (\zeta) \cV (0, \zeta)$. 

\begin{lem}
\label{lem76mars}
The entries of $\hat \cG$ satisfy: 
\begin{eqnarray*}
&& \hat \cG^{11}    -  
 \big(  \cG^{11}    +   | \zeta |^{-1} \cG^{12}   \cV^{21}_\rmp   + \cG^{13}   \cV^{31}_\rmp \big)
 \in \Gamma^{-1}  , 
\\ 
& &  \hat \cG^{12}    \in \Gamma^0, \qquad   \hat \cG^{13}  \in \Gamma^{-1},
\qquad
 \hat \cG^{21} \in \Gamma^{-1}, \qquad  \hat \cG^{31} \in \Gamma^{0},
\\
&& 
 \hat \cG^{22} \in \Gamma^0, \qquad   \hat \cG^{23}  - \Id \in \Gamma^{-1},
 \\
 && \hat \cG^{32}  - ( \cG^{32} - V^{31} \cG^{12}) \in \Gamma^0 , 
\qquad   \hat \cG^{33}  -    \cG^{33}  \in \Gamma^0  , 
 \end{eqnarray*}
\end{lem}

\begin{proof}
We first compute the entries of $\cG_I =  \cV_I^{-1}  \cG \cV_I$. 
Direct computations show that 
$$
\begin{aligned}
 &   \cG_I^{11}    = 
 \cG^{11}   +  | \zeta |^{-1} \cG^{12}   \cV^{21}_\rmp   + \cG^{13}   \cV^{31}_\rmp  , 
\qquad 
  \cG_I^{12}  = \cG^{12}, \qquad    \cG^{13}_I  = \cG^{13}
 \\
 &   \cG_I^{32}  =   \cG^{32} - V^{31} \cG^{12} , 
\qquad  \cG_I^{33}  =   \cG^{33} - V^{31} \cG^{13} .  
\end{aligned}
 $$ 
Moreover, 
$$
 \cG_I^{21} =  -  | \zeta |^{-1} \cV^{21}_\rmp \cG^{11}  +  \cV^{31} 
 - | \zeta |^{-1} \cV^{21} \big( | \zeta |^{-1} \cG^{12} \cV^{21}_\rmp + 
 \cG^{13} \cV^{31}_\rmp \big)  . 
$$ 
The first two terms are of  degree zero, and by \eqref{conjW2b}, the sum of their 
principal terms vanishes; the third term is of degree $-1$  
thus $   \cG^{21}_I \in \Gamma^{-1}$. 
Similarly,  $ \cG^{31}_I$ is of degree 1 and its principal part vanishes 
by \eqref{conjW2c}. Thus, 
$$
 \cG^{21}_I  \in \Gamma^{-1}, \qquad    \cG^{31}_I \in \Gamma^{0}. 
 $$
Next   
$$
  \cG_I^{22} = - | \zeta |^{-1} \cV^{21}_\rmp \cG^{12} \in \Gamma^{0}, 
 \qquad 
  \cG_I^{22}  - \Id  = -  |\zeta |^{-1} \cV^{21} \cG^{13}  \in \Gamma^{-1}.
 $$
 
 The computations for $\cG_{II}  = \cV_{II}^{-1} \cG_I \cV_{II}$ are quite similar. 
 This new conjugation annihilates the principal parts of 
$ \cG_I^{12}$ and $\cG^{13}_I$ and contributes to remainder terms in the 
other entries. 

Finally, direct computations show that $ \cV^{-1} \D_z \cV $ only contributes to remainder. 
\end{proof}

 \medbreak 
    The main idea is to consider \eqref{linq8}  as a perturbation of the decoupled system
   \begin{eqnarray}
  \label{linq9}
 &&  \D_z \hat u^1   =     \hat \cG^{11}     \hat u^1  + \hat f_1, 
  \\ 
   \label{linq10}
&&  \D_z \begin{pmatrix}\hat u^2 \\ \hat u^3 \end{pmatrix}    =   
  \begin{pmatrix}0 & \Id  \\ \cG^{32} & \cG^{33}  \end{pmatrix} 
      \begin{pmatrix}\hat u^2 \\ \hat u^3 \end{pmatrix} 
    + \begin{pmatrix}\hat f^2 \\ \hat f^3 \end{pmatrix} . 
  \end{eqnarray}
  Introduce then 
  \begin{equation}
  \label{dec7nn}
  \cG'  = \hat \cG  -  \begin{pmatrix}\hat \cG^{11} & 0 & 0 
  \\ 0 & 0 & \Id  \\ 0 & \cG^{32} & \cG^{33}  \end{pmatrix} 
  \end{equation}
  The next lemma how the estimates are transported by the change of variables
  $u = \cV \hat u$. We use the notations \eqref{nota74} for the scaled norms.

  \begin{lem}
  \label{lemnorms}
  There are constant $C$ and $\rho_1$ such that for all $\zeta$ in the cone $\rC_\delta$ with 
  $| \zeta | \ge \rho_1$, there holds 
  \begin{equation}
\label{inva75nn}
\begin{aligned}
 & \| \cV^{-1}  \hat u\|_{sc}  \le C  \| \hat u \|_{sc}, \quad \quad \ \ 
   \|  \cV  f \|'_{sc}  \le C \| f \|'_{sc} ,
   \\
  &  | \cV^{-1} \hat u(0) |_{sc} \le C | \hat u(0) |_{sc},     \quad 
  |  \cV  u(0) |_{sc} \le C |   u(0) |_{sc},   
  \end{aligned}
\end{equation}
and 
\begin{equation}
\label{inva76nn}
| \hat \Gamma(\zeta)  \hat u (0) |_{sc}  \le   C | \hat u(0) |_{sc}. 
\end{equation}
Moreover, 
\begin{equation}
\label{invafnn}
\|  \cG' \hat u  \|_{sc}  \le   C \Lambda^{-1}   \|  \hat u \|_{sc}. 
\end{equation}
  
  \end{lem}

  \begin{proof}   Direct computations, using \eqref{vanV12},  show that  $ u = \cV \hat u$ satisfies 
  \begin{eqnarray*}
  && u^1 = O( 1) \hat u^1  + O(| \eta| \,   | \zeta |^{-1} ) \hat u^2  + 
  O(    | \zeta |^{-1} ) \hat u^3, 
  \\
  && u^2 = O(  | \zeta |^{-1} ) \hat u^1  + O(1 ) \hat u^2  + 
  O(    | \zeta |^{-1} ) \hat u^3, 
  \\
   && u^3  = O( 1) \hat u^1  + O(1 ) \hat u^2  + 
  O(  1 ) \hat u^3. 
  \end{eqnarray*}
  This implies the first estimate in \eqref{inva75nn}, using 
  the inequalities 
  $$
  (1+ \gamma) | \eta | / | \zeta |  \ls  \Lambda  , \quad  (1+ \gamma) /  | \zeta | \ls 1, \quad 
  \Lambda / | \zeta | \ls 1. 
  $$
  The proof of the other estimates of \eqref{inva75nn} is similar, using in particular for the 
  traces the inequality $  (1+ \gamma)^\mez | \eta | / | \zeta |  \ls  \Lambda^\mez$.

  The inequality \eqref{inva76nn} follows from the second line of \eqref{inva75nn} 
  and the estimate $| \Gamma u(0) |_{sc} \le | u(0) |_{sc}$ which is a direct 
  consequence of the form \eqref{linviscbcd} of the boundary conditions.  
  
  Finally, Lemma~\ref{lem76mars} implies that $ \hat f =  \cG' \hat u $ satisfies
   \begin{eqnarray*}
  && \hat f^1 =   O(1 ) \hat u^2  + 
  O(    | \zeta |^{-1} ) \hat u^3, 
  \\
  && \hat f^2 = O(  | \zeta |^{-1} ) \hat u^1  + O(1 ) \hat u^2  + 
  O(    | \zeta |^{-1} ) \hat u^3, 
  \\
   && f^3  = O( 1) \hat u^1  + O(1 ) \hat u^2  + 
  O(  1 ) \hat u^3. 
  \end{eqnarray*}
and \eqref{invafnn} follows. 
  \end{proof}

The parabolic bloc \eqref{linq10} is studied in \cite{MZ1}. We now focus on 
the hyperbolic block \eqref{linq9}, recalling and extending the  analysis of 
\cite{GMWZ4}. 

 
 \subsection{Analysis of the hyperbolic block. }

\subsubsection{The genuine coupling condition}

For $u \in \cU^*$, denote by $\lambda_j(u, \xi)$ the distinct eigenvalues 
of $ \overline A^{11} (u, \xi)$, which are real and have constant multiplicity $\nu_j$
by Assumption (H5). 
 Assumption (H6) implies the following: 
 
 \begin{lem}
 \label{lemn76}  For all $u\in \cU^*$,  all $\xi \in \RR^d$ and all $j$, there holds 
 $\D_{\xi_d} \lambda_j (u , \xi ) \ne 0$, and all these derivatives have the same sign.
 \end{lem}

\begin{proof}
If  $\D_{\xi_d} \lambda_j (u , \underline \eta, \underline \xi_d  ) = 0$, then the equation 
$\tau + \lambda(\underline \eta, \xi_d ) = 0$ would have complex roots 
in $\xi_d$ for some $\tau $ close to $\underline \tau = - \lambda_j (u , \underline \eta, \underline \xi_d  )$ (recall that $\lambda_j$ is  real analytic).  Thus hyperbolicity in the normal direction prevents glancing. Moroever, by continuity the sign of 
$\D_{\xi_d} \lambda_j (u , \eta, \xi_d )$ is  constant for all $\xi_d \in \RR$ when 
$\eta \ne 0$. Thus the functions $\xi_d \mapsto \lambda_j (u,  \eta, \xi_d)$ are monotone
and tend to infinity as $\xi_d$ tends to $\pm \infty$. Since 
$\lambda_j \ne \lambda_k$ when $j\ne k$, they must  be all increasing or all decreasing.
This remains true for $\eta = 0$ by continuity.  
\end{proof}

According to the terminology of Section 4, we will say that the 
hyperbolic block $L^{11}$ is \emph{incoming}  [resp. \emph{outgoing}] 
when the derivatives   $ \D_{\xi_d} \lambda_j (u , \xi ) $  
are positive [resp. negative]. 

\begin{cor}
\label{corn77} 
i)  The matrix $G_\rmp^{11} (u, \zeta) $  has  no purely imaginary eigenvalues when 
$\gamma > 0$.  They are all lying in $\{ \re \mu > 0 \}$  if the $11$-block is outgoing 
 and in $\{ \re \mu < 0 \}$ if it is incoming. 

ii) Near points $\underline \zeta$ with  $ \underline \gamma =0$, 
 $G^{11}_\rmp (u, \zeta) $ has semi-simple eigenvalues $\mu_j (u, \zeta)$ of constant 
 multiplicity $\nu_j$, which are purely imaginary when $\gamma = 0$. Moreover, 
 $\D_ \gamma \re \mu_j  > 0$  when the $11$-block is outgoing and 
 $\D_ \gamma \re \mu_j  <  0$  when the $11$-block is incoming.  
\end{cor}
\begin{proof}
Note that  $ \mu $ is an eigenvalue of $G_\rmp^{11} (u, \zeta) $ if and only if $ -  \tau + i \gamma$ is an eigenvalue of $\overline A^{11}(u, \eta, \xi) $ 
with $\xi  = - i \mu $. 

Consider the equations in $\xi_d$ : 
$\tau + \lambda_j (u, \eta, \xi_d) = 0$. Since $\lambda_j$ is strictly monotone and tends 
to infinity at both infinity, it always have a unique solution, $\psi_j (u, \eta , \tau)$
and $\D_\tau \psi_j$ has the same sign as $ - \D_{\xi_d} \lambda_j$. 
This solution  extends  analytically for $\im \tau $ small. 
This yields  distinct eigenvalues $\mu_j (u, \zeta) = i \psi_j (u, \eta, \tau - i \gamma)$ of $G^{11}_\rmp$  for $\zeta$ close to the real domain. 
In particular   $  \D_ \gamma \mu_j  = \D_\tau \psi_j$ and the eigenvalues all lie in 
 $\{ \re \mu > 0 \}$  if the $11$-block is outgoing 
 and in $\{ \re \mu > 0 \}$ if it is incoming. 
 
 The kernel of $G^{11}_\rmp - \mu_j$ is the kernel of 
 $ \overline A^{11} - \lambda_j$, thus has dimension equal to the multiplicity of
$ \lambda_j$.  Since  these dimensions add up to $N^1$, this shows that 
$G^{11}_\rmp$ has only semi-simple eigenvalues of constant multiplicity, which all 
lie in a given half space when $\gamma > 0$.

Hyperbolicity of $L^{11}$ implies that   $G_\rmp^{11} (u, \zeta) $  has  no purely imaginary eigenvalues when 
$\gamma \ne  0$ and by continuity they all lie in the same half space.
\end{proof}

Next we need more information on the zero-th order correction of $\hat \cG^{11}$. 
From \eqref{coefflin} \eqref{Einfty} and \eqref{deccG} we see that
\begin{equation}
\label{}
\hat \cG^{11} (z, \zeta)  - (\cV^{-1} \D_z \cV)^{11}  \big)   = G^{11}_\rmp (w(z), \zeta)   
+  \cE (z, \zeta),  
\end{equation}
where  $\cE \in \Gamma^0$. Denote its principal part by $\cE_\rmp$. 
Its  limit at $z = \infty $ is 
\begin{equation}
E_\rmp (p, \zeta) =    | \zeta |^{-1} G^{12}_\rmp (p,  \zeta)  V^{21}_\rmp (p  , \zeta)    
 G^{13}_\rmp (p ,  \zeta) V^{31}_\rmp (p , \zeta) 
\end{equation}
where $p =\lim_{z \to +\infty}  w(z)$ and $V_\rmp^{21}(p, \zeta)$, $V_\rmp^{31}(p, \zeta)$ 
denote the end points of $\cV_\rmp^{21}$ and $\cV_\rmp^{31}$, 
that is the solutions of the intertwining relations \eqref{conjW2b} \eqref{conjW2c}
with matrices $\cG^{ab}_\rmp$ replaced bay their endpoint values 
$G^{ab}_\rmp(p, \zeta)$. 
The next result is crucial and follows from the genuine coupling condition 
 (H4). 

\begin{prop}
\label{propn711}
Fix   $\underline \zeta$   with $|\underline \zeta | = 1$ and $\underline \gamma = 0$. 
For $\zeta$ in a neighborhood of $\underline \zeta$, consider a basis  where  
 $G^{11}(u, \zeta)$ has the block diagonal form 
 $\diag (\mu_j  \Id_{\nu_j})$. Denote by $E_{j, k} (u, \zeta)$ the corresponding blocks 
 of $E$ is this basis. 
 Then, for $u \in \cU$ the eingenvalues of the diagonal blocks    $\re E_{j,j}$  
 have a positive [resp. negative]  real part if the 
 $11$-block is outgoing  [resp. incoming]. 
\end{prop}

\begin{proof}
It is sufficient to prove the positivity  at 
$\underline \zeta$. 
Suppose that $\gamma = 0$, denote by $\vp_{j, p}  $ with $p \in \{ 1, \ldots, \nu_j\}$ a basis of eigenvectors 
of $G^{11}(u,   \zeta)$. 
Fix  $j$ and set  $   \xi_d  = - i \mu_j (u,   \zeta) \in \RR$, 
$\xi = (\eta, \xi_d)$. Then the $\vp_{j, p}$ are right eigenvectors of 
$\overline A^{11} (u, \xi)$ associated to the eigenvalue $- \tau = \lambda_j(u, \xi)$. 

Consider left eigenvectors $\ell_{j, p} $ of  $\overline A^{11} (u, \xi)$, dual to 
the $\vp_{j, p}$. Then, 
the left eigenvectors of $G^{11}_\rmp(u, \zeta)$ associated 
to $\mu_j$ are 
$\frac{1}{\beta_j} \ell_j \overline A^{11}_d $ with 
$\beta_j = \D_{\xi_d} \lambda_j(u, \eta, \xi$, wee Lemma~\ref{lem32}. The entries of the block $E_{j, j}$   are 
\begin{equation}
\label{coefEjj}
\frac{1}{\beta_j}  \ell_{j, p}  \overline A_d^{11} E_\rmp  (u, \zeta)  \vp_{j, p'}. 
\end{equation}
Computing the eigenvalues of order $\eps$ of 
$ \overline B (u, \xi) +  i  \eps \overline A(u, \xi)$, leads to consider the matrix  
\begin{equation}
  i  \overline A^{11} +   \eps \overline A^{12} (\overline B^{22})^{-1}
\overline A^{21} .
\end{equation}
The genuine coupling condition (H4) implies that for $u \in \cU$, its spectrum lies in 
$\re \mu > c \eps$ for $\eps$ small, and this implies that  the matrix
$F_{j, j}$ with entries 
\begin{equation}
\label{coeffFjj}
\ell_{j, p} \overline A^{12} (\overline B^{22})^{-1}
\overline A^{21}   \vp_{j, p'}
\end{equation}
has its eigenvalues in the right half plane $\{ \re \mu > 0 \}$. 

Because $G^{11}_\rmp \vp_{j, p'} = i \xi_d \vp_{j, p'}$, the relation 
\eqref{conjW2b} implies  
$$
V^{31}_\rmp  \vp_{j, p'} =  | \zeta |^{-1}  V_\rmp ^{21} G^{11}_\rmp
\vp_{j, p'} =  i \xi_d | \zeta |^{-1}  V_\rmp ^{21} 
\vp_{j, p'} ,  
$$
and, using the expressions of the matrices $G^{a, b}$ yields
$$
(| \zeta |^{-1} G^{12}_\rmp V^{21}_\rmp + G^{13}_ \rmp V^{31}_\rmp )\vp_{j, p'} 
= -  i | \zeta |^{-1} (\overline A_d^{11})^{-1} \overline A^{12}(\eta, \xi)  V^{21}_\rmp \vp_{j, p'} 
$$
and  
$$
( | \zeta |^{-1} G^{32}_\rmp V^{21}_\rmp + G^{33}_\rmp V^{31}_\rmp - 
V^{31}_{\rmp}  G^{11}_\rmp ) \vp_{j, p'} =  | \zeta |^{-1} 
(\overline B^{22}_{dd}) ^{-1}  \overline B_{22} (\eta,\xi)  V^{21}_\rmp \vp_{j, p'} 
$$
By \eqref{conjW2c} this is equal to 
$$
- G^{31}_\rmp \vp_{j, p'} = -  i  (\overline B^{22}_{dd}) ^{-1}   
\overline A^{21} (\eta, \xi) \vp_{j, p'}. 
$$
Thus 
$$
 | \zeta |^{-1}   V^{21}_\rmp \vp_{j, p'}  =  - i 
 \big( \overline B_{22} (\eta,\xi) \big)^{-1}   
\overline A^{21} (\eta, \xi) \vp_{j, p'}. 
$$
and 
$$
E_\rmp \vp_{j, p'}  =  -  (\overline A_d^{11})^{-1}  \overline A^{12}(\eta, \xi) 
 \big( \overline B_{22} (\eta,\xi) \big)^{-1}   
\overline A^{21} (\eta, \xi) \vp_{j, p'}. 
$$
Multiplying on the left by $\ell_j \overline A^{11}_d$, this shows that the 
coefficients in \eqref{coefEjj} and \eqref{coeffFjj} only differ by the factor 
$- 1/ \beta_j$, and the proposition follows. 
\end{proof}


\subsubsection{Estimates}

We are now in position to prove maximal estimates for the solutions of the equation 
\eqref{linq9}.

\begin{prop} 
\label{propn79}
 There are constants $C$ and $\rho_1 \ge 1$ such that 
for all $\zeta $ in the cone $\rC_\delta$ with $| \zeta | \ge \rho_1$ and all 
$\hat u^1$ and $\hat f^1$ in $L^2(\RR_+)$ satisfying  $\eqref{linq9}$, 
there holds
\begin{equation}
\label{estHF11}
\begin{aligned}
(1+ \gamma) \| \hat u^1 \|_{L^2} +   (1+&  \gamma)^\mez |  \hat u^{1 +}   (0) | 
 \\
 &  \le   C \Big(  \| \hat f^1 \|_{L^2}  +  (1+ \gamma)^\mez |  \hat u^{1 -}   (0) |  \Big)
\end{aligned}
\end{equation}
where $\hat u^{1 +} = \hat u^1$ and  $\hat u^{1 -} = 0 $ if the $11$-block is outgoing 
and  $\hat u^{1 +} = 0$ and  $\hat u^{1 -} = \hat u^1 $ if it is incoming. 
\end{prop}

\begin{proof}
{\bf a) } 
Fix $\underline \zeta \in \overline S^{d+1}_+$.  We prove the estimate for 
$\zeta$ in a conical neighborhood of 
$\underline \zeta$. Suppose first that $\underline \gamma = 0$ (the most difficult case). 
By Corollary~\ref{corn77} there is a matrix $\cV^{11} (z, \zeta)$ homogeneous 
of degree $0$ such that $(\cV^{11})^{-1} \cG^{11}_\rmp \cV^{11} 
= \diag ( \mu_j (w(z), \zeta) \Id_{\nu_j} $.  Setting  $\hat u^1 = \cV^{11} u^1$ 
transforms the equation to 
\begin{equation}
\D_z  u^1  =     (\diag ( \mu_j (w(z), \zeta) \Id_{\nu_j}  + \widetilde  \cE ) u^1 + f^1 
\end{equation}
with $\widetilde \cE = \cE - ( \cV^{11})^{-1} \D_z \cV^{11} \in \Gamma^0$, whose principal 
part $\widetilde \cE_\rmp$ has the same end point $E_\rmp (p, \zeta)$ as $\cE_\rmp$. 

As usual, since the $\mu_j$ are pairwise distinct, there is a new change $u^1 = (\Id + \cV_{-1}) \tilde u^1$ with 
$\cV^{11}_{-1} \in \Gamma^{-1}$,  such that the resulting system has the same form 
with the additional property that the zero-th order part is also block diagonal, so that 
$\widetilde \cE_\rmp  = \diag (\cE_{j,j} ) $  and the end points of the 
blocks $\cE_{j,j}$ are $E_{j, j} $ introduced in Proposition~\ref{propn711}. 

The term $(\widetilde \cE - \widetilde \cE_\rmp) u$ is 
$O(| \zeta |^{-1} | u| )$, is incorporated to $f^1$ and finally absorbed from the right 
to the left of the inequality by choosing $| \zeta |$ large enough. 
This reduces the proof to the case where the equation reads
\begin{equation}
\label{linq11} 
\D_z \hat u^1 = \mu_j(w(z), \zeta)  \hat u^1  + E_{j, j}  (\zeta) \hat u^1 
+ F_{j, j} (z, \zeta) \hat u^1  + \hat f^1
\end{equation}
with $ |  F_{j,j} |  \le   C_0 e^{ - \theta z } $.

Consider the outgoing case. Then, Corollary~\ref{corn77} implies that  there is a constant 
$c > 0$ such that 
 $  \re \mu_j (u, \zeta)  \ge c \gamma $.  Moreover, Proposition~\ref{propn711} implies that 
 the eigenvalues of $E_{j, j} $ have a positive real part. Thus, there is a positive
 definite (constant) matrix $S (\zeta)\ ge \Id $ such that $\re S E_{j, j} $  is definite positive, 
 say $\re S E_{j, j} \ge \Id  $. 
 Introduce   $ a = C_0 | S |  \int_0^z e^{ - \theta s } ds$ such that 
$\D_z a \ge  | S F_{j, j}| $ and $ a$ is  bounded in $L^\infty$ uniformly with respect to 
$\zeta$. 
 Therefore, multiplying the equation by  
 $e ^{ 2 a(z) } S $  and taking the $L^2$ scalar product with $\hat u^1$  implies that 
 $$
( 1 + c \gamma)  \| e^{ a} \hat u^1 \|_{L^2}^2  + | \hat u^1(0) |^2 \le  C 
 \| e^{ a} \hat u^1 \|_{L^2}    \| e^{ a} \hat f^1 \|_{L^2}
 $$ 
 which implies \eqref{estHF11}. The proof in the incoming case is similar. 
 
 \medbreak
 {\bf b) } Suppose next that $\underline \gamma = 0$. Consider again the outgoing case. 
 Then, the eigenvalues of $G^{11}_\rmp$ satisfy 
 $\re \mu \ge c | \zeta |$ in a conical neighborhood of $\underline \zeta$. 
 This is the classical ``elliptic'' case. There is a symmetric definite positive 
 matrix $S(u, \zeta)\in \Gamma^0$ such that $\re  S G^{11} \ge c | \zeta | \Id$ and 
 usual integrations by parts imply that
  $$
 c | \zeta |   \|   \hat u^1 \|^2_{L^2}  + | \hat u^1(0) |^2 \le 
 C   \|   \hat u^1 \|_{L^2}    \|   \hat f^1 \|_{L^2}  +   C_1  \|   \hat u^1 \|^2_{L^2}  
 $$ 
 where $C_1$ involve estimates of the zero-th order terms, which include 
 $\D_z S(w(z) , \zeta)  $.  This term   is eliminated choosing $| \zeta | $ large enough. 
 The proof in the incoming case is similar.
\end{proof}

\begin{rem}
\label{remH5}
\textup{The proof above contains two ingredients.  First,  the $11$-block is totally incoming 
or totally outgoing, in analogy with the terminology of Section $4$. Thus the decoupling
incoming/outgoing is trivial. More generally, this could be replaced by a decoupling condition in the spirit of Section 4. For instance, for shocks, such a decoupling is immediate in \cite{GMWZ4} corresponding to equations on each side of the front. 
Next, we construct symmetrizers for the  incoming and outgoing components.  
There we use the genuine coupling condition.  If the eigenvalues are  not of constant 
multiplicity one can introduce adapted bases or use symmetry also in the spirit 
of Section $4$. }
\end{rem}

\subsubsection{About Assumption (H6)}

We show on an example that hyperbolicity in the normal direction is crucial in 
the proof of estimates of the form \eqref{estHF11}.  Suppose that the $L^{11}$- block 
reads 
\begin{equation}
\label{cexple731}
\left\{ \begin{array}{ll}
\D_t u - \D_y u  + \D_x v    , \\ 
\D_t v +  \D_y v  + \D_x u     .   
\end{array}\right. 
\end{equation}  
Then,    on the Fourier side, the $11$ equation will be of the form 
\begin{equation}
\label{cexple732}
\left\{ \begin{array}{ll}
( i (\tau - \eta)    + \gamma ) u  + \D_z v  +  a (z)  u  = f , \\ 
( i (\tau + \eta)    + \gamma ) v  + \D_x u   +   a (z)  v  =   g ,  
\end{array}\right. 
\end{equation}
and the only information we have from the genuine coupling condition is that 
$a $ is positive at $z = + \infty$. Suppose that 
$a (z_0 ) < 0$ for some $z_0 > 0$. 
Then glancing waves for \eqref{cexple731} will propagate parallel to the boundary and 
thus may remain in a region where $a$ is negative and thus may never be damped. 
This is illustrated by choosing 
$ \tau = \eta $, large,  $\gamma = - a(z_0) > $ and 
$$
u_\tau(z) =   \chi ( \tau^{\frac{1}{3}} (z - z_0)) , \quad 
v_\tau (z) =  \frac{- \D_z u_\tau  }{  2 i  \tau     + \gamma + a } 
$$
with $\chi \in C^\infty_0(\RR)$. 
Then \eqref{cexple732} is satisfied with 
$f =  (a(z) - a(z_0) ) u_\tau +\D_z v _\tau$ and $g=0$.
Moreover,  
$ \| f \|_{L^2} = O (\tau^{ - \frac{1}{3}})  \| u \|_{L^2}$  and $u(0) = v(0) =0$, 
showing that no estimate of the form \eqref{estHF11} can be valid.


\subsection{Proof of Theorem~\ref{tm7n1} }

\subsubsection{In the cone $\rC_\delta$} 

We consider now the equation \eqref{linq10} and briefly recall the results from
\cite{MZ1}. It is natural to rescale  the problem using the parabolic weights: 
with $v^2 =  \hat u^2$ and $v^3 = \Lambda^{-1}  \hat u^3$ 
and $g^2 =  \hat f^2$ and $g^3 =\Lambda^{-1}  \hat f^3$the system reads
\begin{equation}
\label{linq22}
\D_z \begin{pmatrix} v^2 \\ v^3 \end{pmatrix} =  \cG_P \begin{pmatrix} v^2 \\ v^3 \end{pmatrix}  + \begin{pmatrix} g^2 \\ g^3 \end{pmatrix}
\end{equation}
with 
$$
\cG_P = 
\begin{pmatrix} 0 & \Lambda \Id  \\ \Lambda^{-1} \cG^{32}  & \cG^{31}  \end{pmatrix} 
\in \PG^1
$$
of quasi-homogenenous degree one and principal part 
$ 
G_P(w(z), \zeta)
$ 
with 
\begin{equation}
G_P(u, \zeta) = 
\begin{pmatrix} 0 & \Lambda \Id  \\ \Lambda^{-1}\big( 
 (i \tau + \gamma) ( \overline B^{22})^{-1}  +    G^{32}_{\rmp} (u, \eta) \big)  & 
G_{\rmp} ^{31} (u, \eta)  \end{pmatrix} 
\end{equation}

\begin{lem}[\cite{MZ1}] 
There is $c >0$ such that the spectrum of $G_P$ lies in $\{ | \re \mu | \ge c \Lambda \}$, 
with $N'$ eigenvalues, counted with their multiplicity, of positive real part.  
There is a smooth change of  variables 
$\cW \in \PG^0$ such that 
$$
\cW^{-1} \cG_P \cW  = 
\begin{pmatrix} \cP_+   & 0 \\0  &\cP_-  \end{pmatrix} 
$$
with $\cP_\pm \in \PG^1$ having their eigenvalues satisfying 
$\pm \re \mu \ge c \Lambda$. 
\end{lem} 

Introduce 
$$ 
\begin{pmatrix} v^+    \\ v^-   \end{pmatrix}  = \cW^{-1} \begin{pmatrix} v^2    \\ v^3    \end{pmatrix}  
$$

\begin{cor}[ \cite{MZ1}] 
There are $C$ and $\rho_1$ such that for all 
$\zeta \in \rC_{\delta}$ with $| \zeta | \ge \rho_1$, there holds 
\begin{eqnarray*}
&& \Lambda  \| v^+ \|_{L^2} + \Lambda^\mez | v^+(0) | 
\le C \| (\D_z - \cP^+) v^+ \|_{L^2} , 
\\
&& \Lambda  \| v^- \|_{L^2}  \le C \| (\D_z - \cP^-) v^- \|_{L^2} + 
C \Lambda^\mez | v^-(0) | . 
\end{eqnarray*}
\end{cor}

Scaling back, introduce
\begin{equation}
 \begin{pmatrix} \hat u^{2, +}     \\ \hat u^{3, +}  \end{pmatrix} = 
 \begin{pmatrix} \Id    & 0 \\0  &\Lambda   \end{pmatrix} \cW  
 \begin{pmatrix} v^+    \\ 0   \end{pmatrix} , 
 \quad 
  \begin{pmatrix} \hat u^{2, -}     \\ \hat u^{3, -}  \end{pmatrix} = 
 \begin{pmatrix} \Id    & 0 \\0  &\Lambda   \end{pmatrix} \cW  
 \begin{pmatrix} 0    \\  v^-   \end{pmatrix} . 
\end{equation}
Because,  $\cW^{-1} \D_z \cW$ is uniformly bounded, the Corollary implies the
following estimate: 

\begin{prop}  
There are $C$ and $\rho_1$ such that for all 
$\zeta \in \rC_{\delta}$ with $| \zeta | \ge \rho_1$, there holds 
\begin{equation*}
\begin{aligned}
 \Lambda  \| u^{2,+} \|_{L^2} +&  \| u^{3,+} \|_{L^2} + 
\Lambda^\mez | u^{2, +} (0) |  + \Lambda^{-\mez } | u^{3, +} (0) |
\\
& \le C  \|  \hat f^2 \|_{L^2}  + C \Lambda^{-1}  \|  \hat f^3 \|_{L^2}
+ \|  \hat u^2 \|_{L^2}  + C \Lambda^{-1}  \|  \hat u^3 \|_{L^2}, 
\end{aligned} 
\end{equation*}
\begin{equation*}
\begin{aligned}
 \Lambda  \| u^{2,-} \|_{L^2} +&  \| u^{3,-} \|_{L^2}  \le  
 C \Lambda^\mez | u^{2, -} (0) |  + C  \Lambda^{-\mez } | u^{3, -} (0) |
\\
& +  C  \|  \hat f^2 \|_{L^2}  + C \Lambda^{-1}  \|  \hat f^3 \|_{L^2}
+ \|  \hat u^2 \|_{L^2}  + C \Lambda^{-1}  \|  \hat u^3 \|_{L^2} . 
\end{aligned} 
\end{equation*}
\end{prop}

Finally, with $\hat u^{1, \pm}$ as in Proposition~\ref{propn79},  introduce 
\begin{equation}
\label{splitmna}
\hat u^\pm = {}^t (\hat u^{1, \pm}, \hat u^{2, \pm}, \hat u^{2, \pm}) 
\end{equation}
Adding up the various estimates and using \eqref{invafnn}, 
one obtains the following estimates. 

 \begin{prop}  
 \label{prop716}
There are $C$ and $\rho_1$ such that for all 
$\zeta \in \rC_{\delta}$ with $| \zeta | \ge \rho_1$
and all $\hat u \in H^1(\overline \RR_+)$:
 \begin{eqnarray}
\label{estu+c}
&&\| \hat u^+ \|_{sc}   + | \hat u^+(0) |  \le C \| (\D_z - \cG) \hat u  \|_{sc} + 
\Lambda^{-1}  \| \hat u \|_{sc} , 
\\
&&\label{estu-c}
\| \hat u^- \|_{sc}     \le C \|  (\D_z - \cG) \hat u   \|_{sc} + 
\Lambda^{-1}   \| \hat u \|_{sc}    +  C  | \hat u^-(0) |.   
\end{eqnarray}
 \end{prop}

As indicated at the end of Section 7.1, these estimates imply the maximal estimates of Theorem~\ref{tm7n1} provided that the boundary conditions are  uniformly  spectral stable.


\subsubsection{Analysis in the central zone}
 We now consider the remaining cone where 
 \begin{equation}
 \label{internzone}
 \zeta \in \RR^{d+1}, \quad  \gamma \ge \delta | \zeta | \quad
 \mathrm{ and } \    | \eta | \ge \delta | \zeta |. 
 \end{equation} 
 We consider the rescaled $\tilde \cG$ matrix \eqref{dectG1}, for the rescaled unknows
 $\tilde u =  \rmh _{|\zeta |} u  :=  (  u^1,   u^2,  | \zeta |^{-1}  u^3)$, 
 $\tilde f =  \rmh _{|\zeta |} f  :=  (  f^1,   f^2,  | \zeta |^{-1}  f^3)$. 
 We note that in the region under consideration we now have 
 $(1+ \gamma) \approx \Lambda \approx |\zeta |$, so that the 
 rescaled norms \eqref{nota74} are equivalent to 
 \begin{equation}
 \label{resc761} 
 \begin{aligned}
 & \| u \|_{sc} \approx  | \zeta | \| \tilde u \|_{L^2}, 
 \\
 & | u(0) |_{sc} \approx | \zeta |^\mez | \tilde u(0) | , 
 \\
  & \| f \|'_{sc} \approx   \| \tilde f \|_{L^2}. 
 \end{aligned}
 \end{equation}

By Lemma~\ref{lemspec2},  there is a smooth matrix $\cV \in \Gamma^0  $ such that 
\begin{equation*}
\cV^{-1} (z \zeta)  \cG_\rmp (z, \zeta) \cV  (z, \zeta)   =  \begin{pmatrix}
     \cG^{+}_\rmp &  0  \\
     0 &    \cG^{-}_\rmp
\end{pmatrix} := \cG^{\diag}_{\rmp}
\end{equation*}
where the spectrum of $\cG^\pm_{\rmp}\in \Gamma^1 $ is contained in 
$\{ \pm \re \mu \ge c | \zeta | \}$. 
 We use the notations
 \begin{equation}
\hat u :=  \cV \tilde u  = \begin{pmatrix}  \hat  u^+ \\ \hat  u^- \end{pmatrix}. 
 \end{equation}
 $\tilde u^+ $ has dimension $N+ N' - N_b$ and $u^-$ has dimension $N_b$. 
 The equation for $\hat u$ reads
 \begin{equation}
 \label{eq771}
 \D_z  \hat u = \hat \cG \hat u + \hat f ,  
 \end{equation}
 with  $\hat \cG = \cG^\diag + O(1)$. The ellipticity of $\cG^\diag$ immediately 
 implies the following estimates.

 \begin{prop}
 There are constants $C$ and $\rho_1$ such that for all 
 $\zeta$ satisfying $\eqref{internzone}$ and $| \zeta | \ge \rho_1$
 and all $\tilde u \in H^1(  \RR_+)$ satisfying $\eqref{eq771}$, there holds
 \begin{eqnarray}
\label{199a}
&& | \zeta | \, \| \hat u^+ \|_{L^2}  +  | \zeta |^\mez | u^+ (0) |   \le   C  \| \hat f \|_{L^2}
+ C \| \hat u \|_{L^2}, 
 \\
\label{199b}
&& | \zeta | \, \| \hat u^-  \|_{L^2}          \le   C  \| \hat f \|_{L^2}
+ C \| \hat u \|_{L^2}  + C 
  | \zeta |^\mez | \hat u^- (0) |^2   . 
\end{eqnarray}

 \end{prop} 

 Thanks to \eqref{resc761}, this is the exact analogue of Proposition~\ref{prop716}
 and these estimates imply the maximal estimates of Theorem~\ref{tm7n1} provided that the boundary conditions are  uniformly  spectral stable, as explained in Section~7.1.



     %
     %

\section{Application to magnetohydrodynamics}
\label{SecMHD}

We now apply our results to the equations of isentropic
magnetohydrodynamics (MHD), for which the inviscid case was
treated in \cite{MZ2}.  The full (nonisentropic) inviscid equations have
been treated in \cite{Kw}, and have essentially the
same symbolic structure as the isentropic inviscid equations. 

\subsection{The equations} 

The equations of isentropic magnetohydrodynamics (MHD) appear in 
basic form as
\begin{equation}
\label{mhdeq}
\left\{ \begin{aligned}
 & \D_t \rho +  \div (\rho u) = 0
 \\
 &\D_t(\rho  u) + \div(\rho u^tu)+ \na p + H \times \curl H =  \eps \nu \Delta u 
 \\
 &  
 \D_t H + \curl (H \times u) = \eps \mu \Delta H 
 \end{aligned}\right.
\end{equation}
\begin{equation}\label{divfree}
\div H=0,
\end{equation}
where $\rho\in \RR$ represents density, $u\in \RR^3$ fluid velocity,
$p=p(\rho)\in \RR$ pressure, and $H\in \RR^3$ magnetic field.
When  $H\equiv 0$, \eqref{mhdeq} reduces to the equations of isentropic
fluid dynamics.
We assume that $\nu $ and $\mu $ are positive. 

Equations \eqref{mhdeq} may be put in conservative form using identity
\begin{equation}
H\times \curl H= (1/2)\div (|H|^2I-2H^tH)^\trans + H\div H
\end{equation}
together with constraint \eqref{divfree} to express the second equation
as
\begin{equation}\label{cons2}
 \D_t(\rho  u) + \div(\rho u^tu)+ \na p + (1/2)\div (|H|^2I-2H^tH)^\trans= \eps  \nu \Delta u.
\end{equation}
They may be put in symmetrizable (but no longer conservative) 
form by a further change, using identity
\begin{equation}
\curl (H \times u) = 
 (\div u) H+ (u\cdot \na)H -(\div H) u- (H\cdot \na)u
\end{equation}
together with constraint \eqref{divfree} to express the third equation as
\begin{equation}\label{symm3}
 \D_t H + 
 (\div u) H+ (u\cdot \na)H - (H\cdot \na)u= \mu \eps \Delta H.
\end{equation}
Forgetting  the constraint equation, we get a $7 \times 7$ symmetric system.

Neglecting zero-th order  terms,   
  the linearized equations of \eqref{mhdeq} about $(\rho, u, H)$ are
\begin{equation}\label{linearized}
\left\{ \begin{aligned}
 & D_t  \dot \rho + \rho \div \dot u    
 \\
 & \rho  D_t \dot u +     c^2 \na \dot \rho  +   H \times \curl \dot H  - \eps \nu \Delta \dot u 
 \\
 &  
 D_t \dot H + (\div \dot u) H   - H \cdot \na  \dot u  - \eps  \mu \Delta \dot H 
 \end{aligned}\right.
\end{equation}
with $D_t = \D_t + u \cdot \na$ and $c^2  = dp/d\rho$ which we assume to be positive. 
This system is hyperbolic symmetric, 
with symmetrizer 
$S = \bdiag( c^2, \rho \Id, \Id ) $.  
It enters the general framework  of linearized equations studied in this paper, with parameters
 $  (\rho, u, H)$.

 
 \subsection{Eigenvalues and eigenvectors}

Eigenvalues and eigenvectors    of the symbol solve  
\begin{equation}
\label{Asymbol}
\left\{ \begin{aligned}
 & \tilde \tau   \dot \sigma  + \rho (\xi \cdot \dot u)   = 0 
 \\
 &    \tilde \tau  \dot u +   c^2  \dot \sigma \xi    +   v \times (\xi \times \dot  v  )  = i 
  \nu | \xi^2 | \dot u / \rho 
 \\
 &  
\tilde  \tau \dot v +  ( \xi \cdot \dot u) v  - (v \cdot \xi)  \dot u   = i \mu | \xi |^2 \dot v. 
 \end{aligned}\right.
\end{equation}
with 
\begin{equation}
\label{nwvar}
\tilde \tau  = \tau + u \cdot \xi, \quad 
v = H \sqrt \rho, \quad \dot \sigma = \dot \rho / \rho , \quad 
 \dot v = \dot H / \sqrt \rho . 
 \end{equation}  
 The structure condition \eqref{struc1} is satisfied with $N'= 6$. The kernel of 
 $B(\xi)$ is generated by ${}^\trans ( 1, 0, \ldots 0)$ which is never an eigenvector of 
 $A(\xi)$ when $\xi \ne 0$. Thus the Assumptions (H1'), (H1) and (H2) are satisfied. 

\medbreak

Consider next the inviscisd problem. 
  The seven  eigenvalues of 
  $A( \rho, u,H ,  \xi) $   are (see e.g. \cite{MZ2}): 
  \begin{equation}
\label{7.21}
\left\{\begin{aligned} 
  \lambda_0  & =  u\cdot \xi , 
\\
 \lambda_{\pm s} & =  \lambda_0 \pm c_s  | \xi |, 
\\
 \lambda_{\pm 2} & = \lambda_0  \pm  v  \cdot \xi , 
\\
  \lambda_{\pm f} & = \lambda_0 \pm c_f  | \xi | , 
\end{aligned}\right.
\end{equation}
  with 
$$ 
  \begin{aligned}
  &  c_f^2  := \frac{1}{2} 
  \Big( c^2 + |v |^2   + \sqrt{ (c^2 - | v |^2)^2 + 4 b^2 c^2  } \Big)
  \\
 &   c_s^2  := \frac{1}{2} 
  \Big( c^2 + | v| ^2   -  \sqrt{ (c^2 - | v | ^2)^2 +  4 b^2  c^2   } \Big), 
  \\
  &c^2 = p'(\rho) > 0 , \quad  v = H/ \sqrt{\rho},  \quad  b= | \hat \xi \times v |,  \quad 
  \hat \xi = \xi/ | \xi |.  
\end{aligned}
$$

  The first eigenvalue corresponds to the transport of the 
  constraint. It can be decoupled from the system :  there is a smooth one dimensional subspace, 
  $\EE_0$ such that $A( \xi) = \lambda_0$  on this space 
  and $\EE_0^\perp$ is stable for $A(  \xi)$. 
  The other eigenvalues are in general simple.  
    
 \begin{lem} [\cite{MZ2}] 
 \label{lem71}
  Assume that $ 0 < \vert v \vert^2  \ne  c^2$. Consider 
  $\xi  \in \RR^3\backslash\{0\}$. 
  
 \quad i)  When $ \xi  \cdot v \ne 0$ and 
 $    \xi  \times H \ne 0$, the eigenvalues are simple.

 \quad ii) On the manifold $  \xi  \times  v =  0$, $\lambda_0$ is simple. 
 When $\vert v \vert^2 <   c^2$  [resp. $\vert v \vert^2 >  c^2$ ], 
 $\lambda_{\pm f}$ [resp. $\lambda_{\pm s}$] are simple, the other eigenvalues  
 $\lambda_{\pm 2} = \lambda_{\pm s}$  [resp. 
  $\lambda_{\pm 2} = \lambda_{\pm f}$]  are double, algebraically regular but not geometrically
  regular. Moreover,   
  \begin{equation}
  \label{DL2}
    \lambda_{\pm 2} -    \lambda_{\pm s} = O\big( | \xi  \times  v|^2) \quad  [resp. \ 
  \lambda_{\pm 2} - \lambda_{\pm f} = O\big( | \xi  \times  v|^2)] . 
  \end{equation}

   \quad iii) On the manifold $ \xi  \cdot v =  0$  the eigenvalues $\lambda_{\pm f}$ 
 are simple and the multiple eigenvalue $\lambda_0 = \lambda_{\pm s} = \lambda_{\pm 2}$
 is geometrically regular.  More precisely, there are smooth $\lambda_{\pm 1}$
 such that 
  $\{ \lambda_{  s}, \lambda_{- s} \} = \{ \lambda_{ 1}, \lambda_{-1} \}$. 
 Moreover,  
 \begin{equation}
\label{DL1}
\lambda_{\pm1}   = u \cdot \xi  \pm  \delta v \cdot \xi +  O\big((   v  \cdot \xi )^2 \big), 
\quad \delta = \frac {c}{\sqrt{c^2 + h^2}} . 
\end{equation}
 One can choose smooth eigenvectors $e_0, e_{\pm 1}, e_{\pm2}$  such that,
on the manifold  $ \xi \cdot v = 0$,  
\begin{equation}
\label{VP1}
e_0 = \begin{pmatrix}
     0  \\ 0       \\  \hat \xi   
\end{pmatrix}, \quad
e_{\pm 1}  =  \frac{\delta}{\sqrt 2 | v| }  \begin{pmatrix}
  -      | v |^2 / c^2   \\ \mp v / \delta    \\   v 
\end{pmatrix}, 
\quad 
e_{\pm 2}  = 
 \frac{1 }{\sqrt 2 | v| } \begin{pmatrix}
    0   \\  \mp   w \\ w 
\end{pmatrix}, 
\end{equation}
  with $ w = \hat \xi  \times v $. 

   \end{lem}


\subsection{Glancing and  viscous coupling}

   The boundary $\{   x_3  = 0 \} $ is noncharacteristic for the hyperbolic part if and only if 
  \begin{equation}
\label{7.23}
u_3   \notin \big\{ 0, \  v_3  , \ \pm c_s(n), \ \pm c_f(n)
\big\}  
\end{equation}
  where $c_s(n)$ and $c_f(s)$ are the slow and fast speed computed in the 
  normal direction $n = (0, 0, 1)$.

  \begin{lem} 
  \label{lem72}
   Assume that $ 0 < \vert v \vert   \ne  c $.
  
  i)  On the manifold $\xi \times v = 0$,  the multiple eigenvalues are 
  nonglancing  if and only if $ u_3 \ne \pm v_3$. In this case, they are totally nonglancing.
  
  ii) On the manifold $\xi \cdot v = 0$, the multiple eigenvalues are non glancing 
  if and only if $ u_3 \ne 0$, $u_3 \ne \pm v_3 $  and  $u_3 \ne \pm \delta v_3$.
  They are totally nonglancing when $ | u_3 | > | v_3 |$. 
  
  \end{lem} 
  
  \begin{proof}
By  \eqref{DL2}, on $\xi \times v = 0$,  with $j = s$ when 
$| v | < c$ and $j = f $ when $| v | > c$,  there holds
\begin{equation*}
\D_{\xi_3}  \lambda_{\pm j  } = \D_{\xi_3}  \lambda_{\pm 2  } = u_3 \pm v_3 . 
\end{equation*}
 This implies $i)$. 
  
  In addition,  $\D_{\xi_3}  \lambda_{0 } = u_3 $, 
  $ \D_{\xi_3}  \lambda_{\pm 2  } = u_3 \pm v_3 $, and by \eqref{DL1}
  $ \D_{\xi_3}  \lambda_{\pm 1  } = u_3 \pm \delta v_3 $ on the manifold 
  $\xi \dot v = 0$. This implies $ii)$.   
 \end{proof} 
  
  Next we study the viscous coupling of  vectors  $e_j$ at geometrically regular 
  modes. 
  In the variables $(\dot \rho/ \rho, \dot u, \dot v)$, the system
  \eqref{linearized} is symmetric, with symmetrizer $S = \mathrm{diag} (c^2, \Id, \Id)$,  
  and the viscosity matrix is  $B (\xi) = | \xi |^2 \mathrm{diag} (0,  \nu \Id/ \rho, \mu \Id)$.  
  The basis \eqref{VP1} is orthonormal  for $S$. Therefore, 
  according to the general rule \eqref{bssym},   the 
  matrix 
  $B^\sharp$ is symmetric with nondiagonal entries  
 \begin{equation}
\label{Bsex}
\begin{aligned}
& B^\sharp_{0, \pm 1} = B^\sharp_{0, \pm 2} =   B^\sharp_{\pm1 , \pm 2}= 0 , 
\\
& B^\sharp_{1, -  1} =     \frac{\delta^2 \mu}{2}  - \frac{\nu}{2 \rho},  \qquad 
B^\sharp_{2, -2} =  \frac{  \mu}{2}  - \frac{\nu}{2 \rho}. 
\end{aligned}
\end{equation}
  When $ | u_3| <  | v_3 | $, then   one of the eigenvalue $\lambda_{\pm 2}$ is incoming 
  and the other one outgoing (depending on the sign of $v_3$).  Therefore, 
  if $\mu - \nu / \rho$,  the coupling coefficient $B^\sharp_{2, -2} $ does not vanish. 
   Summing up, we have proved: 
  
  \begin{lem}
  If $ | u_3| <  | v_3 | $,  and $\nu \ne \rho \mu$, then the decoupling condition 
  $\eqref{331x}$ is not satisfied at modes where $\xi \cdot v = 0$. 
  \end{lem} 
  
  \begin{rem}
  \textup{The decoupling of the mode $\lambda_0$ from the other ones  reflects 
  that the constraint \eqref{divfree} is propagated by the viscous equation as well. 
  The other partial decoupling observed above depend on the particular choice
  of the viscosity matrices and disappear for  general $B$.  }
  \end{rem}


\subsection{Shocks}

Consider an inviscid planar shock. We suppose that the front is 
$x_3 = \sigma t$ and denote by $(\rho^-, u^-, H^-)$ and 
$(\rho^+ , u^+, H^+)$ the states on the left and on the right respectively. 
All the analysis of the preceding section is valid, if we change 
$u_3$ to $u_3 - \sigma$.  

The jump conditions  are deduced from  the conservative form 
of the equations: 
\begin{equation}\label{RH}
\left\{ \begin{aligned}
 & [\rho (u_3- \sigma) ] = 0, 
 \\
 &[\rho u (u_3 - \sigma) ] +  r_3   \left[  p + \frac{1}{2} |H|^2\right] - [ H_3 H ]=0, 
 \\
 &  [(u_3 - \sigma ) H ] - [ H_3 u ] = 0, 
 \\
 & [ H_3 ] = 0,    
\end{aligned}\right. 
\end{equation}
where $r_3 = {}^{t} (0, 0, 1)$.  The last jump condition comes from the constraint equation \eqref{divfree}. Apparently
this system of 8 scalar equations is too large. However, projecting the third equation 
in the normal direction yields $\sigma [H_3 ] = 0$ which is implied by the last equation. 
This shows that \eqref{RH} is made of 7 independent equations, as expected. 
Denoting by $u_{\rm tg}$  and $H_{\rm tg}$ the tangential part of $u$ and $H$, 
that is their orthogonal projection on $r_3^\perp$, \eqref{RH} is equivalent to 
\begin{equation}\label{RHm}
\left\{ \begin{aligned}
 & [\rho (u_3- \sigma) ] = 0, 
 \\
 &[\rho u (u_3 - \sigma) ] +  r_3  \left[  p + \frac{1}{2} |H|^2\right] - [ H_3 H ]=0, 
 \\
 &  [(u_3 - \sigma ) H_{\rm tg}  ] - [ H_3 u_{\rm tg}  ] = 0,  \quad  [ H_3 ] = 0.    
\end{aligned}\right. 
\end{equation}

\bigbreak
\noindent 1. \emph{Fast Lax' shocks.}   Consider an extreme shock. 
Changing $x$ to $-x$ if necessary, the Lax condition read: 
\begin{equation}
\label{fLax}
\begin{aligned}
        & u^-_3 +   | v^-_3 | <     \sigma  <  u^-_3 + c^-_f,  \\
       & u^+_3 +    c^+_f  <     \sigma   . 
\end{aligned}
\end{equation}
In particular, this implies that the front is not characteristic on both side, and that 
the nonglancing conditions in Lemma~\ref{lem72} are also satisfied on both side
and the multiple modes are totally nonglancing. 
Therefore : 

\begin{prop}
For extreme Lax shocks, the assumptions of Theorem~$\ref{th11} $ are satisfied. 
\end{prop} 

\bigbreak
\noindent 2. \emph{Slow Lax' shocks.}  
 Consider a shock associated to one of the middle eigenvalue $\lambda_{\pm s}$. 
Changing $x$ to $-x$ if necessary, the Lax condition read: 
\begin{equation}
\label{sLax}
\begin{aligned}
        & u^-_3 -     c^-_s  <     \sigma  <  u^-_3 + c^-_s   \\
       & u^+_3 +    c^+_s  <     \sigma    < u^+_3 + | v^+_3 |.  
\end{aligned}
\end{equation}
On both side we have 
$| u_3 - \sigma | < | v_3 |$, therefore

\begin{prop}
For slow  Lax shocks, the decoupling condition is never  satisfied. 
\end{prop}

 
\subsection{The $H \to 0$ limit}

When $H = 0$, the system \eqref{mhdeq} reduces to isentropic Euler's equations
and \eqref{RHm} to the corresponding Rankine Hugoniot condition.

When $H = 0$, the eigenvalues are 
\begin{equation}
\label{7.24}
\lambda_0 = \lambda_{\pm1} = \lambda_{\pm 2} =  u \cdot \xi , \quad 
 \lambda_{\pm 3} = \lambda_0 \pm c \vert   \xi \vert. 
\end{equation} 
In particular $\D_{\xi_3} \lambda_0 = u_3$.  
Moreover,  at $ \underline U = (\rho, u, 0)$, the tangent characteristic polynomial
$\underline \Delta$ in \eqref{Deltatan} is $ (\tau + u \cdot \xi )^5$. Therefore, 
 if $u_3 \ne 0$, the eigenvalue $\lambda_0$ is totally nonglancing.

  \begin{lem}
  \label{lem74}
  Consider a state $\underline U =  (\underline \rho, \underline u, 0)$. 
 Suppose that 
   \begin{equation}
  \label{7.25}
\underline u_3   \notin \{  - \underline c,  0, + \underline c \}. 
  \end{equation}
  Then, for $U$ in a neighborhood of $\underline U$,  the boundary $x_3 = 0$ is non characteristic for the hyperbolic linearized equation and the eigenvalues $\lambda_{\pm 3}$ are simple.
  Moreover, for all $\xi \ne 0$, 
the multiple eigenvalue $\lambda_0$ is totally 
  nonglancing at $\underline U$.  
    \end{lem}


 


\end{document}